\documentclass[letterpaper]{amsart}
\usepackage{subfiles}

\usepackage[utf8]{inputenc}
\usepackage{amsmath}
\usepackage{amsfonts}
\usepackage{amsthm}
\usepackage{mathrsfs}
\usepackage{enumitem}
\usepackage{dsfont}
\usepackage{amssymb}
\usepackage{pifont}
\usepackage{mathtools}
\usepackage{graphicx}
\usepackage{color}
\usepackage{import}
\usepackage{tikz-cd}
\usepackage{hyperref}
\usepackage{import}

\newtheorem{thm}{Theorem}[section]
\newtheorem{lem}[thm]{Lemma}
\newtheorem{prop}[thm]{Proposition}

\newtheorem{cor}[thm]{Corollary}

\newtheorem*{lem*}{Lemma}
\newtheorem*{prop*}{Proposition}

\theoremstyle{definition}

\newtheorem{defn}[thm]{Definition}

\newtheorem{remark}[thm]{Remark}

\numberwithin{equation}{section}

\newenvironment{aside}[1]{\par\smallskip\noindent \textit{#1}.}{\par}

\newcommand{\eps}{\varepsilon}
\newcommand{\dee}{\ensuremath{\partial}}

\newcommand{\mbf}[1]{\ensuremath{\mathbf{#1}}}

\newcommand{\mc}[1]{\ensuremath{\mathcal{#1}}}

\newcommand{\mr}[1]{\ensuremath{\textrm{#1}}}

\newcommand{\beq}{\begin{equation}}
\newcommand{\eeq}{\end{equation}}

\newcommand{\pair}[2]{\ensuremath{{\langle #1 , #2 \rangle}}}

\newcommand{\N}{\ensuremath{\mathbb{N}}}
\newcommand{\Z}{\ensuremath{\mathbb{Z}}}

\newcommand{\R}{\ensuremath{\mathbb{R}}}

\renewcommand{\H}{\ensuremath{\mathbb{H}}}
\renewcommand{\P}{\ensuremath{\mathbb{P}}}

\newcommand{\E}{\ensuremath{\mathbb{E}}}

\newcommand{\Rd}{\ensuremath{\R^d}}
\newcommand{\RP}{\ensuremath{\R\mathrm{P}}}

\DeclareMathOperator{\Aut}{Aut}

\DeclareMathOperator{\End}{End}

\DeclareMathOperator{\Hom}{Hom}

\DeclareMathOperator{\Gr}{Gr}

\DeclareMathOperator{\supp}{supp}
\DeclareMathOperator{\Stab}{Stab}

\DeclareMathOperator{\SL}{SL}
\DeclareMathOperator{\GL}{GL}
\DeclareMathOperator{\PGL}{PGL}
\DeclareMathOperator{\SO}{SO}
\DeclareMathOperator{\PO}{PO}

\newcommand{\Heq}[1]{\ensuremath{[#1]_{\mc{H}}}}

\newcommand{\strat}{\ensuremath{F_\Omega}}

\newcommand{\Cor}{\ensuremath{\mathrm{Cor}_\Omega}}
\newcommand{\limset}{\Lambda_\Omega(\Gamma)}
\newcommand{\CH}{\ensuremath{\mathrm{Hull}_\Omega}}

\title[Dynamical properties of convex cocompact actions]{Dynamical
  properties of convex cocompact actions in projective space}

\author{Theodore Weisman} \address{Department of Mathematics,
  University of Michigan, Ann Arbor, MI 48109, USA}
\email{tjwei@umich.edu} \date{\today}

\begin{document}

\begin{abstract}
  We give a dynamical characterization of convex cocompact group
  actions on properly convex domains in projective space in the sense
  of Danciger-Guéritaud-Kassel: we show that convex cocompactness in
  $\RP^d$ is equivalent to an expansion property of the group about
  its limit set, occuring in different Grassmannians. As an
  application, we give a sufficient and necessary condition for convex
  cocompactness for groups which are hyperbolic relative to a
  collection of convex cocompact subgroups. We show that convex
  cocompactness in this situation is equivalent to the existence of an
  equivariant homeomorphism from the Bowditch boundary to the quotient
  of the limit set of the group by the limit sets of its peripheral
  subgroups.
\end{abstract}

\maketitle

\tableofcontents


\section{Introduction}
\label{sec:introduction}

When $G$ is a rank-one semisimple Lie group, a \emph{convex cocompact}
subgroup is a discrete group $\Gamma \subset G$ which acts cocompactly
on some convex set in the Riemannian symmetric space $G / K$, where
$K$ is a maximal compact in $G$. Convex cocompact subgroups in
rank-one have long been objects of great interest, and have a wide
variety of possible characterizations.

More recently, efforts have been underway to understand the
appropriate generalization of convex cocompactness in higher-rank Lie
groups. A key concept is \emph{Anosov representations}: discrete and
faithful representations of word-hyperbolic groups into reductive Lie
groups which generalize a \emph{dynamical} definition of convex
cocompact subgroups in rank one. Anosov representations were first
defined for surface groups by Labourie \cite{labourie2006anosov}, and
the definition was later extended to general word-hyperbolic groups by
Guichard-Wienhard \cite{gw2012anosov}. Guichard-Wienhard also
demonstrated that an Anosov representation $\rho$ of a group $\Gamma$
can be interpreted as the holonomy of a certain geometric structure
associated to $\rho$. Currently, understanding the connection between
Anosov representations and geometric structures is an area of active
research.

In \cite{dgk2017convex}, Danciger, Guéritaud, and Kassel developed a
notion of convex cocompact representations in $\PGL(d,\R)$ that (as in
the rank-one setting) have concrete and transparent \emph{convex}
geometric objects associated to them---in this case, a compact
manifold (or orbifold) with a \emph{convex projective structure}.

Recall that a subset $\Omega$ of projective space $\RP^{d-1}$ is
\emph{convex} if it is contained in some affine chart in $\RP^{d-1}$,
and $\Omega$ is a convex subset of that affine chart. The set $\Omega$
is \emph{properly convex} if its closure is also contained in an
affine chart, and it is a \emph{properly convex domain} if it is also
open. A \emph{convex projective orbifold} is a quotient of a convex
set in $\RP^{d-1}$ by a discrete subgroup of $\Aut(\Omega)$, where
\[
  \Aut(\Omega) := \{g \in \PGL(d, \R) : g \cdot \Omega = \Omega\}.
\]

The Danciger-Guéritaud-Kassel definition of convex cocompactness in
$\RP^{d-1}$ says that a group $\Gamma \subset \PGL(d,\R)$ is convex
cocompact when it is the holonomy of a compact convex projective
orbifold satisfying certain conditions.

\begin{defn}
  \label{defn:full_orbital_limit_set}
  Let $\Omega$ be properly convex domain in $\RP^{d-1}$, and let
  $\Gamma \subseteq \Aut(\Omega)$.
  \begin{itemize}
  \item The \emph{full orbital limit set} $\Lambda_\Omega(\Gamma)$ is
    the set of accumulation points in $\dee \Omega$ of $\Gamma$-orbits
    in $\Omega$, i.e. the union
    \[
      \bigcup_{x \in \Omega}(\overline{\Gamma \cdot x} \cap \dee
      \Omega).
    \]
  \item The \emph{convex core} of $\Gamma$ in $\Omega$, denoted
    $\Cor(\Gamma)$, is the convex hull in $\Omega$ of the full orbital
    limit set $\Lambda_\Omega(\Gamma)$.
  \end{itemize}
\end{defn}

\begin{defn}[{\cite[Definition 1.11]{dgk2017convex}}]
  \label{defn:dgk_convex_cocompact}
  Let $\Omega$ be a properly convex domain in $\RP^{d-1}$, and let
  $\Gamma$ be a discrete group acting by projective transformations on
  $\Omega$. The group $\Gamma$ \emph{acts convex cocompactly on
    $\Omega$} if it acts cocompactly on $\Cor(\Gamma)$.

  A group $\Gamma \subset \PGL(d, \R)$ \emph{acts convex cocompactly
    in $\RP^{d-1}$} if it acts convex cocompactly on \emph{some}
  properly convex domain $\Omega \subset \RP^{d-1}$.
\end{defn}
Note that this definition is strictly stronger than merely asking for
$\Gamma$ to act cocompactly on some $\Gamma$-invariant convex subset
of a properly convex domain $\Omega$.

Danciger-Guéritaud-Kassel prove that when a discrete word-hyperbolic
group $\Gamma \subset \PGL(d, \R)$ preserves a properly convex domain,
$\Gamma$ acts convex cocompactly on some domain
$\Omega \subset \RP^{d-1}$ precisely when the inclusion
$\Gamma \hookrightarrow \PGL(d, \R)$ is \emph{$P_1$-Anosov}; a related
result was independently shown by Zimmer in
\cite{zimmer2017projective}. Moreover, a group acting convex
cocompactly on a domain $\Omega$ is word-hyperbolic precisely when
there are no nontrivial projective segments in its full orbital limit
set (and in this case the definition is equivalent to a similar notion
of convex cocompactness introduced by Crampon-Marquis in
\cite{cm2014finitude}).

However, there are also many examples of non-hyperbolic groups which
have convex cocompact representations as in Definition
\ref{defn:dgk_convex_cocompact}; see Section
\ref{subsec:examples}. These non-hyperbolic convex cocompact groups
are somewhat more mysterious than their hyperbolic
counterparts. Recently, however, there has been significant progress
towards a deeper understanding of them, especially in the case where
the $\Gamma$-action on the entire domain $\Omega$ is cocompact: see
e.g.  \cite{islam2019rank}, \cite{bobb2020codimension},
\cite{zimmer2020higher}. Of particular relevance to this paper is the
description, due to Islam-Zimmer \cite{iz2019flat},
\cite{iz2019convex}, of the domains with a convex cocompact action by
a relatively hyperbolic group relative to a family of virtually
abelian subgroups of rank at least two.

\subsection{Convex cocompactness and Anosov dynamics}

The first main result of this paper (Theorem
\ref{thm:convex_cocompact_equals_expansion} below) is a dynamical
characterization of convex cocompactness that applies even for
non-hyperbolic groups, generalizing the relationship between convex
cocompactness and $P_1$-Anosov representations. This addresses a
question asked by Danciger-Guéritaud-Kassel in \cite{dgk2017convex}.

The main idea behind Theorem
\ref{thm:convex_cocompact_equals_expansion} is to generalize the
dynamical description of convex cocompactness explored by Sullivan
\cite{sullivan1985quasiconformal} in the rank-one setting: when
$\Gamma$ is a discrete subgroup of $\PO(d,1)$, $\Gamma$ is convex
cocompact if and only if $\Gamma$ satisfies an expansion/contraction
property about its limit set in $\dee \H^d$.

More generally, when $\rho:\Gamma \to G$ is a $P$-Anosov
representation (for $P$ a parabolic subgroup of a reductive Lie group
$G$), work of Kapovich-Leeb-Porti \cite{klp2017} shows that
$\Gamma$ satisfies an expansion property on the flag manifold
$G/P$. Kapovich-Leeb-Porti also showed in \cite{klpdomains} (using the
same basic idea as Sullivan) that this expansion property can be used
to find cocompact domains of discontinuity for Anosov representations
in $G/P$.

When $\Gamma$ is hyperbolic with $\rho:\Gamma \to \PGL(d, \R)$ convex
cocompact, $\rho$ is $P_1$-Anosov, yielding an expansion property in
$\RP^{d-1}$.  Theorem \ref{thm:convex_cocompact_equals_expansion} says
that convex cocompact representations $\Gamma \to \PGL(d, \R)$ are
characterized by a similar expansion property on \emph{multiple} flag
manifolds (and some additional technical conditions): different
elements of $\Gamma$ expand neighborhoods in different Grassmannians
$\Gr(k, d)$. That is, convex cocompactness in $\PGL(d,\R)$ is
equivalent to a kind of ``mixed Anosov'' property.

Using work of Cooper-Long-Tillmann \cite{clt2018deforming},
Danciger-Gu\'eritaud-Kassel also observe that convex cocompactness in
the sense of Definition \ref{defn:dgk_convex_cocompact} is
\emph{stable}: if $\Gamma \subset \PGL(d, \R)$ is a convex cocompact
subgroup, then there is an open subset $\mc{U}$ of
$\Hom(\Gamma, \PGL(d,\R))$, containing the inclusion
$\Gamma \hookrightarrow \PGL(d ,\R)$, such that any $\rho \in \mc{U}$
is injective and discrete with $\rho(\Gamma)$ convex cocompact.

Theorem \ref{thm:convex_cocompact_equals_expansion} then implies that
the ``mixed Anosov'' property we consider in this paper is also stable
under small deformations. This suggests one possible route towards a
generalization of Anosov representations for non-hyperbolic groups.

\subsubsection{Structure of the boundary of a convex domain}

The expansion property we use to characterize convex cocompactness in
$\PGL(d, \R)$ is given in terms of the natural decomposition of the
boundary of a convex domain into convex pieces. Each point $x$ in the
boundary of a convex set $\Omega$ lies in a unique \emph{face}
$\strat(x)$: the set of all points in $\dee \Omega$ which lie in a
common open line segment in $\dee \Omega$ with $x$. The \emph{support}
$\supp(F)$ of a face $F$ in $\dee \Omega$ is the projective span of
$F$. We can view the support as an element of the Grassmannian of
$k$-planes $\Gr(k,d)$, for some $1 \le k < d$.

When $\Lambda$ is a subset of $\dee \Omega$, we can ask for it to be
well-behaved with respect to the decomposition of $\dee \Omega$ into
faces.
\begin{defn}
  \label{defn:fully_stratified}
  Let $\Lambda$ be a subset of $\dee \Omega$. If, for all $x \in
  \Lambda$, we have
  \[
    \strat(x) \subset \Lambda,
  \]
  we say that $\Lambda$ \emph{contains all of its faces}.
\end{defn}

\begin{defn}
  \label{defn:boundary_convex}
  Let $\Lambda$ be a subset of $\dee \Omega$. We say that $\Lambda$ is
  \emph{boundary-convex} if any supporting hyperplane of $\Omega$
  intersects $\Lambda$ in a convex set.
\end{defn}

\begin{figure}[h]
  \centering

  \def\svgwidth{1.5in}
\begingroup%
  \makeatletter%
  \providecommand\color[2][]{%
    \errmessage{(Inkscape) Color is used for the text in Inkscape, but the package 'color.sty' is not loaded}%
    \renewcommand\color[2][]{}%
  }%
  \providecommand\transparent[1]{%
    \errmessage{(Inkscape) Transparency is used (non-zero) for the text in Inkscape, but the package 'transparent.sty' is not loaded}%
    \renewcommand\transparent[1]{}%
  }%
  \providecommand\rotatebox[2]{#2}%
  \newcommand*\fsize{\dimexpr\f@size pt\relax}%
  \newcommand*\lineheight[1]{\fontsize{\fsize}{#1\fsize}\selectfont}%
  \ifx\svgwidth\undefined%
    \setlength{\unitlength}{195.18225383bp}%
    \ifx\svgscale\undefined%
      \relax%
    \else%
      \setlength{\unitlength}{\unitlength * \real{\svgscale}}%
    \fi%
  \else%
    \setlength{\unitlength}{\svgwidth}%
  \fi%
  \global\let\svgwidth\undefined%
  \global\let\svgscale\undefined%
  \makeatother%
  \begin{picture}(1,1.52111257)%
    \lineheight{1}%
    \setlength\tabcolsep{0pt}%
    \put(0,0){\includegraphics[width=\unitlength,page=1]{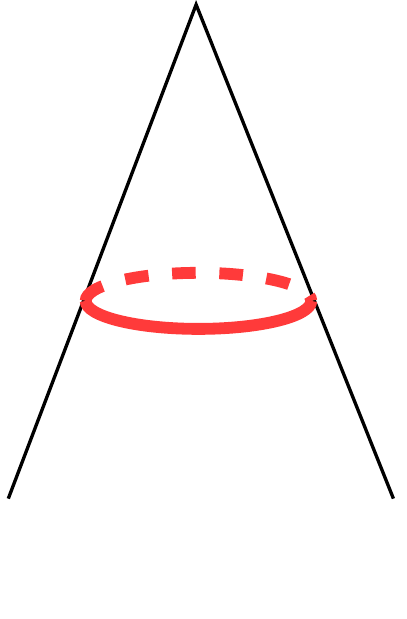}}%
    \put(0.44291935,0.61118013){\color[rgb]{0,0,0}\makebox(0,0)[lt]{\lineheight{1.25}\smash{\begin{tabular}[t]{l}$\Lambda_1$\end{tabular}}}}%
    \put(0,0){\includegraphics[width=\unitlength,page=2]{convexity_conditions_2.pdf}}%
    \put(0.41522079,0.02061817){\color[rgb]{0,0,0}\makebox(0,0)[lt]{\lineheight{1.25}\smash{\begin{tabular}[t]{l}$\Lambda_2$\end{tabular}}}}%
  \end{picture}%
\endgroup%

  \caption{Here $\Omega$ is a cone, the convex hull of a disk and a
    point in $\RP^3$. The set $\Lambda_1$ (red) does not contain all
    of its faces (Definition \ref{defn:fully_stratified}), because
    each point in $\Lambda_1$ is contained in the interior of a line
    segment in $\dee \Omega$ which is not contained in
    $\Lambda_1$. The set $\Lambda_2$ (blue) contains all of its faces,
    but it is not boundary-convex (Definition
    \ref{defn:boundary_convex}): a line segment joining two points of
    $\Lambda_2$ intersects $\dee \Omega - \Lambda_2$.}
\end{figure}

If $\Lambda = \limset$ for a group $\Gamma \subseteq \Aut(\Omega)$
acting convex cocompactly on some domain $\Omega$, then \cite[Lemma
4.1]{dgk2017convex} implies that $\Lambda$ is closed and boundary
convex, and contains all of its faces.

When $\Lambda$ is a subset of $\dee \Omega$ containing all of its
faces, and $\Gamma \subset \PGL(d,\R)$ is a group preserving
$\Lambda$, we say $\Gamma$ is \emph{expanding at the faces of
  $\Lambda$} (Definition \ref{defn:expansion_in_supports}) if for
every face $F$ in $\Lambda$, the group $\Gamma$ has an element which
is expanding in a neighborhood of $\supp(F)$ in $\Gr(k,d)$. When the
expansion constants can be chosen uniformly, we say $\Gamma$ is
\emph{uniformly} expanding at the faces of $\Lambda$. For the full
definitions, see Section \ref{subsec:support_dynamics}.

Here is the precise version of our characterization of convex
cocompactness:
\begin{thm}
  \label{thm:convex_cocompact_equals_expansion}
  Let $\Omega$ be a properly convex domain in $\RP^{d-1}$, and let
  $\Gamma$ be a discrete subgroup of $\Aut(\Omega)$. The following are
  equivalent:
  \begin{enumerate}
  \item $\Gamma$ acts convex cocompactly on $\Omega$.\label{item:main_thm_1}
  \item There is a closed, $\Gamma$-invariant, and boundary-convex
    subset $\Lambda \subset \dee \Omega$ with nonempty convex hull
    such that $\Lambda$ contains all of its faces and $\Gamma$ is
    uniformly expanding at the faces of $\Lambda$.\label{item:main_thm_2}
  \end{enumerate}
  In this case, the set $\Lambda$ is the full orbital limit set
  $\limset$.
\end{thm}

\begin{remark}
  When we prove the implication (\ref{item:main_thm_2}) $\implies$
  (\ref{item:main_thm_1}) of Theorem
  \ref{thm:convex_cocompact_equals_expansion}, we will not actually
  need to assume that the expansion at the faces of $\Lambda$ is
  uniform---only that the expansion occurs with respect to a
  particular choice of Riemannian metric on $\Gr(k,d)$. See Remark
  \ref{remark:support_expanding_well_defined}.
\end{remark}

A special case of convex cocompactness is when a discrete group
$\Gamma \subset \Aut(\Omega)$ acts cocompactly on all of $\Omega$. In
this case, we say that $\Omega$ is \emph{divisible}, and the group
$\Gamma$ \emph{divides} the domain. As $\dee \Omega$ is always
boundary convex and contains all of its faces, when
$\Lambda = \dee \Omega$, Theorem
\ref{thm:convex_cocompact_equals_expansion} can be stated as the
following:
\begin{cor}
  \label{cor:expanding_implies_divides}
  Let $\Gamma$ be a discrete subgroup of $\PGL(d,\R)$ preserving a
  properly convex domain $\Omega$. Then $\Gamma$ divides $\Omega$ if
  and only if $\Gamma$ is uniformly expanding at the faces of
  $\dee \Omega$.
\end{cor}

\subsubsection{Spaces of properly convex domains}
\label{sec:relative_benzecri_comments}

A key technical tool we need for the implication \ref{item:main_thm_2}
$\implies$ \ref{item:main_thm_1} in
Theorem~\ref{thm:convex_cocompact_equals_expansion} is a version of
the \emph{Benz\'ecri cocompactness theorem} for properly convex
domains in $\RP^{d-1}$, which applies relative to a direct sum
decomposition of $\R^d$. This result (proved in
Section~\ref{sec:spaces_of_domains}) may be of independent interest,
so we describe it briefly here.

In \cite{benzecri1960varietes}, Benz\'ecri showed that the group
$\PGL(d, \R)$ acts properly and cocompactly on the space of pointed
properly convex domains in $\RP^{d-1}$. One immediate and useful
consequence of this fact is that if $\Omega \subset \RP^{d-1}$ is a
properly convex domain, and $x_n$ is any sequence of points in
$\Omega$, then after extracting a subsequence, there is a sequence of
projective transformations $g_n \in \PGL(d, \R)$ such that $g_n\Omega$
converges to a fixed properly convex domain $\Omega_\infty$, and
$g_nx_n$ converges to a point in the interior of $\Omega_\infty$.

We prove a result (Proposition~\ref{prop:relative_benzecri}) which
provides some control over the group elements $g_n$ appearing
above. Roughly, our result says that if the sequence $x_n$ lies in a
lower-dimensional ``slice'' $W \subset \Omega$ satisfying some
technical conditions, then the sequence $g_n$ above can be chosen to
preserve \emph{both} the projective subspace spanned by $W$ and a
fixed complementary subspace in $\R^d$.

Proposition~\ref{prop:relative_benzecri} can be compared to earlier
work of Benoist \cite{benoist2003convexes} and Frankel
\cite{frankel89}. These results effectively show that when $x_n$ lies
in a lower-dimensional ``slice'' in $\Omega$ as above, then the
sequence $g_n$ can always be chosen to preserve the span of the
slice. The proposition we prove in this paper has stronger hypotheses,
but also a stronger conclusion. The more precise control we get is
necessary for the intended application---see
Remark~\ref{rem:pseudo_loxodromic_use}.

\subsection{Relative hyperbolicity}

In the second part of this paper, we use the dynamical
characterization of convex cocompactness given by Theorem
\ref{thm:convex_cocompact_equals_expansion} to study convex
cocompactness for a group $\Gamma$ which is hyperbolic relative to a
collection $\mc{H}$ of convex cocompact subgroups. We will give a
necessary and sufficient condition for such a group to act convex
cocompactly in terms of an embedding of the Bowditch boundary
$\dee(\Gamma, \mc{H})$. This strengthens the connection between convex
cocompact groups in $\PGL(d, \R)$ and Anosov representations, since
Anosov representations can also be characterized by the existence of
an equivariant embedding of the \emph{Gromov} boundary of a hyperbolic
group into some flag manifold (see \cite{gw2012anosov},
\cite{klp2017}).

We will give a definition of relatively hyperbolic groups in terms of
convergence dynamics in Section \ref{sec:relative_hyp_background}.

\begin{defn}
  Let $\mc{H} = \{H_i\}$ be a collection of subgroups of
  $\PGL(d, \R)$, each acting convex cocompactly on a fixed properly
  convex domain $\Omega$ with pairwise disjoint full orbital limit
  sets $\Lambda_\Omega(H_i)$.

  We denote by $\Heq{\dee \Omega}$ the space obtained from
  $\dee \Omega$ by collapsing all of the full orbital sets
  $\Lambda_\Omega(H_i)$ to points. Similarly, for $x \in \dee \Omega$,
  or a subset $\Lambda \subseteq \dee \Omega$, we use $\Heq{x}$ and
  $\Heq{\Lambda}$ to denote the images of $x$ and $\Lambda$ in
  $\Heq{\dee \Omega}$.
\end{defn}

When $\mc{H}$ is a conjugacy-invariant collection of subgroups of a
group $\Gamma \subseteq \Aut(\Omega)$, the action of $\Gamma$ on
$\dee \Omega$ descends to an action on $\Heq{\dee \Omega}$. More
generally, if $\Lambda \subseteq \dee \Omega$ is $\Gamma$-invariant,
$\Gamma$ also acts on $\Heq{\Lambda}$.

We show the following:
\begin{thm}
  \label{thm:bowditch_embedding_implies_conv_cocpct}
  Let $\Gamma \subseteq \PGL(d, \R)$ act on a properly convex domain
  $\Omega$, and suppose that $\Gamma$ is hyperbolic relative to a
  family of subgroups $\mc{H} = \{H_i\}$, such that the $H_i$ each act
  convex cocompactly on $\Omega$ with pairwise disjoint full orbital
  limit sets.

  If there is a boundary-convex $\Gamma$-invariant subset
  $\Lambda \subseteq \dee \Omega$ containing all of its faces, and a
  $\Gamma$-equivariant embedding
  $\dee(\Gamma, \mc{H}) \to \Heq{\dee \Omega}$ with image
  $\Heq{\Lambda}$, then $\Gamma$ acts convex cocompactly on $\Omega$
  and $\Lambda$ is the full orbital limit set
  $\Lambda_\Omega(\Gamma)$.
\end{thm}

\begin{remark}
  In Theorem \ref{thm:bowditch_embedding_implies_conv_cocpct}, we do
  not need to assume that $\Gamma$ is discrete in $\PGL(d,\R)$: this
  will also follow from the existence of the equivariant boundary
  embedding.
\end{remark}

There are two special cases of Theorem
\ref{thm:bowditch_embedding_implies_conv_cocpct} worth considering,
which we state separately as corollaries. The first is when the subset
$\Lambda$ is the entire boundary $\dee \Omega$.
\begin{cor}
  Let $\Gamma, \Omega$, and $\mc{H}$ be as in Theorem
  \ref{thm:bowditch_embedding_implies_conv_cocpct}, and suppose that
  $\dee(\Gamma, \mc{H})$ is equivariantly homeomorphic to $\Heq{\dee
    \Omega}$. Then $\Gamma$ divides $\Omega$.
\end{cor}

The second corollary is when the set of peripheral subgroups is empty,
i.e. $\Gamma$ is hyperbolic. 
\begin{cor}
  \label{cor:boundary_embedding_implies_cc}
  Let $\Gamma$ be a word-hyperbolic group in $\PGL(d, \R)$ acting on a
  properly convex domain $\Omega$, and suppose that the Gromov
  boundary of $\Gamma$ embeds equivariantly into $\dee \Omega$ with
  image $\Lambda$.

  If $\Lambda$ is boundary-convex and contains all of its faces, then
  $\Gamma$ acts convex cocompactly on $\Omega$ and
  $\Lambda = \Lambda_\Omega(\Gamma)$.
\end{cor}

When a hyperbolic group acts convex cocompactly on a domain $\Omega$,
its full orbital limit set contains no segments. So in this case,
$\Lambda$ contains all of its faces whenever no point of $\Lambda$
lies in the interior of any segment in $\dee \Omega$.

We also can phrase this corollary in terms of Anosov boundary
maps. Due to \cite[Theorem 1.15]{dgk2017convex} (see also
\cite[Theorem 1.10]{zimmer2017projective}), if a word-hyperbolic group
$\Gamma$ acts convex cocompactly on some domain $\Omega$, then the
inclusion map $\Gamma \hookrightarrow \PGL(d, \R)$ is a $P_1$-Anosov
representation preserving $\Omega$, and in this case the full orbital
limit set is the image of the Anosov boundary map
$\dee \Gamma \to \RP^{d-1}$. Thus Corollary
\ref{cor:boundary_embedding_implies_cc} implies:
\begin{cor}
  Let $\Gamma$ be a word-hyperbolic subgroup of $\PGL(d, \R)$
  preserving a properly convex domain $\Omega$, and suppose that there
  exists a $\Gamma$-equivariant embedding
  $\xi:\dee \Gamma \to \dee \Omega$ whose image is boundary-convex and
  contains all of its faces. Then the inclusion
  $\Gamma \hookrightarrow \PGL(d,\R)$ is a $P_1$-Anosov representation
  with $\RP^{d-1}$ boundary map $\xi$.
\end{cor}

Note that it is not true in general that the $\RP^{d-1}$-boundary map
$\xi$ of a $P_1$-Anosov representation always embeds into the boundary
of some properly convex domain $\Omega \subset \RP^{d-1}$. Moreover
even if $\xi$ does embed into $\dee \Omega$ for some $\Omega$, it does
not necessarily follow that the image of the embedding is
boundary-convex. However, it again follows from \cite[Theorem
1.15]{dgk2017convex} that in this case, there is \emph{some} $\Omega'$
such that $\xi$ embeds $\dee \Gamma$ into $\dee \Omega'$ with a
boundary-convex image containing its faces. In fact it is always
possible to take $\Omega'$ strictly convex with $C^1$ boundary.

\begin{remark}
  Kapovich-Leeb \cite{kl2018relativizing} and Zhu
  \cite{zhu2019relatively} (see also Zhu-Zimmer \cite{ZZ}) have given
  several possible definitions for a relative Anosov representation of
  a relatively hyperbolic group, aiming to generalize geometrical
  finiteness in rank-one in the same way that Anosov representations
  generalize convex cocompactness.

  The non-hyperbolic convex cocompact group actions we consider in
  this paper are not covered by either the Kapovich-Leeb or Zhu
  pictures. For one, not all examples of convex cocompact groups are
  relatively hyperbolic (see Section~\ref{subsec:examples}). But even
  in the relatively hyperbolic setting, the definitions are not
  compatible. Due to \cite[Proposition 10.3]{dgk2017convex}, convex
  cocompact groups in $\PGL(d,\R)$ do not contain weakly unipotent
  elements. In contrast, the relative Anosov subgroups considered by
  Kapovich-Leeb and Zhu always contain weakly unipotent elements if
  the group is not hyperbolic (see section 5 in
  \cite{kl2018relativizing}). However, see \cite{Weisman22} for
  follow-up work which gives a unified approach to studying both
  relatively hyperbolic convex cocompact groups in $\PGL(d,\R)$ and
  relative Anosov representations.
\end{remark}

During the proof of Theorem
\ref{thm:bowditch_embedding_implies_conv_cocpct}, we will see the
following (see Proposition \ref{prop:no_segments_in_non_peripherals}):
\begin{prop}
  In the setting of Theorem
  \ref{thm:bowditch_embedding_implies_conv_cocpct}, every nontrivial
  segment in the set $\Lambda$ is contained in the full orbital limit
  set of some $H_i \in \mc{H}$.
\end{prop}

This leads us to a converse to Theorem
\ref{thm:bowditch_embedding_implies_conv_cocpct}.

\begin{thm}
  \label{thm:no_segment_implies_rel_hyperbolic}
  Let $\Gamma$ be a group acting convex cocompactly on a properly
  convex domain $\Omega$, and suppose that $\Gamma$ has a
  conjugacy-invariant collection of subgroups $\mc{H} = \{H_i\}$, such
  that the groups in $\mc{H}$ lie in finitely many conjugacy classes
  and each $H_i$ acts convex cocompactly on $\Omega$.

  Then $\Gamma$ is hyperbolic relative to $\mc{H}$ if and only if
  \begin{enumerate}[label=(\roman*)]
  \item \label{item:pairwise_disjoint_limit_sets} the full orbital
    limit sets $\Lambda_\Omega(H_i)$, $\Lambda_\Omega(H_j)$ are
    disjoint for distinct $H_i, H_j \in \mc{H}$,
  \item \label{item:no_non_peripheral_segments} every nontrivial
    segment in $\Lambda_\Omega(\Gamma)$ is contained in
    $\Lambda_\Omega(H_i)$ for some $H_i \in \mc{H}$, and
  \item \label{item:self_normalizing} each $H_i \in \mc{H}$ is its own
    normalizer in $\Gamma$.
  \end{enumerate}
  Moreover, in this case, $\dee(\Gamma, \mc{H})$ equivariantly embeds
  into $\Heq{\dee \Omega}$ with image $\Heq{\Lambda_\Omega(\Gamma)}$.
\end{thm}
\begin{remark}
  If conditions \ref{item:pairwise_disjoint_limit_sets} and
  \ref{item:no_non_peripheral_segments} hold for a conjugacy-closed
  collection of subgroups $\mc{H}$ of $\Gamma$, then they also hold
  for the collection of normalizers, since
  $g \cdot \Lambda_\Omega(H_i) = \Lambda_\Omega(gH_ig^{-1})$ for any
  $H_i \in \mc{H}$.
  
  Moreover, condition \ref{item:self_normalizing} of Theorem
  \ref{thm:no_segment_implies_rel_hyperbolic} is always true for the
  peripheral subgroups of a relatively hyperbolic group, because then
  each $H_i \in \mc{H}$ can be exactly realized as the stabilizer of
  its unique fixed point in the Bowditch boundary
  $\dee(\Gamma, \mc{H})$ (see Theorem \ref{thm:yaman_relative_hyp}).
\end{remark}

Islam and Zimmer \cite{iz2019flat, iz2019convex} have previously shown
that when $\Gamma$ is a convex cocompact group which is hyperbolic
relative to a collection $\mc{H}$ of virtually abelian subgroups of
rank at least 2, conditions \ref{item:pairwise_disjoint_limit_sets}
and \ref{item:no_non_peripheral_segments} of Theorem
\ref{thm:no_segment_implies_rel_hyperbolic} hold, and moreover the
assumption that all of the groups in $\mc{H}$ act convex cocompactly
is automatically satisfied. In particular, this implies that the set
$\Heq{\dee \Omega}$ is well-defined.

Thus, in this case, Theorem
\ref{thm:no_segment_implies_rel_hyperbolic} implies the following:
\begin{cor}
  \label{cor:free_abelian_corollary}
  Let $\Omega$ be a properly convex domain, and let $\Gamma$ be a
  group which is hyperbolic relative to a collection $\mc{H}$ of
  virtually abelian subgroups with rank at least 2.

  If $\Gamma$ acts convex cocompactly on $\Omega$, then there is an
  equivariant embedding from $\dee(\Gamma, \mc{H})$ to
  $\Heq{\dee \Omega}$ whose image is $\Heq{\Lambda_\Omega(\Gamma)}$.
\end{cor}

Islam-Zimmer have recently shown in \cite{iz2020convex} that if
$\Gamma$ acts convex cocompactly on some domain $\Omega$, and $\Gamma$
is isomorphic to the fundamental group of an orientable, closed, and
irreducible 3-manifold $M$ which is non-geometric (i.e. does not carry
one of the eight Thurston geometries), then $M$ decomposes into
hyperbolic pieces and therefore $\Gamma$ is hyperbolic relative to
free abelian subgroups. Using different methods, they also give a
proof of Corollary \ref{cor:free_abelian_corollary} in this case.

Theorem \ref{thm:no_segment_implies_rel_hyperbolic} in fact applies to
a strictly larger family of groups than those covered by Corollary
\ref{cor:free_abelian_corollary}: there exist relatively hyperbolic
convex cocompact groups which are \emph{not} hyperbolic relative to
virtually abelian subgroups (and are not themselves hyperbolic). See
the end of Section \ref{subsubsec:non_hyp_convex_cocompact} for
details.

\begin{remark}
  Since this paper first appeared, Islam-Zimmer \cite{iz2022structure}
  have shown that when $\Gamma$ is \emph{any} relatively hyperbolic
  group acting convex cocompactly on a properly convex domain
  $\Omega$, then each peripheral subgroup of $\Gamma$ \emph{always}
  acts convex cocompactly on $\Omega$ as well. This makes it possible
  to drop an assumption from Theorem
  \ref{thm:no_segment_implies_rel_hyperbolic}, and implies that
  Corollary \ref{cor:free_abelian_corollary} holds in the case of
  \emph{arbitrary} peripheral subgroups.

  In the same paper, Islam-Zimmer extend some of our techniques to
  show versions of Theorem \ref{thm:no_segment_implies_rel_hyperbolic}
  and Corollary \ref{cor:free_abelian_corollary} in the more general
  context of \emph{naive} convex cocompact group actions in projective
  space, which we do not discuss in this paper.
\end{remark}

\begin{remark}
  Relatively hyperbolic group actions on convex projective domains
  have previously been studied by Crampon-Marquis \cite{clm2016convex}
  (who provide a notion of a geometrically finite group action on a
  strictly convex domain), and by Cooper-Long-Tillmann
  \cite{clt2015convex} in the context of convex projective cusps. Choi
  \cite{choi2011} has also studied the relationship between the
  projective geometry of strictly convex projective orbifolds with
  ends, and the Bowditch boundaries of their (relatively hyperbolic)
  fundamental groups.

  The relatively hyperbolic group actions we consider in this paper
  are different, however, since we consider projective orbifolds which
  are \emph{not} necessarily strictly convex. In this context, the
  peripheral subgroups do \emph{not} need to be associated to an end
  of the manifold, and act \emph{cocompactly} on a convex subset of
  projective space.
\end{remark}

\subsection{Outline of the paper}

In Section \ref{sec:definitions_notation}, we recall some background
about properly convex domains and their automorphism groups, and
describe some known examples of convex cocompact groups.

Sections \ref{sec:expansion_implies_cocompact} through
\ref{sec:cocompactness_implies_expansion} of this paper are devoted to
the proof of Theorem~\ref{thm:convex_cocompact_equals_expansion}. We
first prove the implication (\ref{item:main_thm_2}) $\implies$
(\ref{item:main_thm_1}) in
Section~\ref{sec:expansion_implies_cocompact}, using a modified
version of an analogous argument in \cite{dgk2017convex}. The proof
proceeds by contradiction: we assume that the space of $\Gamma$-orbits
is not compact in $\Cor(\Gamma)$, which means that the orbits
uniformly accumulate on the boundary of $\Cor(\Gamma)$ in
$\Omega$. But the expansivity of the action on the faces of $\Lambda$
allows us to move points away from the boundary of $\Cor(\Gamma)$ and
we get a contradiction. This basic idea is essentially due to Sullivan
\cite{sullivan1985quasiconformal}, but here the argument is more
complicated because the expansivity occurs in a different space than
the one where we want to find a cocompact action.

\begin{remark}
  In \cite{klpdomains}, Kapovich-Leeb-Porti also adapt Sullivan's
  argument to show that an expanding (Anosov) action in one flag
  manifold can produce a cocompact domain of discontinuity in another
  flag manifold. Their argument does not apply directly to our
  situation, since it relies on the fact that the ``limit set'' where
  the expanding action occurs is \emph{compact}. Typically, if
  $\Lambda$ is the full orbital limit set of a group acting convex
  cocompactly on $\Omega \subset \RP^{d-1}$, the subset of $\Gr(k, d)$
  consisting of supports of $k$-dimensional faces of the full orbital
  limit set $\Lambda$ is \emph{not} compact (see
  e.g. \cite{benoist2006convexes} for examples). Consequently, it is
  necessary for us to work with expansion in several different
  Grassmannians for the same group---and to exploit the fact that the
  full orbital limit set $\Lambda$ actually lies in the boundary of a
  convex subset of $\RP^{d-1}$.
\end{remark}

Section~4 of the paper is mainly devoted to the proof of
Proposition~\ref{prop:relative_benzecri}, the key technical result
mentioned earlier in \ref{sec:relative_benzecri_comments}.  Then, in
Section \ref{sec:cocompactness_implies_expansion}, we use the results
of Section \ref{sec:spaces_of_domains} to prove the implication
(\ref{item:main_thm_1}) $\implies$ (\ref{item:main_thm_2}) of Theorem
\ref{thm:convex_cocompact_equals_expansion}. The basic argument is
again inspired by Sullivan's work \cite{sullivan1979density} in the
case of convex cocompact actions in $\H^d$.

The remainder of the paper (Sections~\ref{sec:relative_hyp_background}
through \ref{sec:relative_hyp_conv_cocompact}) restricts to the
context of relatively hyperbolic groups. In Section
\ref{sec:relative_hyp_background}, we give some background on
relatively hyperbolic groups, and state a dynamical definition of
relatively hyperbolic groups due to Yaman
\cite{yaman2004topological}. Then in Section
\ref{sec:bowditch_bdry_embedding}, we use this characterization to
prove Theorem \ref{thm:bowditch_embedding_implies_conv_cocpct},
relying on the characterization of convex cocompactness given by
Theorem~\ref{thm:convex_cocompact_equals_expansion}.

In Section \ref{sec:relative_hyp_conv_cocompact}, we extend a result
of Islam-Zimmer \cite{iz2019convex}, showing that when $\Gamma$ is a
group acting on a properly convex domain $\Omega$, and $\Gamma$ is
hyperbolic relative to a collection of subgroups acting convex
cocompactly on $\Omega$, conditions
\ref{item:pairwise_disjoint_limit_sets} and
\ref{item:no_non_peripheral_segments} of Theorem
\ref{thm:no_segment_implies_rel_hyperbolic} hold. Then, we prove the
rest of Theorem \ref{thm:no_segment_implies_rel_hyperbolic}, again by
applying Theorem~\ref{thm:convex_cocompact_equals_expansion} and using
the dynamical definition of relative hyperbolicity from Section
\ref{sec:relative_hyp_background}.

\subsection{Acknowledgments}

The author would like to thank his advisor, Jeff Danciger, for
providing endless amounts of guidance and encouragement. We also would
like to thank Fanny Kassel for giving useful feedback, and two
anonymous referees for helpful comments and insights.



\section{Background on properly convex domains}
\label{sec:definitions_notation}

\subsection{Basic definitions}

All real vector spaces in this paper are finite-dimensional.

\subsubsection{Convex cones and convex domains}

\begin{defn}
  \label{defn:convex_cone}
  A \emph{convex cone} in a real vector space $V$ is a convex subset
  of $V - \{0\}$ which is closed under multiplication by positive
  scalars.

  A convex cone is \emph{sharp} if it does not contain any affine
  line.
\end{defn}

The \emph{boundary} of a convex cone $C$ in a real vector space $V$ is
the boundary of $C$ viewed as a cone in its linear span $V'$ in $V$;
this boundary is homeomorphic to a cone over $S^{k-2}$, where
$k = \dim V'$.

\begin{defn}
  \label{defn:convex_domain}
  A subset $\Omega \subset \P(V)$ is \emph{convex} if it is the
  projectivization of some convex cone $\tilde{\Omega} \subset V$, and
  it is \emph{properly convex} if $\tilde{\Omega}$ is sharp
  (equivalently, if $\overline{\Omega}$ is contained in some affine
  chart in $\P(V)$). An open convex set is a \emph{convex domain}.

  The boundary $\dee \Omega$ is the projectivization of
  $\dee \tilde{\Omega} - \{0\}$. A convex set $\Omega$ is
  \emph{strictly convex} if $\dee \Omega$ does not contain a
  nontrivial projective segment.
\end{defn}

\begin{defn}
  \label{defn:supporting_subspace}
  Let $\Omega$ be a convex subset of $\P(V)$. A \emph{supporting
    subspace} of $\Omega$ is a projective subspace which intersects
  $\dee \Omega$ but not $\Omega$. In particular, a \emph{supporting
    hyperplane} is a codimension-1 supporting subspace.
\end{defn}

Using the convex cone over $\Omega$, one can easily verify the
following:
\begin{prop}
  \label{prop:supporting_hyperplanes_exist}
  Let $\Omega$ be a convex subset of $\P(V)$. Every point
  $x \in \dee \Omega$ is contained in at least one supporting
  hyperplane.
\end{prop}

We remark that a convex domain in $\P(V)$ has $C^1$ boundary precisely
when every point in $\dee \Omega$ is contained in \emph{exactly} one
supporting hyperplane.

\subsubsection{Projective line segments}

When $\Omega$ is a properly convex set, and
$x, y \in \overline{\Omega}$, we use $[x,y]$ to denote the unique
(closed) projective line segment joining $x$ and $y$ which is
contained in $\overline{\Omega}$. We similarly use $(x,y)$, $[x,y)$,
$(x,y]$ to denote open and half-open projective line segments.

\subsubsection{Convex hull and ideal boundary}

\begin{defn}
  If $\Omega$ is a properly convex set in $\P(V)$ and
  $\Lambda \subset \dee \Omega$, then the \emph{convex hull} of
  $\Lambda$ is its convex hull in $\Omega$ in any affine chart
  containing $\Omega$. We denote the convex hull of $\Lambda$ by
  $\CH(\Lambda)$.
  
  The \emph{ideal boundary} of a set $C$ in a properly convex set
  $\Omega$ is the set
  \[
    \dee_i(C) := \overline{C} \cap \dee \Omega,
  \]
  where the closure of $C$ is taken in $\P(V)$.
\end{defn}

If $\Lambda$ is a subset of $\dee \Omega$ with nonempty convex hull,
$\Lambda$ is \emph{boundary-convex} (Definition
\ref{defn:boundary_convex}) precisely when
\[
  \Lambda = \dee_i \CH(\Lambda).
\]

\subsection{Faces in $\dee \Omega$}
\label{subsec:stratum_background}

\begin{defn}
  Let $\Omega$ be a properly convex domain. The \emph{face} of
  $\dee \Omega$ at a point $x$, which we denote $\strat(x)$, is the
  set of points $y \in \dee \Omega$ such that $x$ and $y$ lie in an
  open segment $(a,b) \subset \dee \Omega$.

  The \emph{dimension} of a face $F$ is the dimension of a minimal
  projective subspace containing $F$; such a minimal subspace is
  called the \emph{support} of the face and is denoted $\supp(F)$.
\end{defn}

A face is always a convex subset of projective space, open in its
support. A face is a closed subset of $\dee \Omega$ if and only if it
is an extreme point of $\Omega$.

\begin{remark}
  Earlier versions of \cite{dgk2017convex} (and of this paper)
  referred to what we call a \emph{face} as a ``stratum.'' Our current
  definition of face agrees with the definition used by Islam and
  Zimmer. Notably, our faces are \emph{not} the same as the
  \emph{facettes} of Benoist and Benz\'ecri.

  In particular, our definition ensures that every face is relatively
  open in its support, and that every point in the boundary of a
  properly convex domain $\Omega$ is contained in some face.
\end{remark}

\begin{defn}
  When $\Lambda$ is a subset of $\dee \Omega$, a \emph{face} of
  $\Lambda$ in $\Omega$ is a face of $\dee \Omega$ which intersects
  $\Lambda$ nontrivially.
\end{defn}

We warn the reader that this definition of ``face'' in $\Lambda$
depends on both $\Lambda$ and on the domain $\Omega$ whose boundary
contains $\Lambda$. Often, we will not need to worry about this, due
to the following consequence of Lemma 4.1 (1) in \cite{dgk2017convex}:
\begin{lem}
  \label{lem:full_orbital_limit_set_properties}
  Let $\Omega$ be a properly convex domain, and let $\Gamma$ act
  convex cocompactly on $\Omega$. The full orbital limit set
  $\Lambda_\Omega(\Gamma)$ is closed and boundary-convex, and contains
  all of its faces.
\end{lem}

\subsection{The Hilbert metric}

Here we recall the definition of the Hilbert metric, a useful tool for
understanding group actions on properly convex domains. See
e.g. \cite{marquis2013around} for more background.

Given four distinct points $a, b, c, d$ in $\RP^1$ (or four points in
$\RP^{d-1}$ lying on a single projective line), recall that the
\emph{cross-ratio} $[a, b; c,d]$ is given by
\[
  [a, b; c,d] := \frac{|c - a| \cdot |d - b|}{|b - a| \cdot |d - c|},
\]
where the distances are taken in any Euclidean metric on an affine
chart containing $a,b,c,d$.

The cross-ratio is a projective invariant on 4-tuples, and in fact it
parameterizes the space of $\PGL(2, \R)$-orbits of distinct 4-tuples
in $\RP^1$.

\begin{defn}
  Let $\Omega \subset \P(V)$ be a properly convex domain. The
  \emph{Hilbert metric}
  \[
    d_\Omega(\cdot, \cdot):\Omega^2 \to \R_{\ge 0}
  \]
  is given by the formula
  \[
    d_\Omega(x,y) = \frac{1}{2}\log[a, x; y, b],
  \]
  where $a, b$ are the two points in $\dee \Omega$ such that
  $a, x, y, b$ lie on a projective line in that order.
\end{defn}

When the domain $\Omega$ is an ellipsoid of dimension $d$, the Hilbert
metric on $\Omega$ recovers the familiar Klein model for hyperbolic
space $\H^d$. More generally we have the following:
\begin{prop}[{See for example \cite[Section 28]{BK53}}]
  Let $\Omega$ be a properly convex domain. Then:
  \begin{enumerate}
  \item The pair $(\Omega, d_\Omega)$ is a proper metric space.
  \item \label{item:hilb_geodesic} If $x$ and $y$ are in $\Omega$,
    then $[x,y]$ is the image of a geodesic (with respect to
    $d_\Omega$) joining $x$ and $y$.
  \item The group $\Aut(\Omega)$ acts by isometries of $d_\Omega$.
  \end{enumerate}
\end{prop}
This implies that $\Aut(\Omega)$ always acts \emph{properly} on
$\Omega$. In particular, a subgroup of $\Aut(\Omega)$ is discrete in
$\PGL(V)$ if and only if it acts properly discontinuously on $\Omega$.

Part (\ref{item:hilb_geodesic}) of the above Proposition means that
$(\Omega, d_\Omega)$ is always a geodesic metric space. However, in
general it need not be uniquely geodesic---this is one of many ways in
which the geometry on a properly convex domain equipped with its
Hilbert metric can differ from hyperbolic geometry.

The point of the Hilbert metric is that it allows us to understand
many aspects of group actions on convex projective domains in terms of
metric geometry; in particular, we may apply the \v{S}varc-Milnor
lemma when we have a convex cocompact action on a domain.

The Hilbert metric can also be used to characterize faces in
$\dee \Omega$. An easy calculation shows:
\begin{prop}
  \label{prop:strata_hilbert_bounded}
  Let $\Omega$ be a properly convex domain, let $x \in \dee \Omega$,
  and fix points $p_1, p_2 \in \Omega$. For any $y \in \dee \Omega$,
  we have $y \in \strat(x)$ if and only if the Hausdorff distance
  (with respect to $d_\Omega$) between $[p_1, x)$ and $[p_2,y)$ is
  finite.
\end{prop}

Since $[p_1, x)$ and $[p_2, y)$ are the images of geodesic rays in
$(\Omega, d_\Omega)$, the above is equivalent to the condition that,
if $c_x$, $c_y$ are unit-speed geodesic rays in $(\Omega, d_\Omega)$
following projective line segments from $p_1$, $p_2$ to $x, y$,
respectively, then
\[
  d_\Omega(c_x(t), c_y(t)) \le k
\]
for some fixed $k$ independent of $t \in \R_{\ge 0}$.

\subsection{Properly embedded simplices}

\begin{defn}
  A \emph{projective $k$-simplex in $\RP^{d-1}$} is the
  projectivization of the positive linear span of $k+1$ linearly
  independent vectors in $\Rd$.
\end{defn}

A projective $k$-simplex $\Delta$ is an example of a properly convex
set in $\RP^{d-1}$. If $\Delta$ is the span of standard basis vectors
$e_1, \ldots, e_d$, the group $D^+ \subset \PGL(d, \R)$ of
projectivized diagonal matrices with positive entries (isomorphic to
$\R^{d-1}$) acts simply transitively on $\Delta$. Then, any discrete
$\Z^{d-1}$ subgroup of $D^+$ acts properly discontinuously and
cocompactly on $\Delta$, so the \v{S}varc-Milnor lemma implies that
$(\Delta, d_\Delta)$ is quasi-isometric to Euclidean space $\E^{d-1}$.

\begin{defn}
  Let $\Omega$ be a properly convex domain. A convex projective
  simplex $\Delta \subset \Omega$ is \emph{properly embedded} if
  $\dee \Delta$ is contained in $\dee \Omega$.

  A properly embedded simplex in $\Omega$ gives an isometric embedding
  \[
    (\Delta, d_\Delta) \to (\Omega, d_\Omega),
  \]
  which in turn gives a quasi-isometric embedding
  \[
    \E^k \to (\Omega, d_\Omega).
  \]
\end{defn}

\emph{Maximal} properly embedded simplices in $\Omega$ can be thought
of as analogues of \emph{maximal flats} in $\mr{CAT}(0)$ spaces; see
e.g. \cite{benoist2006convexes}, \cite{iz2019flat},
\cite{iz2019convex}, \cite{bobb2020codimension}. However, in general,
the metric space $(\Omega, d_\Omega)$ is not $\mr{CAT}(0)$; in fact
this occurs if and only if $\Omega$ is an ellipsoid
\cite{ks1958curvature}.

\subsection{Duality for convex domains}

Let $V$ be a real vector space. Given a convex set
$\Omega \subset \P(V)$, it is often useful to consider the \emph{dual
  convex set} $\Omega^* \subset \P(V^*)$.

\begin{defn}
  Let $C$ be a convex cone in a real vector space $V$. The \emph{dual
    convex cone} $C^* \subset V^* - \{0\}$ is
  \[
    C^* = \{\alpha \in V^* : \alpha(x) > 0 \textrm{ for all } x \in
    \overline{C} - \{0\}\}.
  \]
\end{defn}

The following is easily verified:
\begin{prop}
  \label{prop:duality_properties}
  Let $C$ be a convex cone in a real vector space $V$.
  \begin{enumerate}
  \item $C^*$ is a convex cone in $V^* - \{0\}$.
  \item $C^{**} = C$, under the canonical identification $V^{**} = V$.
  \item\label{item:sharp_open} $C^*$ is sharp if and only if $C$ has
    nonempty interior.
  \end{enumerate}
\end{prop}

If $\Omega$ is the projectivization of a convex cone in $\P(V)$, the
\emph{dual convex set} is the projectivization $\Omega^*$ of
$\tilde{\Omega}^*$, where $\tilde{\Omega}$ is any cone over
$\Omega$. When $\Omega$ is a properly convex domain in $\P(V)$,
$\Omega^*$ is a properly convex domain in $\P(V^*)$.

In general, $\Omega^*$ need not be projectively equivalent to
$\Omega$. However, the features of $\Omega$ affect the features of
$\Omega^*$. For instance, $\Omega$ is strictly convex if and only if
the boundary of $\Omega^*$ is $C^1$ (and vice versa, since
$\Omega^{**}$ is naturally identified with $\Omega$). We also note
that duality reverses inclusions of convex sets.

If $\Gamma$ is a subgroup of $\PGL(V)$ preserving $\Omega$, the dual
action of $\Gamma$ on $\P(V^*)$ preserves the dual domain
$\Omega^*$. So we can simultaneously view $\Gamma$ as a subgroup of
$\Aut(\Omega)$ and $\Aut(\Omega^*)$.

\begin{remark}
  \label{rem:duality_convention}
  Throughout this paper, we will consistently view the dual projective
  space as a space of projective hyperplanes in $\P(V)$. That is, we
  identify an element $[\alpha] \in \P(V^*)$ with the projective
  hyperplane $\P(\ker \alpha) \subset \P(V)$.
\end{remark}

\subsection{Known examples of convex cocompact groups}
\label{subsec:examples}

Groups with convex cocompact actions in projective space fall into two
main classes: the word-hyperbolic groups and the non-word-hyperbolic
groups. A consequence of \cite{dgk2017convex} is that a group $\Gamma$
acting convex cocompactly on some domain $\Omega$ is word-hyperbolic
if and only if the full orbital limit set of $\Gamma$ does not contain
a nontrivial projective segment. In particular this always holds if
$\Omega$ itself is strictly convex.

Some of the examples we list below are examples of groups
\emph{dividing} domains, meaning that the group $\Gamma$ acts
cocompactly on the entire domain $\Omega$. See
\cite{benoist2008survey} for a survey on the topic of convex divisible
domains.

\subsubsection{Strictly convex divisible examples}
The simplest example of a strictly convex divisible domain is
hyperbolic space $\H^d$. Uniform lattices in $\PO(d,1)$ exist in any
dimension $d \ge 1$, so they are examples of hyperbolic groups acting
cocompactly on the projective model for $\H^d$ (a round ball in
$\RP^d$).

A torsion-free uniform lattice in $\PO(d,1)$ can be viewed as the
image of the holonomy representation of a closed hyperbolic
$d$-manifold $M$. Viewing $\PO(d, 1)$ as a subgroup of $\PGL(d+1,\R)$
allows us to view the hyperbolic structure on $M$ as a convex
projective structure. It is sometimes possible to perturb the subgroup
$\Gamma \simeq \pi_1M$ inside $\PGL(d+1,\R)$ to obtain a new discrete
group $\Gamma' \simeq \pi_1M$ in $\PGL(d+1, \R)$ which is the holonomy
of a different convex projective structure on $M$. The deformed group
$\Gamma'$ acts cocompactly on some properly convex domain $\Omega'$,
which in general is not projectively equivalent to a round ball.

Further examples of groups $\Gamma$ dividing strictly convex domains
have been found by Benoist \cite{benoist2006convexeshyp} in dimension
4, using reflection groups, and Kapovich \cite{kapovich2007convex} in
dimensions $d \ge 4$, by finding convex projective structures on
Gromov-Thurston manifolds \cite{gt1987pinching}. The Benoist and
Kapovich examples share the feature that the dividing group $\Gamma$
is not isomorphic to any lattice in $\PO(d,1)$---while $\Gamma$ is
word-hyperbolic, the quotient orbifold $\Omega / \Gamma$ carries no
hyperbolic structure.

\subsubsection{Non-strictly convex divisible examples}

The simplest examples of non-strictly convex divisible domains are
projective $k$-simplices, which are divided by free abelian groups of
rank $k$. When $k \ge 2$, these simplices are not strictly convex, but
we can still decompose the action into strictly convex pieces---the
cone over the simplex splits as a sum of strictly convex cones, and
each $\Z$ factor acts cocompactly on a summand. So we may wish to find
\emph{irreducible} examples. These exist too: uniform lattices in
$\SL(d,\R)$ act cocompactly on the symmetric space
$\SL(d,\R) / \SO(d)$. This symmetric space can be modeled as the
projectivization of the set of positive definite symmetric matrices
sitting inside the space of $d \times d$ matrices. This set is convex,
but not strictly convex whenever $d > 2$.

Other interesting examples of non-strictly convex divisible domains
have been discovered. In 2006, Benoist \cite{benoist2006convexes}
produced examples of \emph{inhomogeneous} properly convex divisible
domains in dimensions $3-7$, divided by non-hyperbolic groups; other
examples in dimensions $4-7$ were later found by Choi-Lee-Marquis in
\cite{clm2016convex}.

Both of these families of examples essentially come from the theory of
reflection groups. Given a Coxeter group $\Gamma$ acting by
reflections in $\PGL(d,\R)$, there is a fairly straightforward
procedure due to Vinberg \cite{vinberg1971discrete} which determines
whether or not $\Gamma$ acts cocompactly on some convex domain in
projective space. The domain fails to be strictly convex if and only
if $\Gamma$ contains virtually abelian subgroups of rank at least $2$
(which happens only when $\Gamma$ contains a Coxeter subgroup of type
$\tilde{A_n}$).

More recently, Blayac-Viaggi \cite{BV} have constructed additional
examples of irreducible non-strictly convex divisible domains in
$\P(\R^d)$, for any dimension $d \ge 4$. The construction does not use
reflection groups, but rather a combination of arithmetic methods and
a procedure known as \emph{projective bending}.

\subsubsection{Non-hyperbolic convex cocompact groups}
\label{subsubsec:non_hyp_convex_cocompact}

In \cite{bdl2015convex}, Ballas-Danciger-Lee produce examples of
non-hyperbolic groups acting convex cocompactly which do not divide a
properly convex domain. These come from deformations of hyperbolic
structures on certain cusped hyperbolic 3-manifolds.

None of the groups $\Gamma$ in the examples in
\cite{benoist2006convexes}, \cite{clm2016convex}, or
\cite{bdl2015convex} are word-hyperbolic, but they are all
\emph{relatively} hyperbolic, relative to a family of virtually
abelian subgroups of rank at least $2$. This situation was studied
more generally by Islam-Zimmer in \cite{iz2019flat},
\cite{iz2019convex}. Islam-Zimmer show that if $\Gamma$ is hyperbolic
relative to virtually abelian subgroups of rank $\ge 2$, and $\Gamma$
acts convex cocompactly on a properly convex domain $\Omega$, the
peripheral subgroups of $\Gamma$ act cocompactly on properly embedded
projective simplices in $\Omega$. Moreover, in this situation, the
properly embedded maximal simplices in $\Cor(\Gamma)$ of dimension at
least $2$ are \emph{isolated}, and every such maximal simplex has
compact quotient by a free abelian subgroup of $\Gamma$.

Work of Danciger, Guéritaud, Kassel, Lee, and Marquis
\cite{dgklm2021convex} shows that in fact \emph{every} convex
cocompact reflection group is either hyperbolic or relatively
hyperbolic relative to virtually abelian subgroups. But, this is not
true for all non-hyperbolic groups with convex cocompact
actions. Uniform lattices in $\SL(d,\R)$ provide a counterexample,
since the maximal flat subspaces of the Riemannian symmetric space
$\SL(d,\R) / \SO(d)$ are \emph{not} isolated (and thus the properly
embedded maximal simplices in its projective model are not isolated).

There are also examples of convex cocompact relatively hyperbolic
groups which are \emph{not} hyperbolic relative to virtually abelian
subgroups. For instance, for every $d \ge 3$, the construction of
Blayac-Viaggi mentioned earlier yields groups dividing domains in
$\PGL(d, \R)$, which are hyperbolic relative to subgroups which are
virtually the product of an infinite cyclic group and the fundamental
group of a closed hyperbolic $(d-3)$-manifold.

Another construction of relatively hyperbolic convex cocompact groups
uses the following (not yet published) result of
Danciger-Guéritaud-Kassel:
\begin{prop}[{See \cite[Proposition 12.4]{dgk2017convex} for a
  statement}]
  \label{prop:dgk_combination}
  Let $\Gamma_1, \Gamma_2 \subset \PGL(V)$ be groups acting convex
  cocompactly in $\P(V)$, and suppose that $\Gamma_1$, $\Gamma_2$ both
  do not divide any nonempty properly convex open subset in
  $\P(V)$. Then for some $g \in \PGL(V)$, the group generated by
  $\Gamma_1$ and $g\Gamma_2g^{-1}$ is isomorphic to the free product
  $\Gamma_1 * \Gamma_2$, and acts convex cocompactly on $\P(V)$.
\end{prop}

As Danciger-Gu\'eritatud-Kassel observe, it is possible to use this
proposition to construct some exotic examples of convex cocompact
groups. Here we explain the procedure, assuming
Proposition~\ref{prop:dgk_combination} holds. Let $\Gamma_1$ and
$\Gamma_2$ be uniform lattices in $\SL(d, \R)$. The projective model
for the Riemannian symmetric space $\SL(d, \R) / \SO(d)$ is embedded
into $\P(V')$, where $V'$ is the vector space of $d \times d$ real
matrices. We can in turn embed $V'$ into some vector space $V$ so that
$V = V' \oplus V''$ for some complementary subspace $V''$ with
positive dimension. We obtain corresponding embeddings of $\Gamma_1$
and $\Gamma_2$ into $\SL(V)$ by asking for both of these groups to act
via the induced representation $\SL(d, \R) \to \SL(V')$ and trivially
on $V''$.

By \cite[Theorem 1.6(E)]{dgk2017convex}, $\Gamma_1$ and $\Gamma_2$ act
convex cocompactly in $\P(V)$. Further, since $\Gamma_1$ and
$\Gamma_2$ divide a domain in $\P(V')$, their virtual cohomological
dimension is equal to the dimension of $\P(V')$, which prevents either
group from dividing any larger-dimensional domain. So these groups
satisfy the hypotheses of Proposition~\ref{prop:dgk_combination} and
we can find a convex cocompact subgroup $\Gamma$ in $\P(V)$ isomorphic
to $\Gamma_1 * \Gamma_2$.

This (abstract) free product is relatively hyperbolic, relative to the
collection of conjugates of $\Gamma_1, \Gamma_2$. But, it is not
relatively hyperbolic relative to virtually abelian subgroups. One way
to see this is that the convex core of $\Gamma$ in $\P(V)$ must
contain the convex core of $\Gamma_1$, which is a copy of the
projective model for $\SL(d, \R) / \SO(d)$; this domain contains many
maximal properly embedded simplices which are not isolated.



\section{Expansion implies convex cocompactness}
\label{sec:expansion_implies_cocompact}

The goal of this section is to prove the implication
(\ref{item:main_thm_2}) $\implies$ (\ref{item:main_thm_1}) of Theorem
\ref{thm:convex_cocompact_equals_expansion}. First let us specify
exactly what we mean by ``expanding at the faces'' of a subset
$\Lambda \subset \dee \Omega$.

\subsection{Expansion on the Grassmannian}
\label{subsec:support_dynamics}

Recall that a continuous map $f:X \to X$ on a metric space $(X, d_X)$
is said to be \emph{$C$-expanding} on a subset $U \subset X$, for a
constant $C > 1$, if
\[
  d_X(f(x), f(y)) \ge C \cdot d_X(x,y)
\]
for all $x,y \in U$.

\begin{defn}
  \label{defn:expansion_in_supports}
  Let $\Omega$ be a properly convex domain in $\RP^{d-1}$, let
  $\Gamma \subset \PGL(d, \R)$ preserve $\Omega$, and let $\Lambda$ be
  a $\Gamma$-invariant subset of $\dee \Omega$.

  Fix a Riemannian metric $d_k$ on each Grassmannian $\Gr(k,d)$. We
  say that the action of $\Gamma$ on $\Omega$ is \emph{expanding at
    the faces of $\Lambda$} if, for every face $F$ of $\Lambda$, there
  is a constant $C > 1$, an element $\gamma \in \Gamma$, and an open
  subset $U \subset \Gr(k,d)$ with $\supp(F) \in U$ such that $\gamma$
  is $C$-expanding on $U$ (with respect to the metric $d_k$).

  If the constant $C > 1$ can be chosen uniformly for all faces $F$
  of $\Lambda$, then we say the action is \emph{$C$-expanding at the
    faces of $\Lambda$} or just \emph{uniformly expanding at the
    faces}.
\end{defn}

\begin{remark}
  \label{remark:support_expanding_well_defined}
  It is conceivable that a group action could be expanding at the
  faces of $\Lambda$ with respect to some choice of Riemannian metric
  $d_k$ on $\Gr(k,d)$, but not with respect to another.

  However, if $\Gamma$ is $C$-expanding with respect to $d_k$ for a
  uniform constant $C$, the choice of metric does not matter: since
  $\Gr(k,d)$ is compact, all Riemannian metrics on $\Gr(k,d)$ are
  bilipschitz-equivalent, and when $\Gamma$ is $C$-expanding at the
  faces of $\Lambda$, one can apply expanding elements iteratively to
  see that $\Gamma$ is also $C'$-expanding for an arbitrary constant
  $C'$.

  When the set of supports of $(k-1)$-dimensional faces of $\Lambda$
  is compact in $\Gr(k,d)$ for each $k$, then a $\Gamma$-action is
  expanding at faces with respect to some choice of metric $d_k$ if
  and only if it is uniformly expanding at faces with respect to that
  metric (and hence to every metric). For instance, this is the case
  when $\Lambda$ is compact and does not contain any nontrivial
  segments (so the set of faces is the same as the set of points).

  In our context, however, we will not be able to assume this kind of
  compactness. So, when we discuss expansion, we need to either
  specify the metric or assume that the expansion is uniform.
\end{remark}

Any metric on $\RP^{d-1}$ induces a metric on each $\Gr(k,d)$, by
viewing elements of $\Gr(k,d)$ as closed subsets of $\RP^{d-1}$ and
taking Hausdorff distance. From this point forward, we will assume
that $d_\P$ denotes the \emph{angle metric} on projective space, which
is induced by a choice of inner product on $\Rd$. Hausdorff distance
on $\Gr(k, d)$ (with respect to the angle metric) is a Riemannian
metric.

\begin{lem}
  \label{lem:subspace_achieving_dist_exists}
  Let $x \in \RP^{d-1}$, and let $W \in \Gr(k,d)$. There exists
  $V \in \Gr(k,d)$ so that $x \in V$ and
  \[
    d_\P(x,W) = d_H(V,W),
  \]
  where $d_\P$ is the angle metric on projective space, and $d_H$ is
  the metric induced on $\Gr(k,d)$ by Hausdorff distance.
\end{lem}
\begin{proof}
  If $x \in W$, then we can just take $V = W$, so assume that
  $d_\P(x,W) > 0$. The definition of Hausdorff distance immediately
  implies that for any $V$ containing $x$, $d_H(V, W) \ge d_\P(x, W)$,
  so we only need to find some $V$ satisfying the other bound. The
  diameter of projective space in the angle metric is $\pi/2$, which
  gives an upper bound on the Hausdorff distance between any two
  closed subsets of $\RP^{d-1}$. So we only need to consider the case
  where $d_\P(x,W) < \pi/2$.

  In this case, we let $W' = x^\perp \cap W$, and then let
  $V = W' \oplus x$. Let $z$ be the orthogonal projection of $x$ onto
  $W$, so that $d_\P(x, z) = d_\P(x, W)$. Let $\tilde{z}$ and
  $\tilde{x}$ be unit vector representatives of $z$ and $x$,
  respectively, chosen so that if
  \[
    \lambda = \pair{\tilde{x}}{\tilde{z}},
  \]
  then
  \[
    d_\P(x,z) = \cos^{-1}(\lambda).
  \]
  
  Let $v \in V - \{0\}$. We want to show that
  $d_\P([v], W) \le \cos^{-1}(\lambda)$, i.e. that for some $w \in W$,
  \[
    \frac{\pair{v}{w}}{||v|| \cdot ||w||} \ge \lambda.
  \]

  If $v \in W$, then we can choose $w = v$. Otherwise, we can rescale
  $v$ in order to write it as $w' + \tilde{x}$, for $w' \in W'$. Then
  let $w = w' + \tilde{z}$. Note that
  \[
    ||w|| = ||v|| = \sqrt{1 + ||w'||^2}.
  \]
  Now we just compute:
  \begin{align*}
    \frac{\pair{v}{w}}{||v|| \cdot ||w||}
    &= \frac{\pair{w' + \tilde{x}}{w'+ \tilde{z}}}{||v|| \cdot ||w||} = \frac{\pair{\tilde{x}}{\tilde{z}} + \pair{w'}{w'}}{1
      + ||w'||^2}\\
    &\ge \frac{\pair{\tilde{x}}{\tilde{z}} + \pair{\tilde{x}}{\tilde{z}}||w'||^2}{1 + ||w'||^2} = \pair{\tilde{x}}{\tilde{z}} = \lambda.
  \end{align*}
\end{proof}

Most of the work of proving the implication (\ref{item:main_thm_2})
$\implies$ (\ref{item:main_thm_1}) in Theorem
\ref{thm:convex_cocompact_equals_expansion} is contained in the
following:
\begin{prop}
  \label{prop:expansion_implies_cocpct_action}
  Let $\Omega$ be a convex domain preserved by a group
  $\Gamma \subset \PGL(d, \R)$. Let $C$ be a $\Gamma$-invariant subset
  of $\Omega$, closed in $\Omega$, with ideal boundary
  $\dee_iC$. Suppose that $\Gamma$ is expanding at the faces of
  $\dee_i C$, with respect to the metrics on $\Gr(k,d)$ specified in
  Lemma \ref{lem:subspace_achieving_dist_exists}.

  If either
  \begin{enumerate}[label=(\roman*)]
  \item $\Gamma$ is discrete and $\Omega$ is properly convex, or
  \item $\Gamma$ is uniformly expanding at the faces of $\dee_iC$,
  \end{enumerate}
  then $\Gamma$ acts cocompactly on $C$.
\end{prop}

\begin{proof}
  Danciger-Guéritaud-Kassel give a proof of this fact in the case
  where $\dee_iC$ contains no segments (see \cite[Lemma
  8.7]{dgk2017convex}). Their proof is based on work of Kapovich,
  Leeb, and Porti \cite{klpdomains}, which was in turn inspired by
  Sullivan \cite{sullivan1985quasiconformal}. Our proof will use
  essentially the same idea.

  We let $d_\P$ denote the angle metric on projective space, and we
  let $d_H$ denote the metric on $\Gr(k,d)$ induced by Hausdorff
  distance.

  For any $\eps > 0$, the set
  \[
    S_\eps = \{x \in C: d_\P(x, \dee \Omega) \ge \eps\}
  \]
  is compact. So, supposing for a contradiction that the action of
  $\Gamma$ on $C$ is not cocompact, for a sequence $\eps_n \to 0$,
  there exists $x_n$ so that $\Gamma \cdot x_n$ lies in
  $C - S_{\eps_n}$.

  We start by fixing a constant $E \ge 1$. If $\Gamma$ is discrete and
  $\Omega$ is properly convex, then we set $E = 1$; otherwise, we let
  $E$ be less than the uniform expansion constant. In either case, we
  can replace each $x_n$ with an element in its orbit so that
  \begin{equation}
    \label{eq:maximize_boundary_dist}
    d_\P(\gamma x_n, \dee \Omega) \le E \cdot d_\P(x_n, \dee \Omega)
  \end{equation}
  for all $\gamma \in \Gamma$. This is possible if $\Gamma$ is
  discrete and $\Omega$ is properly convex because then
  $\Gamma \cdot x_n$ is a discrete subset of $\Omega$. Otherwise, we
  choose $x_n$ sufficiently close to a point realizing the maximum
  distance between $\Gamma \cdot x_n$ and $\dee \Omega$.

  Up to a subsequence, $x_n$ converges in $\RP^{d-1}$ to some
  $x \in \dee_i C$. Let $F$ be the face of $\dee \Omega$ at $x$, and
  let $V \in \Gr(k,d)$ be the support of $F$.

  Let $U \subset \Gr(k, d)$ be an expanding neighborhood of $V$ in
  $\Gr(k,d)$, with expanding element $\gamma \in \Gamma$ expanding by
  a constant $E(\gamma) > E$ on $U$.

  Since $\dee \Omega$ is compact and $\Gamma$-invariant, there is some
  $z_n \in \dee \Omega$ so that
  \[
    d_\P(\gamma x_n, \gamma z_n) = d_\P(\gamma x_n, \dee \Omega).
  \]

  Since $x_n \to x$, and the distance from $\gamma x_n$ to
  $\gamma z_n$ is at most $\eps_n$, $z_n$ converges to $x$ as well.

  Proposition \ref{prop:supporting_hyperplanes_exist} implies that
  there is \emph{some} supporting hyperplane of $\Omega$ which
  intersects $z_n$. Any such sequence of supporting hyperplanes must
  sub-converge to a supporting hyperplane of $\Omega$ at $x$. This
  supporting hyperplane contains $V$, so there is a sequence
  $V_n \in \Gr(k,d)$ supporting $\Omega$ at $z_n$, which sub-converges
  to $V$.

  Since we know $\gamma z_n$ realizes the distance from $\gamma x_n$
  to $\dee \Omega$, we must have
  \begin{equation}
    \label{eq:expansion_inequality_1}
    d_\P(\gamma x_n, \dee \Omega) \ge d_\P(\gamma x_n, \gamma V_n).
  \end{equation}

  Then, Lemma \ref{lem:subspace_achieving_dist_exists} implies that we
  can choose subspaces $W_n \in \Gr(k,d)$ containing $x_n$ so that
  \begin{equation}
    \label{eq:expansion_inequality_2}
    d_\P(\gamma x_n, \gamma V_n) = d_H(\gamma W_n, \gamma V_n).
  \end{equation}

  Since $d_\P(\gamma x_n, \gamma V_n)$ converges to $0$,
  $d_H(\gamma W_n, \gamma V_n)$ does as well. Since $\gamma$ is fixed,
  and $V_n$ converges to $V$, $W_n$ also converges to $V$. So
  eventually, both $V_n$ and $W_n$ lie in the $E(\gamma)$-expanding
  neighborhood $U$ of $V$, meaning that we have
  \begin{equation}
    \label{eq:expansion_inequality_3}
    d_H(\gamma W_n, \gamma V_n) > E \cdot d_H(W_n, V_n).
  \end{equation}
  
  The trivial bound on Hausdorff distance implies that
  \begin{equation}
    \label{eq:expansion_inequality_4}
    d_H(W_n, V_n) \ge d_\P(x_n, V_n).
  \end{equation}
  Since $x_n \in \Omega$ and $V_n \subset \RP^{d-1} - \Omega$, any
  $d_\P$-geodesic from $x_n$ to $V_n$ must intersect $\dee \Omega$.
  This implies
  \begin{equation}
    \label{eq:expansion_inequality_5}
    d_\P(x_n, V_n) \ge d_\P(x_n, \dee \Omega).
  \end{equation}

  Putting (\ref{eq:expansion_inequality_1}),
  (\ref{eq:expansion_inequality_2}),
  (\ref{eq:expansion_inequality_3}),
  (\ref{eq:expansion_inequality_4}), and
  (\ref{eq:expansion_inequality_5}) together, we see that
  \[
    d_\P(\gamma x_n, \dee \Omega) > E \cdot d_\P(x_n, \dee
    \Omega),
  \]
  which contradicts (\ref{eq:maximize_boundary_dist}) above.
\end{proof}

We need one more lemma before we can show the main result of this
section. The statement is closely related to \cite[Lemma
6.3]{dgk2017convex}, and gives a condition for when a
$\Gamma$-invariant convex subset of a properly convex domain $\Omega$
contains $\Cor(\Gamma)$. (The result in \cite{dgk2017convex} is stated
for a cocompact action of a group $\Gamma$ on a convex set $C$, but
the proof only uses $\Gamma$-invariance.)
\begin{lem}
  \label{lem:gamma_invariant_contains_full_orbital_limit_set}
  Let $C$ be a nonempty convex set in $\Omega$ whose ideal boundary
  contains all of its faces, and suppose that
  $\Gamma \subseteq \Aut(\Omega)$ preserves $C$. Then $\dee_i C$
  contains $\Lambda_\Omega(\Gamma)$, the full orbital limit set of
  $\Gamma$.

  In particular, if $\Gamma$ is discrete, and the $\Gamma$ action on
  $C$ is cocompact, then the action of $\Gamma$ on $\Omega$ is convex
  cocompact and $\dee_i C = \Lambda_\Omega(\Gamma)$.
\end{lem}
\begin{proof}
  We follow the proof of Lemma 6.3 in \cite{dgk2017convex}.
  
  Let $z_\infty \in \Lambda_\Omega(\Gamma)$, which is by definition
  the limit of a sequence $\gamma_n z$ for some $z \in \Omega$ and a
  sequence $\gamma_n \in \Gamma$. Fix $y \in C$, and consider the
  sequence $\gamma_ny$. Since $d(\gamma_nz, \gamma_ny) = d(z,y)$ for
  all $n$, Proposition \ref{prop:strata_hilbert_bounded} implies that
  up to a subsequence, $\gamma_nz$ and $\gamma_ny$ both converge to
  points in the same face of $\dee \Omega$. But any accumulation point
  of $\gamma_ny$ in $\dee \Omega$ lies in $\dee_i C$ and $\dee_i C$
  contains its faces, so $z_\infty \in \dee_i C$.

  Since $\dee_iC$ contains $\Lambda_\Omega(\Gamma)$, $C$ must contain
  $\Cor(\Gamma)$. \cite[Lemma 4.10 (3)]{dgk2017convex} then implies
  that $\Lambda_\Omega(\Gamma) = \dee_iC$ is closed in $C$, which
  means that $\Cor(\Gamma)$ is closed in $C$ and the action on
  $\Cor(\Gamma)$ is cocompact.
\end{proof}

\begin{proof}[Proof of (\ref{item:main_thm_2})
  $\implies$ (\ref{item:main_thm_1}) in Theorem
  \ref{thm:convex_cocompact_equals_expansion}]
  Let $\Omega$ be a properly convex domain, let $\Gamma$ be a discrete
  subgroup of $\Aut(\Omega)$, and $\Lambda$ be a $\Gamma$-invariant,
  closed and boundary-convex subset of $\dee \Omega$ with nonempty
  convex hull, such that $\Lambda$ contains all of its faces and
  $\Gamma$ is uniformly expanding at the faces of $\Lambda$.

  Since $\Lambda$ is boundary-convex and has nonempty convex hull,
  $\Lambda$ is exactly the ideal boundary of $\CH(\Lambda)$. So,
  Proposition \ref{prop:expansion_implies_cocpct_action} implies that
  $\Gamma$ acts cocompactly on $\CH(\Lambda)$. Since $\Lambda$ also
  contains its faces, applying Lemma
  \ref{lem:gamma_invariant_contains_full_orbital_limit_set} with
  $C = \CH(\Lambda)$ completes the proof.
\end{proof}



\section{Actions on spaces of projective domains}
\label{sec:spaces_of_domains}

In this section we recall the statement of Benz\'ecri's cocompactness
theorem for convex projective domains, as well as prove a version of
it (Proposition \ref{prop:relative_benzecri}) that applies relative to
a direct sum decomposition of $\Rd$.

\subsection{The space of projective domains}
Good references for this material include \cite{goldman88projective}
and \cite{marquis2013around}.

Let $V$ be a real vector space. We denote the set of non-empty
properly convex open subsets of $\P(V)$ by $\mc{C}(V)$. We topologize
$\mc{C}(V)$ via the metric:
\[
  d(\Omega_1, \Omega_2) := d_{\mr{Haus}}(\overline{\Omega_1},
  \overline{\Omega_2}),
\]
where $d_{\mr{Haus}}(\cdot, \cdot)$ is the Hausdorff distance
induced by any metric on $\P(V)$ (the choice of metric on $\P(V)$
does not affect the topology on $\mc{C}(V)$).

\begin{defn}
  A \emph{pointed properly convex domain} in $\P(V)$ is a pair
  $(\Omega, x)$, where $\Omega \in \mc{C}(V)$ and $x \in \Omega$. We
  denote the set of pointed properly convex domains in $\P(V)$ by
  $\mc{C}_*(V)$, and topologize $\mc{C}_*(V)$ by viewing it as a
  subspace of $\mc{C}(V) \times \P(V)$.
\end{defn}

$\PGL(V)$ acts on both $\mc{C}(V)$ and $\mc{C}_*(V)$ by
homeomorphisms. We have the following important result:
\begin{thm}[Benz\'ecri, \cite{benzecri1960varietes}]
  \label{thm:benzecri_cocompactness}
  The action of $\PGL(V)$ on $\mc{C}_*(V)$ is proper and cocompact.
\end{thm}

\subsection{Benz\'ecri relative to a direct sum}

We now let $V_a$, $V_b$ be subspaces of $V$ so that
$V_a \oplus V_b = V$. The decomposition induces natural projection
maps $\pi_{V_a}:V \to V_a$ and $\pi_{V_b}:V \to V_b$, as well as a
decomposition of the dual $V^*$ into $V_a^* \oplus V_b^*$. Here, and
throughout this section, we will identify $V_a^*$, $V_b^*$ with the
linear functionals on $V$ which vanish on $V_b, V_a$.

When $\Omega$ is a convex subset of $\P(V)$ which is disjoint from
$\P(V_b)$, we let $\pi_{V_a}(\Omega)$ be the projectivization of
$\pi_{V_a}(\tilde{\Omega})$, where $\tilde{\Omega}$ is a cone over
$\Omega$. A priori this is only a convex subset of $\P(V_a)$, although
we will see (Proposition \ref{prop:projection_sharp}) that if $\Omega$
is properly convex and open, and $\overline{\Omega}$ is disjoint from
$\P(V_b)$, then $\pi_{V_a}(\Omega)$ is properly convex and open in
$\P(V_a)$.

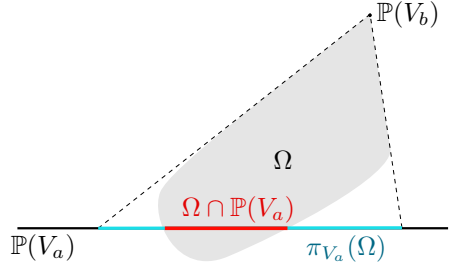
\begin{figure}[h]
  \centering
  \def\svgwidth{2.5in}
\begingroup%
  \makeatletter%
  \providecommand\color[2][]{%
    \errmessage{(Inkscape) Color is used for the text in Inkscape, but the package 'color.sty' is not loaded}%
    \renewcommand\color[2][]{}%
  }%
  \providecommand\transparent[1]{%
    \errmessage{(Inkscape) Transparency is used (non-zero) for the text in Inkscape, but the package 'transparent.sty' is not loaded}%
    \renewcommand\transparent[1]{}%
  }%
  \providecommand\rotatebox[2]{#2}%
  \newcommand*\fsize{\dimexpr\f@size pt\relax}%
  \newcommand*\lineheight[1]{\fontsize{\fsize}{#1\fsize}\selectfont}%
  \ifx\svgwidth\undefined%
    \setlength{\unitlength}{839.78518947bp}%
    \ifx\svgscale\undefined%
      \relax%
    \else%
      \setlength{\unitlength}{\unitlength * \real{\svgscale}}%
    \fi%
  \else%
    \setlength{\unitlength}{\svgwidth}%
  \fi%
  \global\let\svgwidth\undefined%
  \global\let\svgscale\undefined%
  \makeatother%
  \begin{picture}(1,0.53093062)%
    \lineheight{1}%
    \setlength\tabcolsep{0pt}%
    \put(0,0){\includegraphics[width=\unitlength,page=1]{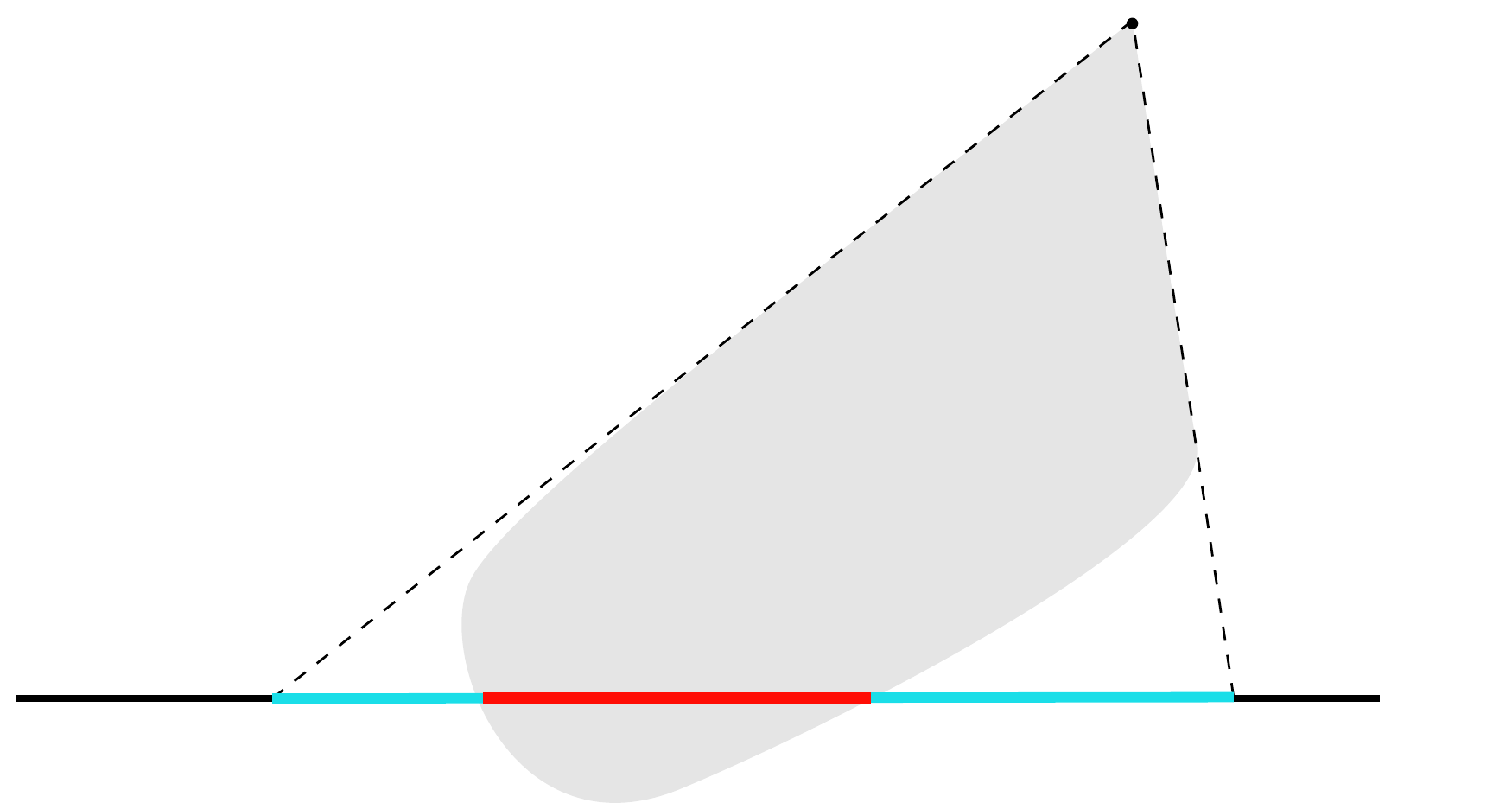}}%
    \put(0.35513043,0.09513217){\color[rgb]{0.88627451,0.01176471,0.01176471}\makebox(0,0)[lt]{\lineheight{1.25}\smash{\begin{tabular}[t]{l}$\Omega \cap \P(V_a)$\end{tabular}}}}%
    \put(0.76371481,0.50378919){\color[rgb]{0,0,0}\makebox(0,0)[lt]{\lineheight{1.25}\smash{\begin{tabular}[t]{l}$\P(V_b)$\end{tabular}}}}%
    \put(-0.00296532,0.01480999){\color[rgb]{0,0,0}\makebox(0,0)[lt]{\lineheight{1.25}\smash{\begin{tabular}[t]{l}$\P(V_a)$\end{tabular}}}}%
    \put(0.61581626,0.01332675){\color[rgb]{0,0.40784314,0.52156863}\makebox(0,0)[lt]{\lineheight{1.25}\smash{\begin{tabular}[t]{l}$\pi_{V_a}(\Omega)$\end{tabular}}}}%
    \put(0.54687188,0.19375711){\color[rgb]{0,0,0}\makebox(0,0)[lt]{\lineheight{1.25}\smash{\begin{tabular}[t]{l}$\Omega$\end{tabular}}}}%
  \end{picture}%
\endgroup%

  \caption{The domains $\Omega \cap \P(V_a)$ and
    $\pi_{V_a}(\Omega)$. In this case, $\pi_{V_a}(\Omega)$ is properly
    convex even though $\overline{\Omega}$ intersects $\P(V_b)$.}
  \label{fig:convex_intersect_projection}
\end{figure}

We remark that if $\Omega \cap \P(V_b)$ is nonempty, then
$\pi_{V_a}(\Omega)$ is not even well-defined. On the other hand, if
$\overline{\Omega} \cap \P(V_b)$ is nonempty, but
$\Omega \cap \P(V_b)$ is empty, then $\pi_{V_a}(\Omega)$ does exist,
and may or may not be a properly convex subset of $\P(V_a)$.

\begin{defn}
  \label{defn:domains_with_compact_slice}
  Let $V = V_a \oplus V_b$, and let $\mc{K}_a$ be a subset of
  $\mc{C}_*(V_a)$. We define the subset $\mc{C}_*(V_a, V_b, \mc{K}_a)$
  by
  \[
    \mc{C}_*(V_a, V_b, \mc{K}_a) := \left\{
      (\Omega, x) \in \mc{C}_*(V) :
      \begin{array}{l}
        \P(V_b) \cap \Omega = \emptyset,\\
        (\Omega \cap \P(V_a), x) \in \mc{K}_a,\\
        (\pi_{V_a}(\Omega), x) \in \mc{K}_a
      \end{array}
    \right\}.
  \]
\end{defn}

The groups $\GL(V_a)$ and $\GL(V_b)$ both have a well-defined action
on $\mc{C}_*(V)$: we take $g \in \GL(V_a)$ and $h \in \GL(V_b)$ to act
by the projectivizations of $g \oplus \mr{id}_{V_b}$,
$\mr{id}_{V_a} \oplus h$ respectively, on $\P(V_a \oplus V_b)$.

Since the $\GL(V_b)$-action on $\P(V)$ fixes $\P(V_a)$ pointwise and
commutes with projection to $\P(V_a)$, for any
$\mc{K}_a \subset \mc{C}_*(V_a)$, $\GL(V_b)$ acts on the subset
$\mc{C}_*(V_a, V_b, \mc{K}_a)$. The main result of this section is the
following:
\begin{prop}
  \label{prop:relative_benzecri}
  Let $V_a$, $V_b$ be subspaces of a real vector space $V$ such that
  $V_a \oplus V_b = V$. For any compact subset
  $\mc{K}_a \subset \mc{C}_*(V_a)$, the action of $\GL(V_b)$ on
  $\mc{C}_*(V_a, V_b, \mc{K}_a)$ is proper and cocompact.
\end{prop}

\subsection{Convex cones in direct sums}

Before proving Proposition \ref{prop:relative_benzecri}, we explore
some of the properties of convex cones in a vector space $V$ which
splits as a direct sum $V = V_a \oplus V_b$.

\subsubsection{Duality}

If $V = V_a \oplus V_b$, a convex cone $C \subset V$ determines two
different convex cones in $V_a$, namely $C \cap V_a$ and
$\pi_{V_a}(C)$. The arguments in this section rely heavily on the fact
that these projection and intersection operations are in some sense
``dual'' to each other. Before explaining how this works, we first
prove a lemma:
\begin{lem}
  \label{lem:orthogonal_proj_dual}
  Let $V$ be a real vector space with $V = V_a \oplus V_b$, and let
  $C$ be a convex cone intersecting $V_b$ trivially. Then
  $C^* \cap V_a^* \subseteq \pi_{V_a}(C)^* \cap V_a^*$ and
  \[
    \overline{\pi_{V_a}(C)^*} \cap V_a^* \subseteq \overline{C^*} \cap
    V_a^*.
  \]
  Moreover, if $\overline{C} \cap V_b = \{0\}$ then in fact
  \[
    \pi_{V_a}(C)^* \cap V_a^* = C^* \cap V_a^*.
  \]
\end{lem}
\begin{proof}
  First let $\alpha \in C^* \cap V_a^*$. Let $v$ be any nonzero
  element of the closure of $\pi_{V_a}(C)$, so that
  $v + v_2 \in \overline{C}$ for some $v_2 \in V_b$. We know that
  $\alpha(v + v_2) \ne 0$ and $\alpha(v_2) = 0$, so $\alpha(v) \ne
  0$. This shows that $\alpha$ is in $\pi_{V_a}(C)^*$.

  Now let $\alpha \in \overline{\pi_{V_a}(C)^*} \cap V_a^* - \{0\}$,
  and let $v \in C$. We can write $v = v_1 + v_2$ for $v_1 \in V_a$,
  $v_2 \in V_b$; since we assume $C$ does not intersect $V_b$, and
  $C \subset V - \{0\}$ by definition, $v_1$ is nonzero. Then since
  $\alpha \in V_a^*$, $\alpha(v) = \alpha(v_1) \ne 0$. So,
  $\alpha \in \overline{C^*}$.

  If we further assume that $\overline{C} \cap V_b = \{0\}$, a similar
  argument shows that any $\alpha \in \pi_{V_a}(C)^* \cap V_a^*$ is
  nonzero on any $v \in \overline{C} - \{0\}$, implying
  $\alpha \in C^*$.
\end{proof}

Now, suppose that $C_a$ is a convex cone in $V - \{0\}$, for
$V = V_a \oplus V_b$. The intersection $C_a^* \cap V_a^*$ consists of
functionals in $C_a^*$ which vanish on $V_b$. If we know that $C_a$
lies inside of $V_a$, then any functional on $V_a$ which does not
vanish anywhere on $\overline{C_a} - \{0\}$ can be extended by zero on
$V_b$ to get an element of $C_a^* \cap V_a^*$. So in this case,
$C_a^* \cap V_a^*$ is canonically identified with the dual of the
cone $C_a$ \emph{viewed as a cone in $V_a$}.

We can say this a different way via the following:
\begin{defn}
  \label{defn:restricted_dual}
  For each subspace $U \subset V$, we define a ``restricted duality''
  operation $D_U$, which takes convex cones in $U$ to convex cones in
  $U^*$ via the dual operation on $U$. Explicitly, if $C \subset U$ is
  a convex cone, we let
  \[
    D_U(C) = \{\alpha \in U^* : \alpha(v) > 0 \textrm{ for all } v \in
    \overline{C} - \{0\}\}.
  \]
\end{defn}
By definition, we have $D_V(C) = C^*$ for any convex cone
$C \subset V$.  The reasoning above tells us that when
$V = V_a \oplus V_b$, then for any cone $C_a \subset V_a$, we have
$D_{V_a}(C_a) = C_a^* \cap V_a^*$.

With this notation, Lemma~\ref{lem:orthogonal_proj_dual} can be
restated as:
\begin{lem}
  \label{lem:restricted_dual_intersect}
  Let $V$ be a real vector space with $V = V_a \oplus V_b$, and let
  $C$ be a convex cone intersecting $V_b$ trivially. Then
  $D_V(C) \cap V_a^* \subseteq D_{V_a}(\pi_{V_a}(C))$ and
  $\overline{D_{V_a}(\pi_{V_a}(C))} \subseteq \overline{D_V(C) \cap
    V_a^*}$. Moreover, if $\overline{C} \cap V_b = \{0\}$, then in
  fact
  \[
    D_V(C) \cap V_a^* = D_{V_a}(\pi_{V_a}(C)).
  \]
\end{lem}

As a consequence of this lemma, we note:
\begin{prop}
  \label{prop:projection_sharp}
  Let $C$ be a sharp (Definition \ref{defn:convex_cone}) open convex
  cone in a vector space $V = V_a \oplus V_b$. If
  $\overline{C} - \{0\}$ intersects $V_b$ trivially, then the
  projection $\pi_{V_a}(C)$ is sharp and open in $V_a$.
\end{prop}
\begin{proof}
  Openness is immediate since projection is an open map. Since $C$ is
  sharp, if $\overline{C}$ does not intersect $V_b$, then there is
  some $\alpha \in V^*$ whose kernel contains $V_b$ and does not
  intersect $\overline{C}$, i.e. $\alpha \in C^* \cap V_a^*$. Since
  non-intersection with $\overline{C}$ is an open condition,
  $C^* \cap V_a^*$ is a nonempty open subset of $V_a^*$. Then Lemma
  \ref{lem:orthogonal_proj_dual} implies that
  $\pi_{V_a}(C)^* \cap V_a^*$ is nonempty and open in $V_a^*$. So its
  dual in $V_a^{**} = V_a$ is sharp by part \ref{item:sharp_open} of
  Proposition~\ref{prop:duality_properties}.
\end{proof}

For the rest of the section we will be working with convex domains in
$\P(V)$, rather than convex cones in $V$. The restricted dual
operation $D_U$ from Definition~\ref{defn:restricted_dual} gives rise
to a restricted dual operation on convex domains contained in
projective subspaces $\P(U) \subseteq \P(V)$; we also denote this by
$D_U$. Also, recall (from Remark~\ref{rem:duality_convention}) that if
$W$ is an element in some dual domain $D_V(\Omega) = \Omega^*$, we
identify $W$ with a projective hyperplane in $\P(V)$.

\subsubsection{Convex hulls}

If $\Omega_1$, $\Omega_2$ are properly convex subsets of $\P(V)$, we
cannot always find a minimal properly convex subset
$\Omega \subset \P(V)$ which contains $\Omega_1 \cup \Omega_2$ (that
is, convex hulls do not always exist). Here we describe some
circumstances under which this is possible.

\begin{defn}
  \label{defn:convex_hull_hyperplane}
  Let $\Omega_1$, $\Omega_2$ be properly convex sets in $\P(V)$. For
  each $W \in \Omega_1^* \cap \Omega_2^*$, we let
  $\mr{Hull}_W(\Omega_1, \Omega_2)$ denote the convex hull of
  $\Omega_1$ and $\Omega_2$ in the affine chart $\P(V) - W$.
\end{defn}

The set $\mr{Hull}_W(\Omega_1, \Omega_2)$ is minimal among all convex
subsets of $\P(V) - W$ containing $\Omega_1 \cup \Omega_2$. However,
it is possible that for some other
$W' \in \Omega_1^* \cap \Omega_2^*$,
$\mr{Hull}_{W'}(\Omega_1, \Omega_2)$ is not contained in $\P(V) -
W$. So, to guarantee minimality among all convex subsets of $\P(V)$,
we need a little more:
\begin{lem}
  \label{lem:convex_hulls_exist}
  If $\Omega_1 \cap \Omega_2$ is nonempty, then for any
  $W \in \Omega_1^* \cap \Omega_2^*$,
  $\mr{Hull}_W(\Omega_1, \Omega_2)$ is the unique minimal properly
  convex subset of $\P(V)$ containing $\Omega_1 \cup \Omega_2$.
\end{lem}
\begin{proof}
  Let $A$ be the affine chart $\P(V) - W$, and let $H$ be any properly
  convex set containing $\Omega_1 \cup \Omega_2$. Since
  $\Omega_1 \cap \Omega_2$ is nonempty, $\Omega_1 \cup \Omega_2$ is a
  connected subset of $A$, so it is contained in a single connected
  component $C$ of $H \cap A$. This component is a convex subset of
  $A$, so by definition $C$ (hence $H$) contains
  $\mr{Hull}_W(\Omega_1, \Omega_2)$.
\end{proof}

Lemma \ref{lem:convex_hulls_exist} allows us to define the convex hull
of a pair of properly convex sets without reference to a particular
affine chart.
\begin{defn}
  \label{defn:convex_hull_unique}
  When $\Omega_1$, $\Omega_2$ are properly convex sets such that
  $\Omega_1 \cap \Omega_2$ and $\Omega_1^* \cap \Omega_2^*$ are both
  nonempty, we let $\mr{Hull}(\Omega_1, \Omega_2)$ denote the minimal
  properly convex set containing $\Omega_1 \cup \Omega_2$.
\end{defn}

\subsection{Proving Benz\'ecri for direct sums}

We can now begin proving Proposition \ref{prop:relative_benzecri}. As
a first step, we consider the case where $\dim V_a = 1$,
i.e. $\P(V_a)$ is identified with a single point in $\P(V)$.
\begin{lem}
  \label{lem:benzecri_1_dim}
  Let $V = V_b \oplus x$ for a point $x \in \P(V)$. Then $\GL(V_b)$
  acts properly and cocompactly on the set of domains
  \[
    \mc{C}_*(x, V_b) := \mc{C}_*(x, V_b, \mc{C}_*(x)) = \{(\Omega, x)
    \in \mc{C}_*(V) : \P(V_b) \cap \Omega = \emptyset \}.
  \]
\end{lem}
Here $\mc{C}_*(x)$ denotes the space of pointed nonempty properly
convex domains in $\P(x) \simeq \RP^0$, so the only nonempty domain in
this space is the singleton $\{x\}$.

Statements similar to this lemma can be found in work of Frankel (see
\cite[Theorem 9.3]{frankel89}) and Benoist (section 2.3 in
\cite{benoist2003convexes}); Benoist notes that the idea already
appears in Benz\'ecri \cite{benzecri1960varietes}.

\begin{proof}
  Properness follows immediately from the standard Benz\'ecri theorem
  (Theorem \ref{thm:benzecri_cocompactness}), since the restriction of
  a proper action of a group $G$ on $X$ to a closed subgroup $H$ and
  an $H$-invariant subset of $X$ is always proper. So, we focus on
  cocompactness.
  
  Let $(\Omega_n, x)$ be a sequence of domains in $\mc{C}_*(x,
  V_b)$. Theorem \ref{thm:benzecri_cocompactness} implies that we can
  find group elements $g_n \in \PGL(V)$ so that the sequence of
  pointed domains
  \[
    (g_n\Omega_n, g_nx)
  \]
  sub-converges to a pointed domain $(\Omega, x')$. We want to show
  that these group elements can be chosen to preserve the
  decomposition $V_b \oplus x$.
  
  We know that $V_b$ lies in $\overline{\Omega_n^*}$, so $g_nV_b$ lies
  in $g_n\overline{\Omega_n^*}$ for all $n$, and a subsequence of
  $g_nV_b$ converges to some $W \in \overline{\Omega^*}$. In
  particular, $W$ does not contain $x'$. This means that we can find a
  sequence of group elements $g_n'$, lying in a fixed compact subset
  of $\PGL(V)$, so that
  \[
    g_n' \cdot g_nV_b = V_b, \quad g_n' \cdot g_nx = x.
  \]
  Since the $g_n'$ lie in a compact subset of $\PGL(V_b)$, the domains
  \[
    g_n'g_n\Omega_n
  \]
  must also sub-converge to some properly convex domain $\Omega'$,
  which contains $x$. So we can replace $g_n$ with $g_n'g_n$ to get
  the desired sequence of group elements.
\end{proof}

Lemma \ref{lem:benzecri_1_dim} gets us partway to proving Proposition
\ref{prop:relative_benzecri}. We see that if $\Omega$ is any domain in
$\mc{C}_*(V_a, V_b, \mc{K}_a)$, we can always find some
$h \in \GL(V_b)$ so that $h \Omega \cap \P(V_b \oplus x)$ lies in a
fixed compact set of domains in $\mc{C}(V_b \oplus x)$. This is almost
enough to ensure that $h\Omega$ itself lies in a fixed compact set of
domains in $\mc{C}(V)$. The exact condition we will need is the
following:
\begin{lem}
  \label{lem:projection_intersection}
  Let $V$ be a real vector space, and suppose
  $V = W_a \oplus V_b \oplus x$, for a point $x \in \P(V)$.

  Let $\Omega_a, \Omega_a'$ be properly convex domains in
  $\P(W_a \oplus x)$, and let $\Omega_b, \Omega_b'$ be properly convex
  domains in $\P(V_b \oplus x)$, such that
  \begin{align*}
    x \in \Omega_a \subset \Omega_a',\\
    x \in \Omega_b \subset \Omega_b'.
  \end{align*}
  
  There exist properly convex domains $\Omega_1 \subset \Omega_2$ in
  $\P(V)$, with $x \in \Omega_1$, such that any $\Omega \in \mc{C}(V)$
  disjoint from $\P(W_a)$ and $\P(V_b)$ which satisfies:
  \begin{enumerate}
  \item $\Omega_a' \supset \pi_{W_a \oplus x}(\Omega)$,
  \item $\Omega_a \subset \Omega \cap \P(W_a \oplus x)$,
  \item \label{item:replaceable_assumption}
    $\Omega_b' \supset \pi_{V_b \oplus x}(\Omega)$
  \item $\Omega_b \subset \Omega \cap \P(V_b \oplus x)$,
  \end{enumerate}
  also satisfies $\Omega_1 \subset \Omega \subset \Omega_2$.
\end{lem}

\begin{figure}[h]
  \centering \def\svgwidth{3in}
\begingroup%
  \makeatletter%
  \providecommand\color[2][]{%
    \errmessage{(Inkscape) Color is used for the text in Inkscape, but the package 'color.sty' is not loaded}%
    \renewcommand\color[2][]{}%
  }%
  \providecommand\transparent[1]{%
    \errmessage{(Inkscape) Transparency is used (non-zero) for the text in Inkscape, but the package 'transparent.sty' is not loaded}%
    \renewcommand\transparent[1]{}%
  }%
  \providecommand\rotatebox[2]{#2}%
  \newcommand*\fsize{\dimexpr\f@size pt\relax}%
  \newcommand*\lineheight[1]{\fontsize{\fsize}{#1\fsize}\selectfont}%
  \ifx\svgwidth\undefined%
    \setlength{\unitlength}{589.73802776bp}%
    \ifx\svgscale\undefined%
      \relax%
    \else%
      \setlength{\unitlength}{\unitlength * \real{\svgscale}}%
    \fi%
  \else%
    \setlength{\unitlength}{\svgwidth}%
  \fi%
  \global\let\svgwidth\undefined%
  \global\let\svgscale\undefined%
  \makeatother%
  \begin{picture}(1,0.87750825)%
    \lineheight{1}%
    \setlength\tabcolsep{0pt}%
    \put(0,0){\includegraphics[width=\unitlength,page=1]{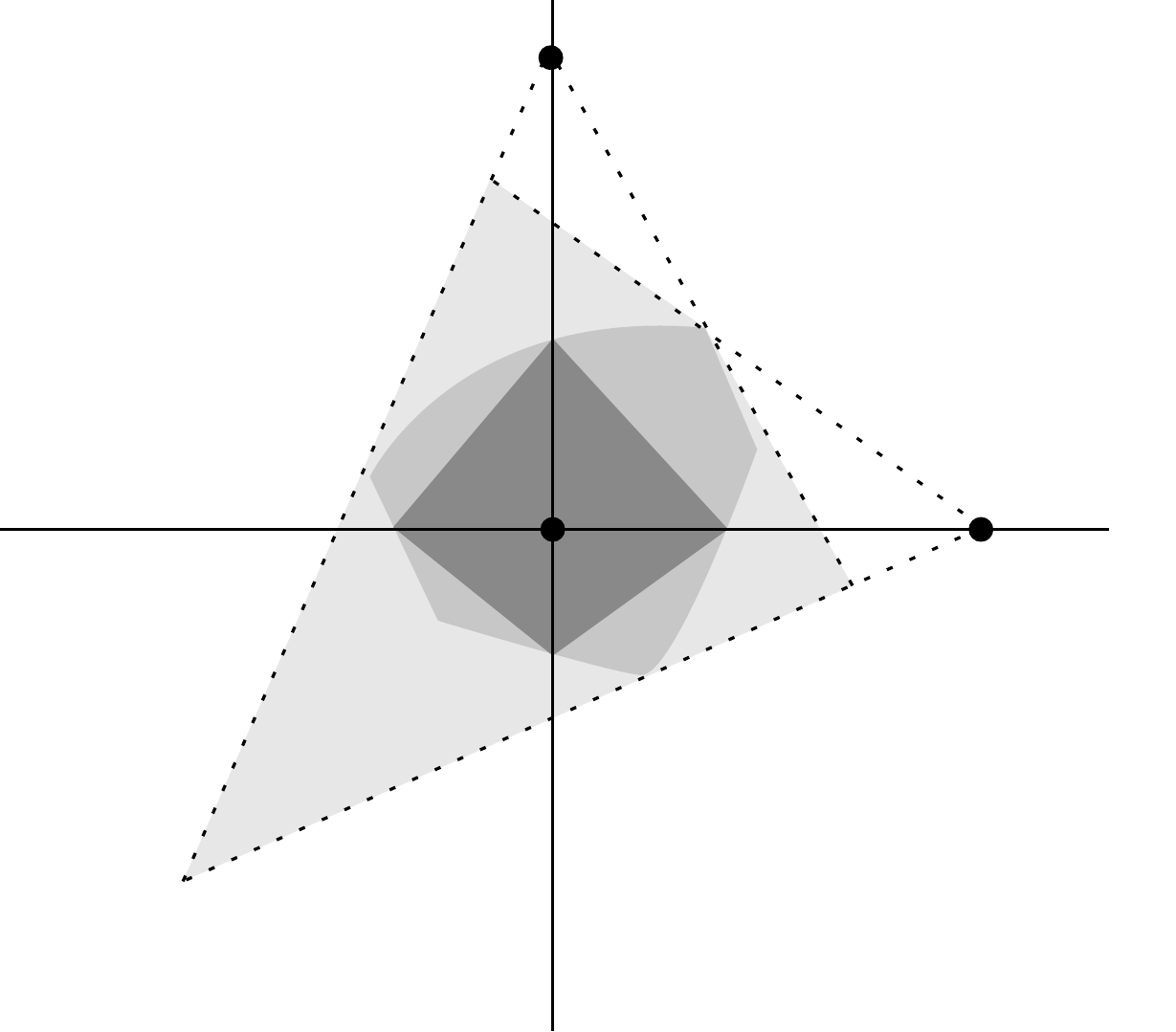}}%
    \put(0.5332346,0.535056){\color[rgb]{0,0,0}\makebox(0,0)[lt]{\lineheight{1.25}\smash{\begin{tabular}[t]{l}$\Omega$\end{tabular}}}}%
    \put(0.49222876,0.45936437){\color[rgb]{0,0,0}\makebox(0,0)[lt]{\lineheight{1.25}\smash{\begin{tabular}[t]{l}$\Omega_1$\end{tabular}}}}%
    \put(0.25831684,0.25114559){\color[rgb]{0,0,0}\makebox(0,0)[lt]{\lineheight{1.25}\smash{\begin{tabular}[t]{l}$\Omega_2$\end{tabular}}}}%
    \put(0.48056346,0.37447317){\color[rgb]{0,0,0}\makebox(0,0)[lt]{\lineheight{1.25}\smash{\begin{tabular}[t]{l}$x$\end{tabular}}}}%
    \put(0.49983342,0.83438241){\color[rgb]{0,0,0}\makebox(0,0)[lt]{\lineheight{1.25}\smash{\begin{tabular}[t]{l}$V_b$\end{tabular}}}}%
    \put(0.83256105,0.35520327){\color[rgb]{0,0,0}\makebox(0,0)[lt]{\lineheight{1.25}\smash{\begin{tabular}[t]{l}$W_a$\end{tabular}}}}%
  \end{picture}%
\endgroup%

  \caption{$\Omega$ fits between a pair of domains $\Omega_1$ and
    $\Omega_2$, which depend only on the intersections and
    projections between $\Omega$ and $\P(W_a \oplus x)$,
    $\P(V_b \oplus x)$.}
\end{figure}
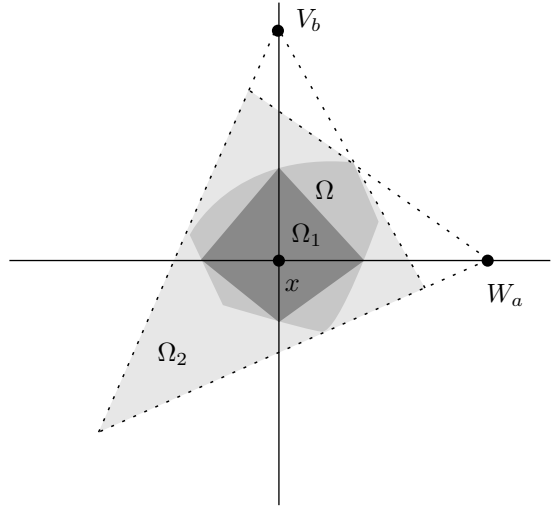

\begin{proof}
  We know $\Omega_a \cap \Omega_b = \{x\}$. If necessary, we can
  slightly shrink $\Omega_a$ and $\Omega_b$ so that
  $\overline{\Omega_a} \cap \P(W_a) = \emptyset$ and
  $\overline{\Omega_b} \cap \P(V_b) = \emptyset$, which means that
  $\P(W_a \oplus V_b)$ can be viewed as an element of
  $\Omega_a^* \cap \Omega_b^*$. So, the convex hull (Definition
  \ref{defn:convex_hull_unique}) $\mr{Hull}(\Omega_a, \Omega_b)$ of
  $\Omega_a, \Omega_b$ exists.

  In any affine chart $A$ containing $\Omega_a \cup \Omega_b$, the
  subspaces $W_a \oplus x$ and $V_b \oplus x$ correspond to pair of
  transverse affine subspaces intersecting at the point $x$, which
  together span all of $A$. We know $\Omega_a$ and $\Omega_b$ are open
  in these subspaces and each contain $x$, so the interior of their
  convex hull in $A$ also contains $x$. Thus, we may define $\Omega_1$
  to be the interior of $\mr{Hull}(\Omega_a, \Omega_b)$.

  To build $\Omega_2$, we consider the ``restricted dual'' domains
  \[
    D_{W_a \oplus x}(\Omega_a') \subset \P((W_a \oplus x)^*), \quad
    D_{V_b \oplus x}(\Omega_b') \subset \P((V_b \oplus x)^*),
  \]
  which are defined using the restricted dual operation $D_U$ from
  Definition~\ref{defn:restricted_dual}. To simplify notation, we
  write $D(\Omega_a') = D_{W_a \oplus x}(\Omega_a')$ and
  $D(\Omega_b') = D_{V_b \oplus x}(\Omega_b')$. Explicitly, we have
  \[
    D(\Omega_a') = (\Omega_a')^* \cap \P((W_a \oplus x)^*), \quad
    D(\Omega_b') = (\Omega_b')^* \cap \P((V_b \oplus x)^*).
  \]
  Since $\Omega_a'$ and $\Omega_b'$ are properly convex subsets of
  $\P(W_a \oplus x)$ and $\P(V_b \oplus x)$, $D(\Omega_a')$ and
  $D(\Omega_b')$ are open in $\P((W_a \oplus x)^*)$ and
  $\P((V_b \oplus x)^*)$ (see part \ref{item:sharp_open} of
  Proposition~\ref{prop:duality_properties}).

  We also know that $x$ lies in $D(\Omega_a')^* \cap
  D(\Omega_b')^*$. So we can define the convex open set $\Omega_2^*$
  to be the interior of
  \[
    \mr{Hull}_x(D(\Omega_a'), D(\Omega_b')),
  \]
  using Definition \ref{defn:convex_hull_hyperplane}. Using similar
  reasoning as for $\Omega_1$, we can see that the interior of this
  hull is nonempty, because $(W_a \oplus x)^*$ and $(V_b \oplus x)^*$
  span $V^*$ and $D(\Omega_a')$ and $D(\Omega_b')$ are open subsets of
  the corresponding projective subspaces.

  Let $\Omega$ be any domain satisfying the hypotheses of the
  lemma. Since duality reverses inclusions, we know
  $D(\Omega_a') \subseteq D_{W_a \oplus x}(\pi_{W_a \oplus
    x}(\Omega))$ and
  $D(\Omega_b') \subseteq D_{V_b \oplus x}(\pi_{V_b \oplus
    x}(\Omega))$. Then, Lemma~\ref{lem:restricted_dual_intersect}
  implies
  \begin{align*}
    &D(\Omega_a')  \subseteq \overline{D_V(\Omega)} \cap \P((W_a
      \oplus x)^*),\\
    &D(\Omega_b') \subseteq \overline{D_V(\Omega)} \cap \P((V_b \oplus x)^*).
  \end{align*}
  In particular, $D(\Omega_a')$ and $D(\Omega_b')$ are both contained
  in $\overline{D_V(\Omega)} = \overline{\Omega^*}$. Since
  $\Omega^{**} = \Omega$ contains $x$, $\overline{\Omega^*}$ is
  contained in the affine chart $\P(V^*) - x$. So,
  $\overline{\Omega^*}$ contains the closure of
  \[
    \mr{Hull}_x(D(\Omega_a'), D(\Omega_b')),
  \]
  meaning $\Omega^*$ contains $\Omega_2^*$ and $\Omega$ is contained
  in the properly convex set $\Omega_2 = \Omega_2^{**}$.
\end{proof}

\begin{remark}
  If $\overline{\Omega}_b$ does not intersect $\P(V_b)$ and
  $\overline{\Omega}_a$ does not intersect $\P(W_a)$, we can work in
  the affine chart $\P(V) - \P(W_a \oplus V_b)$, and Lemma
  \ref{lem:projection_intersection} is equivalent to the fact that if
  a convex subset $C$ of an affine space has open and bounded
  projections to and intersections with a pair of complementary affine
  subspaces, $C$ is itself open and bounded in terms of the size of
  the projections and intersections.

  We do not take this approach because we do \emph{not} want to assume
  that $\overline{\Omega}_b$ and $\P(V_b)$ are disjoint.
\end{remark}

Our next task is to show that we can sometimes replace assumption
(\ref{item:replaceable_assumption}) in Lemma
\ref{lem:projection_intersection} with:
\begin{enumerate}
\item[(3a)] $\Omega_b' \supset \Omega \cap \P(V_b \oplus x)$.
\end{enumerate}

This will be done in Proposition
\ref{prop:cone_dual_intersections_open} below. We start with some
Euclidean geometry.

We endow $\R^d$ with its standard inner product. For a subspace
$W \subseteq \R^d$, we let $\pi_W:\R^d \to W$ denote the orthogonal
projection, and for $R > 0$, let $B(R)$ denote the open ball around
the origin of radius $R$.
\begin{lem}\label{lem:convex_subset_intersect_project}
  Let $\Omega$ be a convex subset of $\R^d$ containing the origin, and
  let $W$ be a subspace of $\R^d$.

  Suppose that there are $R_1, R_2 > 0$ so that:
  \begin{itemize}
  \item $B(R_1) \cap W^\perp \subset \Omega \cap W^\perp$,
  \item $\pi_{W^\perp}(\Omega) \subset B(R_2)$.
  \end{itemize}
  Then there exists a linear map $f:\R^d \to \R^d$, depending only on
  $R_1$ and $R_2$, so that $\pi_{W}(\Omega) \subset f(\Omega \cap W)$.
\end{lem}
\begin{proof}
  Let $p$ be any point in $\pi_W(\Omega)$, and let $z$ be some point
  in $\Omega$ so that $\pi_W(z) = p$. We can write $z = p + y$ for
  $y \in \pi_{W^\perp}(\Omega)$.

  \begin{figure}[h]
    \centering
    \def\svgwidth{2.8in}
    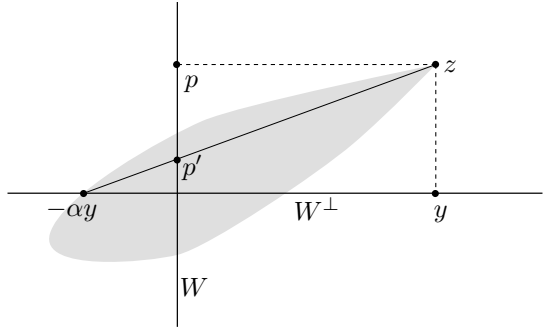
    \caption{Illustration for the proof of Lemma
      \ref{lem:convex_subset_intersect_project}. The ratio
      $||p|| / ||p'||$ is bounded in terms of $\alpha$.}
    \label{fig:euclidean_lemma_proof}
  \end{figure}

  Let $\ell$ be the line through the origin passing through $y$. For
  some $\alpha > 0$, we know that $\ell$ intersects
  $\overline{(\Omega \cap W^\perp)} - B(R_1)$ at $y' = -\alpha y$.
  Note that
  \[
    \alpha = \frac{||y'||}{||y||} > \frac{R_1}{R_2}.
  \]

  Since $\Omega$ is convex and contains $\Omega \cap W^\perp$, it
  contains the open line segment
  \[
    \{t(-\alpha y) + (1 - t)(y + p) : t \in (0,1)\}.
  \]
  This line segment passes through $W$ when
  $t = \frac{1}{1 + \alpha}$, meaning that $\Omega$ must contain the
  point
  \[
    p' = \left(1 - \frac{1}{1 + \alpha}\right)p.
  \]
  
  Since $\Omega$ contains the origin, it also contains
  \[
    \left(1 - \frac{1}{1 + R_1/R_2}\right)p = \frac{R_1}{R_1 + R_2}p.
  \]
  This point lies in $\Omega \cap W$, meaning that $p$ lies in
  $R_3 \cdot (\Omega \cap W)$ where
  \[
    R_3 := \frac{R_1 + R_2}{R_1}.
  \]
  So we can take our map $f$ to be the linear rescaling about the
  origin by $R_3$.
\end{proof}

\begin{prop}
  \label{prop:cone_dual_intersections_open}
  Let $V = W_a \oplus V_b \oplus x$, for $x \in \P(V)$.

  Let $\Omega_a, \Omega_a'$ be properly convex domains in
  $\P(W_a \oplus x)$, and let $\Omega_b''$ be a properly convex domain
  in $\P(V_b \oplus x)$ such that
  \begin{align*}
    x \in \Omega_a \subset \Omega_a', \quad x \in \Omega_b''.
  \end{align*}
  If $\overline{\Omega_a'}$ does not intersect $\P(W_a)$, then there
  exists a properly convex domain $\Omega_b'$ in $\P(V_b \oplus x)$ so
  that any $\Omega \in \mc{C}(V)$ which satisfies
  $\Omega \cap \P(V_b) = \emptyset$ and
  \begin{enumerate}
  \item \label{item:containment_alt_1} $\Omega_a' \supset \pi_{W_a \oplus x}(\Omega)$,
  \item \label{item:containment_alt_2} $\Omega_a \subset \Omega \cap
    \P(W_a \oplus x)$,
  \item[(3a)] \label{item:containment_alt_3a}
    $\Omega_b'' \supset \Omega \cap \P(V_b \oplus x)$
  \end{enumerate}
  also satisfies
  \begin{enumerate}
  \item[(3)] \label{item:containment_alt_3}
    $\Omega_b' \supset \pi_{V_b \oplus x}(\Omega).$
  \end{enumerate}
\end{prop}
\begin{proof}
  Let $H = W_a \oplus V_b$, and consider the affine chart
  $A = \P(V) - \P(H)$. We can choose coordinates and a Euclidean
  metric on this affine chart so that $W_a \oplus x$ and
  $V_b \oplus x$ map to complementary orthogonal subspaces $W_a$,
  $V_b$ of $A$, meeting at the origin. In these coordinates, the
  projectivizations of the projection maps $\pi_{W_a \oplus x}$,
  $\pi_{V_b \oplus x}$ correspond to the orthogonal projections to
  $W_a$ and $V_b$, respectively.

  Since $\P(W_a)$ does not intersect $\overline{\Omega_a'}$, the
  images of $\Omega_a'$ and $\Omega_a$ in $A$ are both bounded open
  convex subsets of $W_a$.

  Let $\Omega$ be a properly convex domain not intersecting $\P(V_b)$
  and satisfying assumptions (\ref{item:containment_alt_1}),
  (\ref{item:containment_alt_2}), (3a). Since
  $\pi_{W_a \oplus x}(\Omega)$ is contained in $A$, $\Omega$ cannot
  intersect $\P(W_a \oplus V_b)$, so $\Omega$ is contained in the
  affine chart $A$ (although its closure need not be).

  In particular, $\Omega \cap \P(V_b \oplus x)$ is contained in the
  unique connected component of $\Omega_b'' \cap A$ which contains
  $x$. So, by replacing $\Omega_b''$ with this connected component, we
  may assume that the image of $\Omega_b''$ in $A$ is a convex open
  subset of $V_b$.

  Lemma \ref{lem:convex_subset_intersect_project} then implies that
  there is an affine map $f:A \to A$, depending only on $\Omega_a$ and
  $\Omega_a'$, so that
  \[
    \pi_{V_b \oplus x}(\Omega) \subseteq f(\Omega_b'').
  \]
  So, we can take $\Omega_b'$ to be the properly convex domain
  $f(\Omega_b'')$.
\end{proof}

We are now ready to prove Proposition \ref{prop:relative_benzecri}.

\begin{proof}[Proof of Proposition \ref{prop:relative_benzecri}]
  As in the proof of Lemma~\ref{lem:benzecri_1_dim}, properness is
  immediate from the Benz\'ecri cocompactness theorem, so we just need
  to show cocompactness. We let $V_a$, $V_b$, and
  $\mc{K}_a \subset \mc{C}_*(V_a)$ be as in the statement of the
  theorem. Let $(\Omega_n, x_n)$ be a sequence of properly convex
  domains in $\mc{C}_*(V_a, V_b, \mc{K}_a)$. We can choose a
  subsequence so that $x_n \to x$. Our goal is to find a pair of
  properly convex domains $\Omega_1, \Omega_2$ (with $x \in \Omega_1$)
  and $h_n \in \GL(V_b)$, so that up to a subsequence,
  \[
    \Omega_1 \subset h_n \cdot \Omega_n \subset \Omega_2.
  \]

  This will be sufficient, because $x_n \in \P(V_a)$, so
  $h_nx_n = x_n$ converges to $x$ and $h_n\Omega_n$ sub-converges to
  some properly convex domain $\Omega$ containing $\Omega_1 \ni x$.
  
  \begin{figure}[h]
    \centering
    \def\svgwidth{4in}
    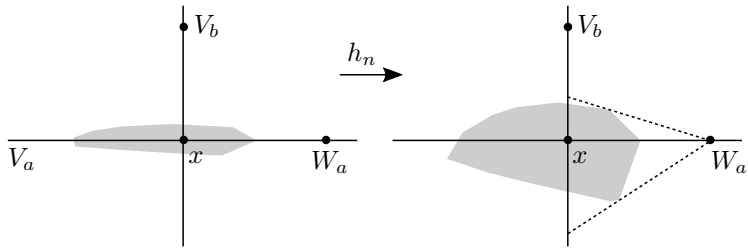
    \caption{Applying an element $h_n \in \GL(V_b)$ ``rescales'' in
      $\P(V_b \oplus x)$ about $x$; if the size of the intersection
      $\Omega \cap \P(V_b \oplus x)$ is bounded, then the size of the
      projection to $V_b \oplus x$ (with respect to the decomposition
      $V = W_a \oplus V_b \oplus x$) is also bounded (Proposition
      \ref{prop:cone_dual_intersections_open}).}
    \label{fig:rescaling_argument}
  \end{figure}

  Consider the sequence of domains
  $\Omega_n' = \Omega_n \cap \P(V_b \oplus x)$. We know $\P(V_b)$ is
  disjoint from $\Omega_n$ for all $n$. So, Lemma
  \ref{lem:benzecri_1_dim} implies that we can find $h_n \in \GL(V_b)$
  so that the domains $h_n\Omega_n'$ sub-converge in
  $\mc{C}(V_b \oplus x)$ to some domain $\Omega'$ in
  $\P(V_b \oplus x)$. In particular, up to a subsequence, we can find
  fixed domains $\Omega_b, \Omega_b'' \subset \P(V_b \oplus x)$ such
  that for all $n$,
  \[
    x \in \Omega_b \subset h_n\Omega_n' \subset \Omega_b''.
  \]
  
  Since the intersections $\Omega_n \cap \P(V_a)$ and projections
  $\pi_{V_a}(\Omega_n)$ both lie in a fixed compact set in
  $\mc{C}(V_a)$, we can also assume that there are domains
  $\Omega_a, \Omega_a' \in \mc{C}(V_a)$ so that for all $n$,
  \[
    \Omega_a \subset \Omega_n \cap \P(V_a), \quad \Omega_a' \supset
    \pi_{V_a}(\Omega_n).
  \]
  Since the action of any $h_n \in \GL(V_b)$ fixes $V_a$ pointwise and
  commutes with projection to $V_a$, this immediately implies that for
  all $n$,
  \[
    \Omega_a \subset h_n\Omega_n \cap \P(V_a), \quad \Omega_a' \supset
    \pi_{V_a}(h_n\Omega_n).
  \]

  Fix a subspace $W_a \subset V_a$ so that $V_a = W_a \oplus x$ and
  $\P(W_a)$ does not intersect the closure of $\Omega_a'$. This allows
  us to define a projection map
  $\pi_{V_b \oplus x}:V \to V_b \oplus x$, whose kernel is
  $W_a$. Proposition \ref{prop:projection_sharp} implies that
  $\pi_{V_b \oplus x}(h_n\Omega_n)$ is a properly convex open subset
  of $\P(V_b \oplus x)$, and Proposition
  \ref{prop:cone_dual_intersections_open} implies that for all $n$,
  $\pi_{V_b \oplus x}(h_n\Omega_n)$ is contained in a properly convex
  domain $\Omega_b' \subset \P(V_b \oplus x)$, depending only on
  $\Omega_a$, $\Omega_a'$, and $\Omega_b''$. Then we can apply Lemma
  \ref{lem:projection_intersection} to the domains
  $\Omega_a, \Omega_a', \Omega_b$, $\Omega_b'$ to finish the proof.
\end{proof}



\section{Cocompactness implies expansion}
\label{sec:cocompactness_implies_expansion}

The main goal of this section is to prove the implication
(\ref{item:main_thm_1}) $\implies$ (\ref{item:main_thm_2}) of Theorem
\ref{thm:convex_cocompact_equals_expansion}. In fact we will prove a
slightly more general statement:

\begin{prop}
  \label{prop:naive_expanding_in_supports}
  Let $C$ be a convex subset of a properly convex domain $\Omega$, and
  suppose that $\Gamma \subseteq \Aut(\Omega)$ acts cocompactly on
  $C$. Then $\Gamma$ is uniformly expanding at the faces of the ideal
  boundary of $C$.
\end{prop}

Afterwards, we will use some of the ideas arising in the proof to show
that a version of ``north-south dynamics'' holds for certain sequences
of elements in a convex cocompact group (Proposition
\ref{prop:line_translation_implies_north_south}).

\subsection{Pseudo-loxodromic elements}
\label{subsec:pseudo_loxodromics}

Our main inspiration comes from an observation in Sullivan's study
\cite{sullivan1979density} of conformal densities on $\H^d$: if
$x_0 \in \H^d$ is a basepoint defining a visual metric on $\dee \H^d$,
and $\gamma$ is any isometry of $\H^d$ not fixing $x_0$, then $\gamma$
expands a small ball in $\dee \H^d$ at the endpoint of the geodesic
ray from $x_0$ to $\gamma^{-1} x_0$, with expansion constant related
to $d(x_0, \gamma^{-1}x_0)$.

This observation relies on the fact that, given distinct points
$x, y \in \H^d$, there is a loxodromic isometry taking $x$ to $y$
whose axis is the geodesic joining $x$ and $y$. The exact analogue of
this fact for properly convex domains does not hold in general, since
there is no reason to expect even the full automorphism group of a
properly convex domain to act transitively on the domain. However,
instead of looking for actual automorphisms of the domain, we can
instead look for elements of $\PGL(d, \R)$ that do not perturb the
domain ``too much.'' We make this precise below.

\begin{defn}
  Let $\Omega \subset \RP^{d-1}$ be a properly convex domain, and let
  $\mc{K}$ be a compact subset of $\mc{C}(\Rd)$ containing
  $\Omega$. An element $g \in \PGL(d, \R)$ is a
  \emph{$\mc{K}$-pseudo-automorphism} of $\Omega$ if
  $g \Omega \in \mc{K}$.
\end{defn}

\begin{defn}
  Let $\Omega \subset \RP^{d-1}$ be a properly convex domain. For a
  compact subset $\mc{K} \subset \mc{C}(\Rd)$ containing $\Omega$, we
  say that a $\mc{K}$-pseudo-automorphism $g \in \PGL(d, \R)$ is
  \emph{$\mc{K}$-pseudo-loxodromic} if there is a $g$-invariant direct
  sum decomposition
  \[
    \Rd = V_- \oplus V_0 \oplus V_+,
  \]
  where:
  \begin{enumerate}[label=(\roman*)]
  \item the subspaces $V_-, V_+$ are positive eigenspaces of $g$ and
    supporting subspaces of $\Omega$,
  \item the convex hull of $\P(V_+) \cap \dee \Omega$ and
    $\P(V_-) \cap \dee \Omega$ has nonempty intersection with
    $\Omega$, and
  \item the projective subspace $\P(V_- \oplus V_+)$ intersects every
    $\Omega'$ in $\mc{K}$.
  \end{enumerate}
  The subspaces $V_-$ and $V_+$ are referred to as \emph{endpoints} of
  $g$. The projective subspace $\P(V_- \oplus V_+)$ is the \emph{axis}
  of the pseudo-loxodromic, and $V_0$ is the \emph{neutral subspace}.
\end{defn}

A pseudo-loxodromic element preserves its axis $\P(V_- \oplus
V_+)$. When $V_-$ and $V_+$ are points in $\RP^{d-1}$, this axis is an
actual projective line.

We do \emph{not} assume that an individual pseudo-loxodromic element
attracts points on its axis towards either of its endpoints, since we
are only interested in the dynamics of \emph{sequences} of
pseudo-loxodromics.

If $g_n$ is a sequence of $\mc{K}$-pseudo-loxodromic elements with
common endpoints, then, up to a subsequence, the domains $g_n\Omega$
converge to a domain $\Omega_\infty$ in $\mc{K}$ which intersects the
common axis. In fact, we observe:
\begin{prop}
  \label{prop:pseudo_loxodromic_limit_contains_hull}
  Let $g_n$ be a sequence of $\mc{K}$-pseudo-loxodromic elements with
  common endpoints $V_+$, $V_-$. If $g_n\Omega$ converges to
  $\Omega_\infty$, then $\Omega_\infty$ contains the relative interior
  of the convex hull (in $\Omega$) of $\P(V_+) \cap \dee \Omega$ and
  $\P(V_-) \cap \dee \Omega$.
\end{prop}
\begin{proof}
  Let $W$ denote the subspace $V_+ \oplus V_-$. We know that the
  intersection $\P(W) \cap \overline{\Omega_\infty}$ is either
  contained in a face of $\Omega_\infty$, or else its relative
  interior is contained in $\Omega_\infty$. It must be the latter,
  since we know $\Omega_\infty$ is in $\mc{K}$ and by definition
  $\Omega' \cap \P(W)$ is nonempty for every $\Omega' \in \mc{K}$.

  We let $C$ denote the convex hull of $\P(V_+) \cap \dee \Omega$ and
  $\P(V_-) \cap \dee \Omega$. It now suffices to show that the
  relative interior of $C$ is contained in the relative interior of
  $\P(W) \cap \Omega_\infty$. First, suppose that every subsequence of
  the restriction of $g_n$ to $W$ has a further subsequence which
  converges to some $g \in \PGL(W)$. In this case, we know
  $g\overline{\Omega} \cap \P(W) = \overline{\Omega_\infty} \cap
  \P(W)$; then $gC$ lies in the relative interior of $\Omega_\infty$
  because $C$ lies in the relative interior of $\Omega \cap \P(W)$ by
  assumption.

  Otherwise, the ratio of the eigenvalues of $g_n$ on $V_+$ and $V_-$
  is unbounded. We let $W' \subset W$ be the vector space whose
  projectivization $\P(W')$ is the projective span of
  $\Omega \cap \P(V_+)$ and $\Omega \cap \P(V_-)$. The unboundedness
  of the eigenvalue ratio of $g_n$ implies that
  $g_n\overline{\Omega} \cap \P(W)$ must converge to a subset of
  $\P(W')$. But this limit is $\overline{\Omega_\infty} \cap \P(W)$,
  which has nonempty relative interior in $\P(W)$ because
  $\P(W) \cap \Omega_\infty$ is nonempty. This is only possible if
  $W' = W$, which means that $C$ also has nonempty relative interior
  in $\P(W)$. Since
  $\overline{C} \subset \overline{\Omega_\infty} \cap \P(W)$ we must
  therefore have $C \subset \Omega_\infty \cap \P(W)$.
\end{proof}

\begin{defn}
  \label{defn:repelling_endpoint}
  Let $\Omega$ be a properly convex domain, and let $g_n$ be a
  sequence of $\mc{K}$-pseudo-loxodromic elements with common
  endpoints $V_+$, $V_-$ and common neutral subspace $V_0$. We say
  that $V_-$ is a \emph{repelling endpoint} of the sequence $g_n$ if
  there is a sequence
  \[
    x_n \in \Omega \cap \P(V_- \oplus V_+)
  \]
  such that $g_nx_n = x$ for some $x \in \Omega$, and
  $x_n \to x_- \in \dee \Omega$ with
  \[
    V_- = \supp\strat(x_-).
  \]
\end{defn}

\subsection{Existence of repelling pseudo-loxodromics}

We will use pseudo-loxodromics to state an analogue (Lemma
\ref{lem:pseudo_loxodromics_exist}) of the fact that any two points in
$\H^d$ can be joined by the axis of a loxodromic isometry. First, we
need a lemma:
\begin{lem}
  \label{lem:projection_subspace_exists}
  Let $x_+$, $x_-$ be a pair of points in the boundary of a properly
  convex domain $\Omega \subset \P(V)$ such that
  $(x_-, x_+) \subseteq \Omega$. Let $\P(H_+)$, $\P(H_-)$ be
  supporting hyperplanes of $\Omega$ at $x_+$, $x_-$. Let
  $V_- = \supp(\strat(x_-))$, and let $W = V_- \oplus x_+$.

  There exists a (possibly trivial) subspace
  $H_0 \subset H_+ \cap H_-$ such that
  \begin{enumerate}
  \item $H_- = H_0 \oplus V_-$, and
  \item $\pi_{W}(\Omega)$ is properly convex, where $\pi_{W}:V \to W$
    is the projection with kernel $H_0$.
  \end{enumerate}
\end{lem}
Note that while $\P(H_0)$ does not intersect $\Omega$, the
intersection $\P(H_0) \cap \overline{\Omega}$ may be nonempty.
\begin{proof}
  First suppose that $V_- = x_-$. In this case, we take
  $H_0 = H_+ \cap H_-$, and $\pi_{W}(\Omega)$ is exactly the line
  segment $(x_-, x_+)$. On the other hand, if $V_-$ has codimension
  one, then $W = V_- \oplus x_+ = V$ and we can take $H_0$ to be
  trivial, so $\pi_W$ is the identity map.
  
  So now suppose that $V_-$ is neither a single point nor a hyperplane
  in $V$. Consider the properly convex set
  $\Omega_- = \dee \Omega \cap H_-$.  We know $\P(H_+ \cap H_-)$ is a
  codimension-one projective subspace of $\P(H_-)$. Because
  $(x_-, x_+) \subseteq \Omega$, $H_+ \cap H_-$ does not contain
  $V_-$.

  We also know $\P(H_+ \cap H_-)$ intersects $\overline{\Omega_-}$ in
  a (possibly empty) properly convex set. We know the projective
  subspace $V_-$ has dimension $k \ge 1$ and positive codimension in
  $H_-$, so there exists a codimension-$k$ projective subspace of
  $H_+ \cap H_-$ which does not intersect $\overline{\Omega_-}$ or
  $V_-$. Let $H_0$ be such a subspace; since $H_0$ is disjoint from
  $\overline{\Omega}$, we are done by Proposition
  \ref{prop:projection_sharp}.
\end{proof}

The following lemma is the main technical result in this section. It
implies in particular that every face in the boundary of a
properly convex domain is the repelling endpoint of \emph{some}
sequence of $\mc{K}$-pseudo-loxodromics.

\begin{lem}\label{lem:pseudo_loxodromics_exist}
  Let $\Omega$ be a properly convex domain, let $x_- \in \dee \Omega$,
  and let $L$ be a projective line intersecting $\Omega$, joining
  $x_-$ with some $x_+ \in \dee \Omega$, $x_+ \ne x_-$. Let
  $F_- = \strat(x_-)$.

  For any sequence $\{x_n\} \subset L$, with $x_n \to x_-$, up to a
  subsequence, there exists a compact set
  $\mc{K} \subset \mc{C}(\Rd)$, a subspace $H_0 \subset \Rd$, and a
  sequence of $\mc{K}$-pseudo-loxodromic elements $g_n$ in
  $\PGL(d, \R)$, with endpoints $\supp(F_-)$ and $x_+$ and neutral
  subspace $H_0$, such that $g_nx_n = x$ for a fixed
  $x \in L \cap \Omega$.
\end{lem}
\begin{proof}
  Our strategy is to start with the case that $F_-$ is
  codimension-one (so the neutral subspace $H_0$ is trivial), and
  then use Proposition \ref{prop:relative_benzecri} to extend to the
  general case.
  
  \subsubsection*{$F_-$ is codimension-one}
  Let $V_-$ be the support of $F_-$. For each $n$, we let
  $s_n \in \GL(d, \R)$ be the diagonal map
  \[
    \lambda_n\mr{id}_{x_+} \oplus \mr{id}_{V_-} = \begin{bmatrix}
      \lambda_n\\
      & \mr{id}_{V_-}
    \end{bmatrix}
  \]
  acting on $x_+ \oplus V_-$, where $\lambda_n \to \infty$ is chosen
  so that $s_nx_n = x$ for a fixed $x \in L \cap \Omega$.
  
  \begin{figure}[h]
    \centering
    \def\svgwidth{2in}
    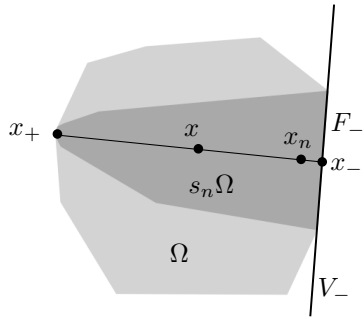
    \caption{Since $s_n$ attracts towards $x_+$ and repels from
      $V_-$, $s_n\Omega$ converges to the convex hull of $F_-$ and
      $x_+$.}
    \label{fig:codim_1_pseudo_loxodromic}
  \end{figure}
  
  The sequence of domains $s_n \cdot \Omega_n$ converges to a cone
  over $F_-$, with a cone point at $x_+$ (see Figure
  \ref{fig:codim_1_pseudo_loxodromic}). Since $F_-$ is a
  codimension-one face of $\Omega$, this cone is a properly convex
  domain containing $x$ in its interior.

  \subsubsection*{The general case}
  Let $V_-$ be the support of $F_-$, and let $H_+$, $H_-$ be
  supporting hyperplanes of $\Omega$ at $x_+$, $F_-$. Let
  $W = V_- \oplus x_+$. We choose a subspace
  $H_0 \subset H_+ \cap H_-$ as in Lemma
  \ref{lem:projection_subspace_exists} so that $H_- = V_- \oplus H_0$
  and $\pi_{W}(\Omega)$ is properly convex, where $\pi_{W} :V \to W$
  is the projection with kernel $H_0$.
  
  \begin{figure}
    \def\svgwidth{4in}
    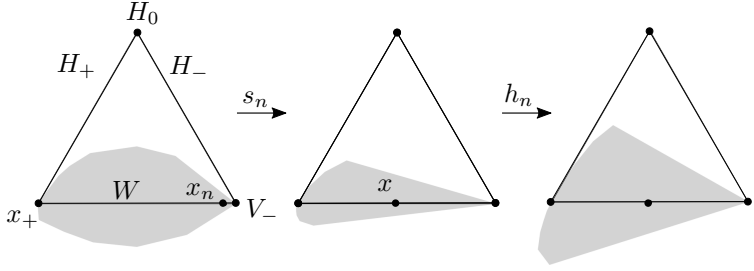
    \caption{To build the sequence of pseudo-loxodromic elements $g_n$,
      we push $x_n$ away from $x_-$ with $s_n \in \GL(W)$, ensuring that
      $s_n\Omega \cap \P(W)$ and $\pi_W(s_n\Omega)$ converge, and then
      use a ``correcting'' element $h_n \in \GL(H_0)$ to keep the domain
      from degenerating. Both $s_n$ and $h_n$ preserve the decomposition
      $\Rd = x_+ \oplus H_0 \oplus V_-$.}
    \label{fig:pseudo_loxodromic_elements}
  \end{figure}
  
  The domains
  \[
    \Omega \cap \P(W), \quad \pi_W(\Omega)
  \]
  are both properly convex open subsets of $\P(W)$ containing $F_-$ as
  a codimension-one face in their boundaries. Using the argument
  from the previous case, we can find group elements $s_n \in \GL(W)$
  so that
  \[
    s_n \cdot (\Omega \cap P(W)), \quad s_n \cdot \pi_W(\Omega)
  \]
  \emph{both} converge to properly convex domains in $\P(W)$
  containing a fixed $x = s_nx_n$ in $\Omega$.

  We extend $s_n$ linearly to the map $s_n \oplus \mr{id}_{H_0}$ on
  $W \oplus H_0$. Consider the sequence of properly convex domains
  \[
    \Omega_n = (s_n \oplus \mr{id}_{H_0}) \cdot \Omega.
  \]
  Since $s_n \oplus \mr{id}_{H_0}$ commutes with projection to $W$ and
  intersection with $W$, the sequences of pointed properly convex
  domains
  \[
    (\Omega_n \cap \P(W), x), \quad (\pi_W(\Omega_n), x)
  \]
  both converge in $\mc{C}_*(W)$. In particular, both of these
  sequences are contained in a fixed compact
  $\mc{K}_W \subset \mc{C}_*(W)$, and the pointed domains
  $(\Omega_n, x)$ all lie in the subset
  \[
    \mc{C}_*(W, H_0, \mc{K}_W)
  \]
  from Definition \ref{defn:domains_with_compact_slice}.

  Then, Proposition \ref{prop:relative_benzecri} (applied to the
  decomposition $\Rd = W \oplus H_0$) tells us that there is a
  sequence of group elements $h_n \in \GL(H_0)$ such that the pointed
  properly convex domains
  \[
    (\mr{id}_{W} \oplus h_n) \cdot (\Omega_n, x)
  \]
  lie in a fixed compact $\mc{K}$ in $\mc{C}_*(\Rd)$.

  Then, we can take our sequence of $\mc{K}$-pseudo-loxodromic
  elements $g_n$ to be the projectivizations of
  $(\mr{id}_W \oplus h_n) \cdot (s_n \oplus \mr{id}_{H_0}) = (s_n
  \oplus h_n)$.
\end{proof}

Next we examine some of the dynamical behavior of pseudo-loxodromic
sequences that have a repelling endpoint. Let $V$ be a normed vector
space. For any $g \in \GL(V)$, recall that the \emph{norm} and
\emph{conorm} of $g$ on $V$ are defined by
\[
  ||g|| = \sup_{v \in V - \{0\}}\frac{||gv||}{||v||}, \quad \mbf{m}(g)
  = \inf_{v \in V - \{0\}} \frac{||gv||}{||v||}.
\]

\begin{prop}
  \label{prop:pseudo_loxodromics_sv_gap}
  Let $g_n$ be a sequence of $\mc{K}$-pseudo-loxodromic elements with
  common endpoints $V_+, V_-$ and common neutral subspace $V_0$, and
  suppose that $V_-$ is a repelling endpoint (Definition
  \ref{defn:repelling_endpoint}) of the sequence $g_n$. Let
  $U_+ = V_+ \oplus V_0$. The sequence $g_n$ satisfies
  \begin{equation}
    \label{eq:sv_gap}
    \frac{\mbf{m}(g_n|_{U_+})}{||g_n |_{V_-}||} \to \infty.
  \end{equation}
\end{prop}
The ratio in (\ref{eq:sv_gap}) can be computed by fixing a norm on
$\Rd$, and then choosing a lift of each $g_n$ in $\GL(d,\R)$. The
value of this ratio does not depend on the choice of lift, and the
asymptotic behavior of the ratio does not depend on the choice of
norm.
\begin{proof}
  We can fix lifts $\tilde{g}_n$ of $g_n$ in $\GL(d, \R)$ which
  restrict to the identity on $V_-$. Our goal is then to show that
  \[
    \mbf{m}(\tilde{g}_n|_{U_+}) \to \infty,
  \]
  or equivalently, that 
  \[
    ||\tilde{g}_n^{-1}|_{U_+}|| \to 0.
  \]

  Suppose otherwise, so that for a sequence $v_n \in U_+$ with
  $||v_n|| = 1$, there is some $\eps > 0$ so that
  \[
    ||\tilde{g}_n^{-1} \cdot v_n|| \ge \eps.
  \]
  
  Let $x_n \in \Omega \cap \P(V_+ \oplus V_-)$ be a sequence so that
  $g_nx_n = x$ for some $x \in \Omega$ and $x_n \to x_-$, where $V_-$
  is the support of $\strat(x_-)$. We can choose a subsequence so that
  $g_n\Omega$ converges to some properly convex domain
  $\Omega_\infty$. We know $\Omega_\infty$ contains $x$ by Proposition
  \ref{prop:pseudo_loxodromic_limit_contains_hull}, so let $U$ be an
  open neighborhood of $x$ whose closure is contained in
  $\Omega_\infty$. We can find a lift $\tilde{x}$ of $x$ in $\Rd$ so
  that the projectivizations of each vector
  \[
    \tilde{x} \pm v_n
  \]
  lie in $U$, and thus in $\Omega_n$ for all sufficiently large
  $n$. Since $\tilde{g}_n$ restricts to the identity on $V_-$, the
  sequence $\tilde{g}_n^{-1} \tilde{x}$ converges to a lift
  $\tilde{x}_-$ of $x_-$.

  Then, up to a subsequence, the sequence of pairs of vectors
  \[
    \tilde{g}_n^{-1} \cdot (\tilde{x} \pm v_n)
  \]
  lies in a lift $\tilde{\Omega}$ of $\Omega$, and converges in $\Rd$
  to $\tilde{x}_- \pm v_\infty$, where $v_\infty \in U_+$ has norm at
  least $\eps$. This pair of points spans a nontrivial projective line
  segment in $\overline{\Omega}$ whose interior intersects the face
  $\strat(x_-)$ only at $x_-$, contradicting the definition of
  $\strat(x_-)$.
\end{proof}

Proposition \ref{prop:pseudo_loxodromics_sv_gap} implies in particular
that a sequence of $\mc{K}$-pseudo-loxodromic elements with repelling
subspace $V_-$ attracts generic points in $\RP^{d-1}$ to the
projective subspace $\P(U_+)$, and repels generic points away from
$\P(V_-)$. In fact, because the subspaces $U_+$ and $V_-$ are
transverse, and $g_n$ preserves each of them, the proposition also
implies that the sequence $g_n$ has \emph{expansion} behavior on the
Grassmannian $\Gr(k,d)$ in a neighborhood of $V_-$. To see this, we
need an estimate relating the ratio appearing in \eqref{eq:sv_gap} to
the metric behavior of $g_n$ on the Grassmannian.

To state the estimate, we choose an inner product on $\R^d$, and endow
$\RP^{d-1}$ with the metric $d_s$ obtained by setting $d_s(x, y)$ to
be the sine of the minimum angle between lifts of $x$ and $y$ in
$\R^d$. Then we let $d_k$ denote the metric on $\Gr(k, d)$ induced by
Hausdorff distance with respect to $d_s$.
\begin{lem}[See the appendix in \cite{bps2019anosov}, specifically
  Lemma A.10]
  \label{lem:sv_expansion_estimate}
  Let $U_+ \in \Gr(d-k,d)$ and $U_- \in \Gr(k,d)$ be transverse
  subspaces of $\R^d$, and let $g \in \PGL(d, \R)$. Suppose that for
  some $\alpha > 0$, we have $\angle(U_+, U_-) > \alpha$ and
  $\angle(gU_+, gU_-) > \alpha$.

  Then, for constants $b > 0$ and $\delta > 0$ (depending only on
  $\alpha$), if $W_1, W_2 \in \Gr(k, d)$ satisfy $d_k(W_i, V_-) <
  \delta$, then:
  \[
    d_k(gW_1, gW_2) \ge b\frac{\mbf{m}(g|_{U_+})}{||g|_{U_-}||}d(W_1,
    W_2).
  \]
\end{lem}

We can use this lemma to estimate the expansion behavior of elements
in our sequence $g_n$:
\begin{cor}
  \label{cor:pseudo_loxodromics_expand}
  Let $g_n$ be a sequence of $\mc{K}$-pseudo-loxodromic elements with
  common endpoints $V_+$, $V_-$ and common neutral subspace $V_0$, and
  suppose that $V_-$ is a repelling endpoint of the sequence, lying in
  $\Gr(k,d)$.

  Then for any Riemannian metric on $\Gr(k,d)$, and any $E > 1$, there
  exists $N \in \N$ such that if $n \ge N$, $g_n$ is $E$-expanding on
  some neighborhood of $V_-$ in $\Gr(k,d)$.
\end{cor}
\begin{proof}
  This follows directly from
  Proposition~\ref{prop:pseudo_loxodromics_sv_gap} and
  Lemma~\ref{lem:sv_expansion_estimate}, taking $U_+ = V_+ \oplus V_0$
  and $U_- = V_-$ and exploiting the fact that each $g_n$ preserves
  the decomposition $\R^d = V_- \oplus U_+$.
\end{proof}

\begin{remark}
  \label{rem:pseudo_loxodromic_use}
  In order to apply Lemma~\ref{lem:sv_expansion_estimate} to our
  situation, we need to know that our sequence of group elements $g_n$
  actually preserves the decomposition $V_- \oplus U_+$ (or at least
  that the sequence of decompositions $g_nV_- \oplus g_nU_+$ does not
  degenerate). This is why it is useful to have the additional control
  afforded by the pseudo-loxodromic sequences we have constructed---it
  is not enough to know merely that the sequence of domains
  $g_n\Omega$ does not degenerate.
\end{remark}

\subsection{Expansion}

Before we proceed, we fix some additional terminology:

\begin{defn}
  \label{defn:limit_along_line}
  Given a properly convex domain $\Omega$ and a point
  $x \in \dee \Omega$, we say that a sequence $x_n \in \Omega$
  \emph{limits to $x$ along a line $L$} if $x_n \to x$ in $\RP^{d-1}$,
  $L$ is an open projective line segment $(x, x') \subseteq \Omega$,
  and there exists a constant $R > 0$ such that
  \[
    d_\Omega(x_n, L) < R
  \]
  for all $n$.

  If the specific line $L$ is implied (or not relevant), we will just
  say that $x_n$ limits to $x$ along a line.

  If $F$ is some face of $\dee \Omega$, we say that $x_n$
  \emph{limits to $F$ along a line $L$} if every subsequence of $x_n$
  has a subsequence limiting to some $x \in F$ along $L$.
\end{defn}

\begin{remark}
  If $\Gamma$ is a group acting on a properly convex domain $\Omega$,
  and there are $\gamma_n \in \Gamma$ so that $\gamma_nx_0$ limits to
  $x$ along a line for some $x_0 \in \Omega$, the point $x$ is often
  referred to as a \emph{conical limit point} for the action of
  $\Gamma$ on $\dee \Omega$. We will avoid this terminology, since we
  will need to discuss conical limit points later in a way that is not
  exactly equivalent.
\end{remark}

\begin{prop}
  \label{prop:line_translation_implies_expansion}
  Let $\Omega$ be a properly convex domain and let
  $\Gamma \subseteq \Aut(\Omega)$. Let $F_-$ be a face of
  $\dee \Omega$, and let $x_n$ be a sequence in $\Omega$ limiting to
  $F_-$ along a line.

  If there exists $\gamma_n \in \Gamma$ so that $\gamma_n x_n$ is
  relatively compact in $\Omega$, then:
  \begin{enumerate}[label=(\alph*)]
  \item \label{item:compact_decomposition} There exists a compact set
    $K \subseteq \PGL(d, \R)$ such that $\gamma_n = k_ng_n$, where
    $k_n \in K$ and $g_n \in \PGL(d, \R)$ is a sequence of
    $\mc{K}$-pseudo-loxodromics with repelling endpoint $\supp(F_-)$.
  \item \label{item:uniform_expansion} For any Riemannian metric $d$
    on $\Gr(k,d)$, and any $E > 1$, for all sufficiently large $n$
    there is a neighborhood $U$ of $\supp(F_-)$ in $\Gr(k,d)$ such
    that $\gamma_n$ is $E$-expanding (with respect to $d$) on $U$.
  \end{enumerate}
\end{prop}

\begin{proof}
  Fix a compact $C \subset \Omega$ so that $\gamma_n x_n \in C$ for
  all $n$. We can move each $x_n$ by a bounded Hilbert distance so
  that it lies on a fixed line segment $L$ with an endpoint on
  $F_-$. So, by enlarging $C$ if necessary, we can assume that the
  points $x_n$ actually lie on the line $L$.
  
  Let $\mc{K}' \subset \mc{C}_*(\Rd)$ be the compact set
  $\{\Omega\} \times C$. By assumption we know that for all $n$, we
  have
  \[
    (\Omega, \gamma_nx_n) \in \mc{K}'.
  \]
  
  Using Lemma \ref{lem:pseudo_loxodromics_exist}, we can find a
  compact subset $\mc{K} \subset \mc{C}(\Rd)$ and a sequence $g_n$ of
  $\mc{K}$-pseudo-loxodromic elements with repelling endpoint
  $\supp(F_-)$ taking $x_n$ to $x$, for some $x \in \Omega \cap
  L$. The $g_n$ can be chosen so that the axis contains $L$, implying
  that the set
  \[
    \mc{K} \times \{x\} \subset \mc{C}_*(\Rd)
  \]
  is compact.

  Each group element $k_n = \gamma_ng_n^{-1}$ takes a pointed domain
  in the compact set $\mc{K} \times \{x\}$ to a pointed domain in the
  compact set $\mc{K}'$. But then, because $\PGL(d,\R)$ acts properly
  on $\mc{C}_*(\Rd)$, the $k_n$ lie in a fixed compact subset of
  $\PGL(d, \R)$. This proves part \ref{item:compact_decomposition}.
  
  Let $V_-$ be the support of $F_-$, and let $k = \dim V_-$. The
  elements $k_n$ can be viewed as lying in a compact subset of the
  diffeomorphisms of the compact manifold $\Gr(k, d)$. So, for any
  fixed Riemannian metric $d$ on $\Gr(k, d)$, there is a constant
  $M > 0$ so that for all $n$ and all $W_1, W_2 \in \Gr(k, d)$,
  \[
    d(k_n W_1, k_n W_2) > M \cdot d(W_1, W_2).
  \]
  
  Fix $E > 1$. Since $g_n$ has repelling endpoint $V_-$, Corollary
  \ref{cor:pseudo_loxodromics_expand} implies that for some
  sufficiently large $n$, there is a neighborhood $U$ of $V_-$ in
  $\Gr(k, d)$ so that $g_n$ satisfies
  \[
    d(g_n W_1, g_n W_2) > \frac{E}{M} \cdot d(W_1, W_2)
  \]
  for all $W_1, W_2 \in U$. But then we have
  \[
    d(\gamma_n W_1, \gamma_n W_2) > E \cdot d(W_1, W_2)
  \]
  giving us the required expansion.
\end{proof}

\begin{proof}[Proof of Proposition
  \ref{prop:naive_expanding_in_supports}]
  Let $\Gamma$ act cocompactly on some convex $C \subset \Omega$. Fix
  a Riemannian metric on $\Gr(k,d)$ and a constant $E > 1$.

  For every face $F$ of $\dee_iC$, there is a sequence $x_n$ in $C$
  limiting to $F$ along a line. Then part \ref{item:uniform_expansion}
  of Proposition \ref{prop:line_translation_implies_expansion} implies
  that if $\gamma_n x_n$ is relatively compact in $C$ for
  $\gamma_n \in \Gamma$, $\gamma_n$ is $E$-expanding on a neighborhood
  of $\supp(F)$ for sufficiently large $n$.
\end{proof}

\begin{proof}[Proof of (\ref{item:main_thm_1})
  $\implies$ (\ref{item:main_thm_2}) in Theorem
  \ref{thm:convex_cocompact_equals_expansion}]
  We apply Proposition \ref{prop:naive_expanding_in_supports} to
  $\Cor(\Gamma)$, whose ideal boundary is the full orbital limit set
  $\limset$. Lemma \ref{lem:full_orbital_limit_set_properties} implies
  that $\limset$ contains all of its faces and is closed and
  boundary-convex, so it is the $\Gamma$-invariant subset required by
  the theorem.
\end{proof}

\subsection{North-south dynamics}

In Section \ref{sec:relative_hyp_conv_cocompact}, it will be useful to
apply a consequence of part \ref{item:compact_decomposition} of
Proposition \ref{prop:line_translation_implies_expansion}. The
following can be thought of as a kind of weak version of north-south
dynamics on the limit set of a group acting on a convex projective
domain.

\begin{prop}
  \label{prop:line_translation_implies_north_south}
  Let $\Omega$ be a properly convex domain, let
  $\Gamma \subset \Aut(\Omega)$, and let $\Lambda$ be a closed
  $\Gamma$-invariant subset of $\dee \Omega$. Let $F$ be a face of
  $\Lambda$, and let $x_n$ be a sequence limiting to $F$ along a line.

  For any sequence $\gamma_n$ such that $\gamma_nx_n$ is relatively
  compact in $\Omega$, there exist subspaces $E_+$ and $E_-$, with
  $E_+ \oplus E_- = \Rd$, so that:
  \begin{enumerate}
  \item \label{item:supporting_subspaces} $\P(E_+)$, $\P(E_-)$ are
    supporting subspaces of $\Omega$, intersecting $\Lambda$,
  \item \label{item:compacts_converge} for every compact
    $K \subset \dee \Omega - \overline{F}$, a subsequence of
    $\gamma_n K$ converges uniformly to a subset of $\P(E_+)$, and a
    subsequence of $\gamma_n F$ converges uniformly to a subset of
    $\P(E_-)$,
  \item \label{item:segments_converge} for every $x \in F$ and every
    $z \in \dee \Omega - \overline{F}$, the sequence of line segments
    \[
      \gamma_n \cdot [x, z]
    \]
    sub-converges to a line segment intersecting $\Omega$.
  \end{enumerate}
\end{prop}
We emphasize again that the subspaces $E_{\pm}$ above are
\emph{complementary} in $\R^d$. Without this additional condition, the
proposition follows easily from (for example) \cite[Proposition
5.7]{iz2019flat}.
\begin{proof}
  Using Proposition \ref{prop:line_translation_implies_expansion}, we
  decompose each $\gamma_n$ as $k_ng_n$, for a sequence $g_n$ of
  $\mc{K}$-pseudo-loxodromic elements with repelling endpoint
  $V_- = \supp(F)$, and $k_n$ lying in a fixed compact in
  $\PGL(d, \R)$. Taking a subsequence, we may assume that $k_n$
  converges to $k \in \PGL(d, \R)$, so that
  \[
    \gamma_nV_- = k_ng_nV_- = k_nV_- \to kV_-.
  \]

  Let $E_- = kV_-$. We let $V_+$ be the other endpoint of the sequence
  $g_n$, let $V_0$ be the neutral subspace, and let
  $E_+ := k(V_+ \oplus V_0)$. Since $\Lambda$ is closed and
  $\Gamma$-invariant, both $\P(E_+)$ and $\P(E_-)$ intersect
  $\Lambda$.

  Fix a compact subset $K$ in $\dee \Omega - F$. Proposition
  \ref{prop:pseudo_loxodromics_sv_gap} implies that $g_nK$ converges
  uniformly to a subset of $\P(V_+ \oplus V_0)$. So, $k_ng_nK$
  converges uniformly to a subset of $\P(E_+)$.

  This shows parts (\ref{item:supporting_subspaces}) and
  (\ref{item:compacts_converge}). To see part
  (\ref{item:segments_converge}), let $L$ be the line segment $[x, z]$. By
  Proposition \ref{prop:strata_hilbert_bounded}, we can find $R > 0$
  and $x_n' \in L$ such that
  \[
    d_\Omega(x_n, x_n') \le R.
  \]
  We know that $\gamma_n x_n$ lies in a fixed compact subset $C$ of
  $\Omega$. So, $\gamma_nx_n'$ lies in a closed and bounded Hilbert
  neighborhood of $C$. This is also a compact subset of $\Omega$, so
  up to a subsequence, $\gamma_nx_n'$ converges to some
  $x_0' \in \Omega$.

  The limit of the line segment $[\gamma_nx_-, \gamma_nx_n']$ is
  nontrivial, intersects $\Omega$, and is a sub-segment of the limit of
  $[\gamma_nx_-, \gamma_n z]$, so this implies the desired result.
\end{proof}



\section{Background on relative hyperbolicity}
\label{sec:relative_hyp_background}

\subsection{A definition using convergence groups}

Relatively hyperbolic groups, like word-hyperbolic groups, have a wide
variety of possible definitions. Here we are most interested in the
dynamical properties of relatively hyperbolic groups, so we will use a
dynamical characterization due to Yaman \cite{yaman2004topological}.

Yaman's characterization uses the language of \emph{convergence group
  actions}, which we review below. Convergence groups were originally
studied in the context of group actions on spheres in $\Rd$ by Gehring
and Martin \cite{gm1987discrete}, and for general group actions on
compact Hausdorff spaces by Freden and Tukia
\cite{freden1997properties, tukia1998conical}.

\begin{defn}
  A group $\Gamma$ acting on a topological space $X$ is said to act on
  $X$ as a \emph{convergence group} if, for every sequence of distinct
  elements $\gamma_n \in \Gamma$, there exist (not necessarily
  distinct) points $a, b \in X$ and a subsequence $\gamma_n'$ of
  $\gamma_n$ such that the restriction of $\gamma_n'$ to $X - \{a\}$
  converges to the constant map $b$.
\end{defn}

When $X$ is a compact Hausdorff space, $\Gamma$ acts on $X$ as a
convergence group if and only if $\Gamma$ acts properly
discontinuously on the space of \emph{pairwise distinct triples} in
$X$ \cite{bowditch1999convergence}.

\begin{defn}
  Let $\Gamma$ act as a convergence group on $X$.
  \begin{itemize}
  \item We say that $x \in X$ is a \emph{conical limit point} if there
    exist \emph{distinct} points $a, b \in X$ and an infinite sequence
    of elements $\gamma_n \in \Gamma$ such that $\gamma_nx \to a$ and
    $\gamma_n y \to b$ for all $y \ne x$ in $X$.

  \item An infinite subgroup $H$ of $\Gamma$ is a \emph{parabolic
      subgroup} if it fixes a point $x \in X$ and each infinite-order
    element of $H$ has exactly one fixed point in $X$.

  \item A point $x \in X$ is a \emph{parabolic point} if its
    stabilizer is a parabolic subgroup. A parabolic point $x$ is
    \emph{bounded} if $\Stab_\Gamma(x)$ acts cocompactly on
    $X - \{x\}$.
  \end{itemize}
\end{defn}

When $\Gamma$ acts as a convergence group on a space $X$ with no
isolated points, and every point in $X$ is a conical limit point, we
say that $\Gamma$ acts as a \emph{uniform convergence group} on
$X$. This can be shown to be equivalent to the condition that $\Gamma$
act \emph{cocompactly} on the space of distinct triples in $X$
\cite{tukia1998conical}.

An important theorem of Bowditch \cite{bowditch1998topological} says
that if $\Gamma$ is a non-elementary group (i.e. not finite or
virtually cyclic), $\Gamma$ acts on a perfect metrizable compact space
$X$ as a uniform convergence group if and only if $\Gamma$ is
word-hyperbolic and $X$ is equivariantly homeomorphic to the Gromov
boundary of $\Gamma$. Yaman later proved an analogous result for
relatively hyperbolic groups:

\begin{thm}[\cite{yaman2004topological}]
  \label{thm:yaman_relative_hyp}
  Let $\Gamma$ be a non-elementary group, and let $\mc{H}$ be the
  collection of all conjugates of a finite collection of
  finitely-generated proper subgroups of $\Gamma$.

  Then $\Gamma$ is hyperbolic relative to $\mc{H}$ if and only if
  $\Gamma$ acts on a compact, perfect, and metrizable space $X$ as a
  convergence group, every point in $X$ is either a conical limit
  point or a bounded parabolic point for the $\Gamma$-action, and the
  parabolic points in $X$ are exactly the fixed points of the groups
  in $\mc{H}$.

  In this case, the Bowditch boundary $\dee(\Gamma, \mc{H})$ is
  equivariantly homeomorphic to $X$.
\end{thm}

We will use Theorem \ref{thm:yaman_relative_hyp} as our definition of
both relative hyperbolicity and the Bowditch boundary of a relatively
hyperbolic group. For other definitions, see
e.g. \cite{bowditch2012relatively}, \cite{ds2005tree}. The groups in
$\mc{H}$ are referred to as the \emph{peripheral subgroups}.

\begin{remark}
  Here we are adopting the convention that a group is hyperbolic
  relative to a conjugacy-closed collection of subgroups lying in
  finitely many conjugacy classes.

  The alternative convention would be to fix a finite set
  $\mc{P} = \{P_i\}$ of representatives for these conjugacy classes,
  and say that the group $\Gamma$ is hyperbolic relative to
  $\mc{P}$. We avoid this since we will work with the collection
  $\mc{H}$ of conjugates more often than we work with $\mc{P}$---the
  main exception is section \ref{subsec:nonperipheral_segments}.
\end{remark}

\section{Embedding the Bowditch boundary}
\label{sec:bowditch_bdry_embedding}

Our goal here is to prove Theorem
\ref{thm:bowditch_embedding_implies_conv_cocpct}. Our first step is
the following:
\begin{prop}
  \label{prop:no_segments_in_non_peripherals}
  Let $\Omega$ be a properly convex domain, and let
  $\Gamma \subset \Aut(\Omega)$ be hyperbolic relative to a collection
  of subgroups $\mc{H} = \{H_i\}$ each acting convex cocompactly on
  $\Omega$ with disjoint full orbital limit sets $\Lambda_\Omega(H_i)$.

  Suppose $\Lambda$ is a $\Gamma$-invariant subset of $\dee \Omega$
  containing all of its faces and containing $\Lambda_\Omega(H_i)$ for
  every $H_i$. If $\Heq{\Lambda}$ is the image of a
  $\Gamma$-equivariant embedding
  $\phi:\dee(\Gamma, \mc{H}) \to \Heq{\dee \Omega}$, then the set
  \[
    \Lambda_c = \Lambda - \bigcup_{H_i \in \mc{H}} \Lambda_\Omega(H_i)
  \]
  contains only extreme points in $\dee \Omega$.
\end{prop}
\begin{proof}
  The equivariant homeomorphism $\phi:\dee(\Gamma, \mc{H}) \to \Heq{\Lambda}$
  means that $\Gamma$ acts on $\Heq{\Lambda}$ as a convergence group
  as in Theorem \ref{thm:yaman_relative_hyp}. In particular, we can
  classify the points of $\Heq{\Lambda}$ as either bounded parabolic
  points or conical limit points, where the parabolic points are
  exactly the points corresponding to $\Lambda_\Omega(H_i)$.

  So, if $x$ is a point in $\Lambda_c$, it represents a conical limit
  point in $\Heq{\Lambda}$. Suppose for a contradiction that $x$ is
  not an extreme point, i.e. $x$ lies in the interior of a nontrivial
  segment $[a,b] \subset \dee \Omega$. Since $\Lambda$ contains all of
  its faces, $(a,b) \subset \Lambda$, and we can find
  $w,z \in \Lambda$ such that $w, x, z$ are pairwise distinct points
  lying on $(a,b)$ in that order.

  Lemma \ref{lem:full_orbital_limit_set_properties} tells us that each
  $\Lambda_\Omega(H_i)$ contains its faces, so we know that $w$ and
  $z$ cannot lie in any $\Lambda_\Omega(H_i)$. So $w$, $x$, and $z$
  represent three distinct points in $\Heq{\Lambda}$.

  This means that there exist group elements $\gamma_n \in \Gamma$ so
  that $\gamma_n\Heq{x} \to a$, and $\gamma_n\Heq{z}$,
  $\gamma_n\Heq{w}$ both converge to some $b \in \Heq{\Lambda}$, with
  $a,b$ distinct.

  This convergence is only in $\Heq{\Lambda}$. However, since the
  Bowditch boundary $\dee(\Gamma, \mc{H})$ is always compact,
  $\Heq{\Lambda}$ is as well, and therefore its preimage $\Lambda$ in
  the compact set $\dee \Omega$ is compact too. So, up to a
  subsequence, we can assume that $\gamma_nx \to u$, and
  $\gamma_nz \to v_1$, $\gamma_nw \to v_2$, with
  \[
    \Heq{u} = a, \; \Heq{v_1} = \Heq{v_2} = b.
  \]

  The line segment $[w,z]$ must converge to the line segment
  $[v_1, v_2]$, which must contain $u$. If $v_1 = v_2$, this is
  clearly impossible without having $u = v_1 = v_2$. If $v_1 \ne v_2$,
  then $v_1, v_2$ both lie in $\Lambda_\Omega(H_i)$ for some
  $H_i$. Since each $\Lambda_\Omega(H_i)$ is boundary-convex (Lemma
  \ref{lem:full_orbital_limit_set_properties} again), $u$ must lie in
  $\Lambda_\Omega(H_i)$ as well, a contradiction.
\end{proof}

The above is important partly because of the following proposition,
which we will use repeatedly in the proof of both Theorem
\ref{thm:bowditch_embedding_implies_conv_cocpct} and its converse.
\begin{prop}
  \label{prop:no_segment_near_limit_set}
  Let $\Omega$ be a properly convex domain, and let $\Lambda$ be a
  boundary-convex subset of $\dee \Omega$ containing all of its
  faces. Let $\mc{H}$ be a collection of subgroups of $\Aut(\Omega)$
  acting convex cocompactly with disjoint full orbital limit sets in
  $\Omega$.

  If every point in
  $\Lambda_c = \Lambda - \bigcup_{H_i \in \mc{H}} \Lambda_\Omega(H_i)$
  is an extreme point, then for any $x, y \in \Lambda$ with
  $\Heq{x} \ne \Heq{y}$, the segment $(x,y)$ lies in $\Omega$.
\end{prop}
\begin{proof}
  We will prove the contrapositive, and show that if $x, y$ are
  distinct points in $\Lambda$ with $(x,y) \subset \dee \Omega$, then
  $\Heq{x} = \Heq{y}$.

  Assume $x, y \in \Lambda$ are distinct. Boundary-convexity tells us
  that if $(x,y) \subset \dee \Omega$, then $(x,y) \subset
  \Lambda$. Since we know $\Lambda_c$ only contains extreme points,
  some $u \in (x,y)$ lies in $\Lambda_\Omega(H_i)$ for some
  $H_i \in \mc{H}$. Since $H_i$ acts convex cocompactly on $\Omega$,
  Lemma \ref{lem:full_orbital_limit_set_properties} implies that
  $[x,y]$ lies in $\Lambda_\Omega(H_i)$, which means that
  $\Heq{x} = \Heq{y}$.
\end{proof}

The following proposition explains why we do not need to assume that
$\Gamma$ is discrete in the statement of Theorem
\ref{thm:bowditch_embedding_implies_conv_cocpct}.
\begin{prop}
  \label{prop:bowditch_embedding_implies_discrete}
  If $\Omega$, $\Gamma$, $\Lambda$ are as in Theorem
  \ref{thm:bowditch_embedding_implies_conv_cocpct}, and $\Gamma$ is
  non-elementary, then $\Gamma$ is discrete.
\end{prop}
\begin{proof}
  $\Gamma$ acts as a convergence group on $\Heq{\Lambda}$, so it acts
  properly discontinuously on the space of pairwise distinct triples
  in $\Heq{\Lambda}$, which we denote $\mc{T}(\Heq{\Lambda})$.

  The map
  \[
    \Gamma \times \mc{T}(\Heq{\Lambda}) \to \mc{T}(\Heq{\Lambda})
  \]
  given by the $\Gamma$-action is continuous, so $\Gamma$ is discrete.
\end{proof}

We are now able to prove Theorem
\ref{thm:bowditch_embedding_implies_conv_cocpct}.

\begin{proof}[Proof of Theorem \ref{thm:bowditch_embedding_implies_conv_cocpct}]
  Let $\Omega$, $\Gamma$, $\Lambda$, $\mc{H}$ be as in the hypotheses
  for Theorem \ref{thm:bowditch_embedding_implies_conv_cocpct}. We can
  assume that $\mc{H} \ne \{\Gamma\}$ and that $\Gamma$ is infinite
  (if not then the theorem is trivial). This means that
  $\dee(\Gamma, \mc{H})$ contains at least two points, and Proposition
  \ref{prop:no_segment_near_limit_set} implies that $\Cor(\Gamma)$ is
  nonempty.

  If $\Gamma$ is virtually infinite cyclic, the hypotheses of the
  theorem imply that the generator $\gamma$ of a finite-index cyclic
  subgroup fixes a pair of points $\{x,y\}$ in $\dee \Omega$ with
  $(x,y) \subset \Omega$; $\gamma$ acts as a translation in the
  Hilbert metric along the axis $(x,y)$. This action is properly
  discontinuous (so $\Gamma$ is discrete) and cocompact. Further,
  since $x$ and $y$ are extreme points, $\gamma^n z$ converges to
  either $x$ or $y$ as $n \to \pm \infty$ for all $z \in \Omega$, so
  $\Lambda_\Omega(\Gamma) = \{x,y\}$.

  So we may assume $\Gamma$ is non-elementary. Owing to Theorem
  \ref{thm:convex_cocompact_equals_expansion}, we only need to show
  that $\Gamma$ is expanding at the faces of $\Lambda$; in fact we
  will show directly that the expansion is uniform.

  Since each $H_i$ acts convex cocompactly on $\Omega$, Theorem
  \ref{thm:convex_cocompact_equals_expansion} means that $\Gamma$ is
  expanding in a neighborhood of the support of any face of
  $\Lambda_\Omega(H_i)$ for some $H_i$. In fact, we can assume that
  the expansion constants are uniform over all $H_i \in \mc{H}$ (see
  Remark \ref{remark:support_expanding_well_defined}), so we only need
  to consider the faces in
  \[
    \Lambda_c = \Lambda - \bigcup_{H_i \in \mc{H}}\Lambda_\Omega(H_i).
  \]
  Proposition \ref{prop:no_segments_in_non_peripherals} implies that
  each of these faces is actually just a point in $\dee \Omega$, whose
  support is equal to itself.

  Let $x$ be a point in $\Lambda_c$. We will build a sequence of
  points $x_n$ in $\Omega$ limiting to $x$ along a line (Definition
  \ref{defn:limit_along_line}), and show that the orbits
  $\Gamma \cdot x_n$ intersect a fixed compact set.

  Since $\Gamma$ is non-elementary, its Bowditch boundary contains at
  least three distinct conical limit points, so we can find
  $y,z \in \Lambda_c$ so that $\Heq{x}, \Heq{y}, \Heq{z}$ are pairwise
  distinct.

  Fix supporting hyperplanes $W, V$ of $\Omega$ at $x$ and $z$,
  respectively. Proposition \ref{prop:no_segment_near_limit_set}
  implies that $W \cap V$ does not contain $x, y$, or $z$, and that
  the line segment $(x,z)$ is in $\Omega$. The projective hyperplane
  \[
    H = (W \cap V) \oplus y
  \]
  intersects $(x,z)$ at a point $w \in \Omega$.

  Since $\Heq{x}$ is a conical limit point, we can find a sequence
  $\gamma_n \in \Gamma$ so that
  \[
    \gamma_n \Heq{x} \to a
  \]
  and
  \[
    \gamma_n \Heq{z}, \; \gamma_n \Heq{y} \to b
  \]
  for $a, b$ distinct. As in the proof of Proposition
  \ref{prop:no_segments_in_non_peripherals}, we can pick subsequences
  so that $\gamma_nx$, $\gamma_ny$, and $\gamma_nz$ all converge to
  points $x_\infty, y_\infty, z_\infty$ in $\Lambda$, and $\gamma_nW$
  and $\gamma_nV$ converge to supporting hyperplanes $W_\infty$,
  $V_\infty$ of $\Omega$ at $x_\infty$ and $z_\infty$.

  Since $x_\infty$ and $z_\infty$ represent distinct points of
  $\Heq{\Lambda}$, Proposition \ref{prop:no_segment_near_limit_set}
  implies that $W_\infty \cap V_\infty$ must not contain $x_\infty$ or
  $z_\infty$; for the same reason $y_\infty$ is also not contained in
  $W_\infty \cap V_\infty$.

  While $\Heq{z_\infty} = \Heq{y_\infty}$, it is not necessarily true
  that $y_\infty = z_\infty$. However, we do know that the segment
  $(y_\infty, x_\infty)$ cannot lie in $\dee \Omega$. So, the sequence
  \[
    \gamma_n (H \cap (x,z)) = \gamma_n w
  \]
  cannot approach $x_\infty$.

  Proposition \ref{prop:bowditch_embedding_implies_discrete} means that
  we know $\Gamma$ is discrete, and so its action on $\Omega$ is
  properly discontinuous. Thus $\gamma_nw$ must accumulate to an
  endpoint of $[x_\infty, z_\infty]$---and therefore to $z_\infty$.

  \begin{figure}[h]
    \centering

    \def\svgwidth{2.9in}
    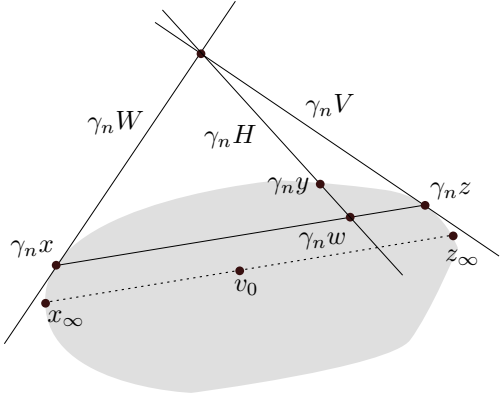
    
    \caption{The sequence $\gamma_n w$ limits to $z_\infty$, so the
      sequence $\gamma_n^{-1} v_0$ limits to $x$ along a line.}
    \label{fig:limiting_along_line}
  \end{figure}
  
  Let $\ell$ be the line segment $[x, z]$. This segment has a
  well-defined total order, where $a < b$ if $a$ is closer to $x$ than
  $b$. If $\ell_n = [\gamma_n x, \gamma_n z]$, then $\gamma_n$ is an
  order-preserving isometry from $\ell$ to $\ell_n$, where the metric
  is the restricted Hilbert metric $d_\Omega$.
  
  Fix a basepoint $v_0$ on the line segment
  $\ell_\infty = [x_\infty, z_\infty]$, and choose $v_n \in \ell_n$
  converging to $v_0$. Since $\gamma_nw$ converges to $z_\infty$, we
  see that $v_n < \gamma_nw$ and
  \[
    d_\Omega(v_n, \gamma_nw) \to \infty.
  \]
  Thus we must have $\gamma_n^{-1}v_n \to x$.

  But now we can apply part \ref{item:uniform_expansion} of
  Proposition \ref{prop:line_translation_implies_expansion} to the
  sequence $\gamma_n^{-1}v_n \subset \ell$ to see that $\gamma_n$ is
  eventually expanding in a neighborhood of $x$ in $\RP^{d-1}$.
\end{proof}

\section{Convex cocompact groups which are relatively hyperbolic}
\label{sec:relative_hyp_conv_cocompact}

The goal of this section is to prove Theorem
\ref{thm:no_segment_implies_rel_hyperbolic}.

\subsection{Non-peripheral segments in the boundary}
\label{subsec:nonperipheral_segments}

We start by showing that conditions
\ref{item:pairwise_disjoint_limit_sets} and
\ref{item:no_non_peripheral_segments} of Theorem
\ref{thm:no_segment_implies_rel_hyperbolic} are satisfied whenever
$\Gamma$ is a convex cocompact group hyperbolic relative to a
collection of convex cocompact subgroups. That is, we will show:
\begin{prop}
  \label{prop:relative_hyp_non_peripheral_segment}
  Let $\Gamma$ be a group hyperbolic relative to a collection $\mc{H}$
  of subgroups, and suppose that $\Gamma$ and each $H_i \in \mc{H}$
  act on a properly convex domain $\Omega$ convex cocompactly.

  Then:
  \begin{enumerate}[label=(\roman*)]
  \item \label{item:relhyp_disjoint_limit_sets} The full orbital limit
    sets $\Lambda_\Omega(H_i)$ are disjoint for distinct
    $H_i, H_j \in \mc{H}$,
  \item \label{item:relhyp_no_segment_limit_sets} Every nontrivial
    segment in $\Lambda_\Omega(\Gamma)$ is contained in the full
    orbital limit set of some peripheral subgroup $H_i$,
  \end{enumerate}
\end{prop}

We will closely follow the proof of a similar result of Islam and
Zimmer \cite[Theorem 1.8 (7)]{iz2019convex}. The main idea is that a
nontrivial segment $\ell$ in the full orbital limit set
$\Lambda_\Omega(\Gamma)$ of a convex cocompact group $\Gamma$ is
accumulated to by segments in the boundary of some maximal properly
embedded simplices in $\Cor(\Gamma)$. When $\Gamma$ is hyperbolic
relative to a collection $\mc{A}$ of virtually abelian subgroups of
rank at least $2$, Islam and Zimmer show that $\mc{A}$ is in
one-to-one correspondence with the set of maximal properly embedded
simplices in $\Cor(\Gamma)$, and then use a coset separation property
due to Dru\c{t}u and Sapir \cite{ds2005tree} to see that these maximal
properly embedded simplices are isolated. This ends up implying that
$\ell$ lies in the boundary of one of the simplices that accumulate to
it.

When we do \emph{not} assume the peripheral subgroups are virtually
abelian, we need to modify this approach slightly. First, we need to
assume that the peripheral subgroups act convex cocompactly on
$\Omega$ (\cite{iz2019convex} implies that this assumption is always
satisfied in the virtually abelian case). Second, in our situation,
the maximal properly embedded simplices in $\Cor(\Gamma)$ do not need
to be isolated. However, it is true that the convex cores $\Cor(H_i)$
of the peripheral subgroups in $\mc{H}$ are isolated. So the desired
result ends up following from the fact that every maximal properly
embedded $k$-simplex ($k \ge 2$) in $\Cor(\Gamma)$ lies in $\Cor(H_i)$
for some $H_i \in \mc{H}$; this is Lemma
\ref{lem:simplices_in_peripheral_subgroups} below.

\subsubsection{Cosets and convex cores of peripheral subgroups}

Let $\Gamma$ be hyperbolic relative to a collection of subgroups
$\mc{H}$, and suppose that $\Gamma$ and each $H_i \in \mc{H}$ act
convex cocompactly on a fixed properly convex domain $\Omega$. We fix
a basepoint $x \in \Omega$, and fix a finite set $\mc{P} = \{P_i\}$ of
conjugacy representatives for $\mc{H}$.

The \v{S}varc-Milnor lemma implies that $\Gamma$ is finitely generated
and that, under the word metric induced by any finite generating set,
$\Gamma$ is equivariantly quasi-isometric to the convex core
$\Cor(\Gamma)$ equipped with the restricted Hilbert metric
$d_\Omega$. The quasi-isometry can be taken to be the orbit map
$\gamma \mapsto \gamma x$.

Since each $P_i$ also acts convex cocompactly on $\Omega$, each $P_i$
is also finitely generated, and $P_i$ is quasi-isometric to
$\Cor(P_i)$, which isometrically embeds into $\Cor(\Gamma)$. We may
assume that the quasi-isometry constants are uniform over all
$P_i \in \mc{P}$, and fix a finite generating set for $\Gamma$
containing generating sets for each $P_i$.

Since $g \cdot \Cor(P_i) = \Cor(gP_ig^{-1})$, if we fix a
$\Gamma$-equivariant quasi-isometry
\[
  \phi:\Cor(\Gamma) \to \Gamma,
\]
we know $\phi$ restricts to a quasi-isometry
\[
  \Cor(gP_ig^{-1}) \to gP_i,
\]
with uniform quasi-isometry constants over all $g \in \Gamma$,
$P_i \in \mc{P}$.

The cosets $gP_i$ have a \emph{separation property}: distinct cosets
cannot stay ``close'' to each other over sets of large diameter. The
precise statement is as follows. For any metric space $X$, and any
$A \subseteq X$, we let
\[
  N_X(A; r)
\]
denote the open $r$-neighborhood of $A$ in $X$ with respect to the
metric $d_X$, and let
\[
  B_X(x; r)
\]
denote the open $r$-ball about $x \in X$.
\begin{thm}[{\cite[Theorem 4.1 ($\alpha_1$)]{ds2005tree}}]
  \label{thm:cosets_separated}
  Let $\Gamma$ be hyperbolic relative to $\mc{H}$, and let $\mc{P}$ be
  a finite set of conjugacy representatives. For every $r > 0$, there
  exists $R > 0$ such that for every distinct pair of left cosets
  $g_1P_1, g_2P_2$, the diameter of the set
  \[
    N_\Gamma(g_1P_1; r) \cap N_\Gamma(g_2P_2; r)
  \]
  is at most $R$.
\end{thm}

In addition, Theorem 1.7 of \cite{ds2005tree} implies that if
$k \ge 2$, any quasi-isometrically embedded $k$-flat in a relatively
hyperbolic group $\Gamma$ is contained in the $D$-neighborhood of a
coset $gP_i$ of some peripheral subgroup $P_i \in \mc{P}$. This allows
us to see the following:

\begin{lem}
  \label{lem:simplices_in_peripheral_subgroups}
  Suppose $\Gamma$ acts convex cocompactly on $\Omega$, and that
  $\Gamma$ is hyperbolic relative to a collection of subgroups
  $\mc{H}$ also acting convex cocompactly on $\Omega$. Every properly
  embedded $k$-simplex ($k \ge 2$) in $\Omega$ with boundary in
  $\Lambda_\Omega(\Gamma)$ is contained in $\Cor(H_i)$ for some
  $H_i \in \mc{H}$.
\end{lem}
\begin{proof}
  Each such embedded $k$-simplex $\Delta$ is a quasi-isometrically
  embedded $k$-flat in $\Cor(\Gamma)$, so $\phi(\Delta)$ is a
  quasi-isometrically embedded $k$-flat in $\Gamma$.
  \cite{ds2005tree}, Theorem 1.7 implies that $\phi(\Delta)$ is
  contained in a uniform neighborhood of $gP$ for some $P \in \mc{P}$.

  Applying a quasi-inverse of $\phi$ tells us that $\Delta$ is in a
  uniform Hilbert neighborhood of $\Cor(gPg^{-1})$ in $\Omega$. So the
  boundary of $\Delta$ is contained in $\dee_i \Cor(gPg^{-1})$, and
  $\Delta$ itself lies in $\Cor(gPg^{-1})$.
\end{proof}

We now quote:
\begin{lem}[{\cite[Lemma 15.4]{iz2019convex}}]
  \label{lem:nbhd_of_segment_in_simplex}
  Let $(u,v)$ be a nontrivial line segment in
  $\Lambda_\Omega(\Gamma)$, let $m \in (u,v)$ and
  $p \in \Cor(\Gamma)$, and let $V$ be the span of $(u,v)$ and
  $p$. For any $r > 0$, $\eps > 0$, there exists a neighborhood $U$ of
  $m$ in $\P(V)$ such that if $x \in U \cap \Cor(\Gamma)$, then there
  is a properly embedded simplex $S_x \subset \Cor(\Gamma)$ such that
  \[
    B_\Omega(x; r) \cap \P(V) \subset N_\Omega(S_x; \eps).
  \]
\end{lem}

Now we can prove Proposition
\ref{prop:relative_hyp_non_peripheral_segment}. The proof of part
\ref{item:relhyp_no_segment_limit_sets} is nearly identical to the
proof of Lemma 15.5 in \cite{iz2019convex}.

\begin{proof}[Proof of Proposition
  \ref{prop:relative_hyp_non_peripheral_segment}]
  \begin{aside}{\ref{item:relhyp_disjoint_limit_sets}}
    Let $H_i, H_j$ be a pair of peripheral subgroups in $\mc{H}$, and
    suppose that $\Lambda_\Omega(H_i) \cap \Lambda_\Omega(H_j)$
    contains a point $x \in \dee \Omega$. We can find a pair of
    projective-line geodesic rays in $\Cor(H_i)$ and $\Cor(H_j)$ with
    one endpoint at $x$. Proposition \ref{prop:strata_hilbert_bounded}
    implies that the images of these rays have finite Hausdorff
    distance.

    Thus, in $\Gamma$, a uniform neighborhood of the coset $g_iP_i$
    corresponding to $H_i$ contains an infinite-diameter subset of the
    coset $g_jP_j$ corresponding to $H_j$. So Theorem
    \ref{thm:cosets_separated} implies that $H_i = H_j$.
  \end{aside}

  \begin{aside}{\ref{item:relhyp_no_segment_limit_sets}}
    Consider any nontrivial segment $[u,v]$ in
    $\Lambda_\Omega(\Gamma)$, and fix $m \in (u,v)$ and
    $p \in \Cor(\Gamma)$. Theorem \ref{thm:cosets_separated} implies
    that for some $R > 0$, there exists $r > 0$ such that the diameter
    of
    \[
      N_{\Omega}(\Cor(H_i); r) \cap N_{\Omega}(\Cor(H_j); r)
    \]
    is less than $R$ whenever $H_i$ and $H_j$ are distinct.

    Let $V$ be the span of $u$, $v$, and $p$. Lemma
    \ref{lem:nbhd_of_segment_in_simplex} implies that for some
    neighborhood $U$ of $m$ in $\P(V)$, for every $x \in U$, there is
    some properly embedded simplex $S_x$ such that
    \[
      B_\Omega(x; R) \cap \P(V) \subset N_\Omega(S_x; r).
    \]
    Lemma \ref{lem:simplices_in_peripheral_subgroups} means that the
    simplex $S_x$ is contained in the convex hull $\Cor(H_x)$ of some
    peripheral subgroup $H_x$, and part
    \ref{item:relhyp_disjoint_limit_sets} implies that this peripheral
    subgroup is unique.
    
    We can shrink $U$ so that it is convex, and claim that in this
    case $H_x = H_y$ for all $x,y \in U \cap \Cor(\Gamma)$. By
    convexity, it suffices to show this when $d_\Omega(x, y) \le
    R/2$. Then
    \begin{align*}
      B_\Omega(x; R/2) \cap \P(V) \subset B_\Omega(y; R) \cap \P(V)
      \subset N_\Omega(S_y; r) 
    \end{align*}
    so the diameter of
    \[
      N_\Omega(\Cor(H_x); r) \cap N_\Omega(\Cor(H_y); r)
    \]
    is at least the diameter of $B_\Omega(x; R/2) = R$. Thus
    $H_x = H_y$.

    Fix $H = H_x$ for some $x \in U \cap \Cor(\Gamma)$. Then, if $x_n$ is
    a sequence in $\Cor(\Gamma)$ approaching $m$, there is a sequence
    $x_n' \in \Cor(H)$ such that
    \[
      d_\Omega(x_n, x_n') \le k,
    \]
    for $k$ independent of $n$. Up to a subsequence, $x_n'$ converges
    to some $x' \in \Lambda_\Omega(H)$. Proposition
    \ref{prop:strata_hilbert_bounded} implies that
    \[
      \strat(x') = \strat(m) \supseteq (u,v).
    \]
    $\Lambda_\Omega(H)$ contains $x'$. It is also closed and contains
    all of its faces (Lemma
    \ref{lem:full_orbital_limit_set_properties}). So
    $[u,v] \subset \Lambda_\Omega(H)$.
  \end{aside}
\end{proof}

\subsection{Convex cocompact and no relative segment implies
  relatively hyperbolic}
\label{sec:no_segment_implies_rel_hyperbolic}

We now turn to the rest of Theorem
\ref{thm:no_segment_implies_rel_hyperbolic}. As in our proof of
Theorem \ref{thm:bowditch_embedding_implies_conv_cocpct}, the main
tool will be Yaman's dynamical characterization of relative
hyperbolicity (Theorem \ref{thm:yaman_relative_hyp}). If $\Gamma$ is
virtually cyclic, Yaman's theorem does not apply, but in this case
$\Gamma$ is hyperbolic and the result follows from
\cite{dgk2017convex}.

Throughout the rest of this section, we assume (as in the hypotheses
to Theorem \ref{thm:no_segment_implies_rel_hyperbolic}) that $\Omega$
is a properly convex domain in $\RP^{d-1}$ preserved by a discrete
non-elementary group $\Gamma$ acting convex cocompactly with full
orbital limit set $\Lambda_\Omega(\Gamma)$, and $\mc{H}$ is a
conjugacy-invariant set of subgroups of $\Gamma$ lying in finitely
many conjugacy classes, with each $H_i \in \mc{H}$ acting convex
cocompactly on $\Omega$. We also assume $\mc{H} \ne \{\Gamma\}$, since
the result is trivial in this case.

We will prove the following:
\begin{prop}
  \label{prop:geom_finite_action}
  Suppose that conditions \ref{item:pairwise_disjoint_limit_sets},
  \ref{item:no_non_peripheral_segments}, and
  \ref{item:self_normalizing} of Theorem
  \ref{thm:no_segment_implies_rel_hyperbolic} hold for the collection
  of subgroups $\mc{H}$. Then:
  \begin{enumerate}
  \item \label{item:convergence_gp_action} $\Gamma$ acts as a
    convergence group on $\Heq{\limset}$,
  \item \label{item:topological_props} $\Heq{\limset}$ is compact,
    metrizable, and perfect,
  \item \label{item:parabolic_pts} the groups $H_i$ are parabolic
    subgroups, and their fixed points are bounded parabolic,
  \item \label{item:conical_limit_pts} every point in
    \[
      \Heq{\limset} - \{\Heq{\Lambda_\Omega(H_i)} : H_i \in \mc{H}\}
    \]
    is a conical limit point for the $\Gamma$-action on
    $\Heq{\limset}$.
  \end{enumerate}
\end{prop}

Since convex cocompact groups are always finitely generated, Theorem
\ref{thm:no_segment_implies_rel_hyperbolic} is a direct consequence of
Proposition \ref{prop:geom_finite_action}, Proposition
\ref{prop:relative_hyp_non_peripheral_segment}, and Theorem
\ref{thm:yaman_relative_hyp}.

\subsubsection{Dynamics of the $\Gamma$-action on $\limset$}

We start by establishing a basic dynamical fact about the action of
$\Gamma$ on $\limset$ and $\Heq{\limset}$. We need to recall some
basic properties of \emph{divergent sequences} (that is, sequences
which leave every compact) in $\PGL(d,\R)$.

\begin{defn}
  \label{defn:attract_repel}
  Let $g_n$ be a divergent sequence of elements in $\PGL(d, \R)$. We
  say that a pair of nontrivial subspaces $E_+, E_- \subset \R^d$ is a
  \emph{pair of attracting and repelling subspaces} for $g_n$ (or just
  an \emph{attracting/repelling pair}) if $\dim(E_+) + \dim(E_-) = d$,
  and there is a subsequence $g_m$ of $g_n$ so that for any compact
  $K \subset \RP^{d-1} - \P(E_-)$, the set $g_mK$ accumulates
  uniformly on $\P(E_+)$.
\end{defn}
Note that neither the subspaces $E_+, E_-$ nor even their respective
dimensions are uniquely determined by the sequence $g_n$ itself. We
also emphasize that while the subspaces in an attracting/repelling
pair have complementary dimension, they do \emph{not} need to be
transverse.

It is always possible to find at least one pair of attracting and
repelling subspaces for a divergent sequence $g_n \in \PGL(d, \R)$,
for instance by embedding $\PGL(d, \R)$ into the compact space
$\P(\End(\R^d))$. Then, if $g \in \P(\End(\R^d))$ is any accumulation
point of $g_n$, the image and kernel of $g$ form a pair of attracting
and repelling subspaces for the sequence. However, $g_n$ may have
other attracting/repelling pairs which do not arise in this way. For
example, if $g_n$ is given by the sequence of matrices
\[
  \begin{pmatrix}
    2^n\\
    & 1\\
    & & 2^{-n}
  \end{pmatrix},
\]
then the limit of $g_n$ in $\P(\End(\R^3))$ is the (projectivized)
matrix $\begin{psmallmatrix} 1 & 0 & 0\\
  0 & 0 & 0\\
  0 & 0 & 0\end{psmallmatrix}$, yielding the attracting subspace
$\langle e_1 \rangle$ and repelling subspace
$\langle e_2, e_3 \rangle$. But in this example, the subspace
$\langle e_1, e_2 \rangle$ is also an attracting subspace, paired with
the repelling subspace $\langle e_3 \rangle$.
  
The result below is certainly well-known, but we include a proof for
completeness.
\begin{lem}
  \label{lem:limiting_subspaces_support}
  Let $\Omega$ be a properly convex domain in $\RP^{d-1}$, let
  $\Gamma$ be a subgroup of $\Aut(\Omega)$, and let $\Lambda$ be any
  closed $\Gamma$-invariant subset of $\dee \Omega$ with nonempty
  convex hull in $\Omega$.

  If $E_+$ and $E_-$ are a pair of attracting and repelling subspaces
  for some divergent sequence $\{\gamma_n\} \subset \Gamma$, then
  $\P(E_+)$ and $\P(E_-)$ are supporting subspaces of $\Omega$ that
  intersect $\Lambda$ nontrivially.
\end{lem}
Note that if $E_+$ and $E_-$ are attracting and repelling subspaces
arising from a limit of $g_n$ in $\P(\End(\R^d))$, then this lemma can
be seen as a consequence of \cite[Proposition 5.6]{iz2019flat} (and in
fact this case of the lemma is already strong enough for our intended
application).
\begin{proof}
  It suffices to show the claim for $E_+$, because replacing
  $\gamma_n$ with $\gamma_n^{-1}$ reverses the role of the attracting
  and repelling subspaces.

  Since $\Omega$ is open, it is not contained in $\P(E_-)$. So, for
  some $x \in \Omega$, the limit of $\gamma_n x$ is contained in
  $\P(E_+)$. Since $\Omega$ is $\Gamma$-invariant, $\P(E_+)$
  intersects $\overline{\Omega}$ nontrivially.

  Let $E_+^*$ be the subspace of $(\R^d)^*$ consisting of functionals
  which vanish on $E_+$. Then $E_+^*$ is an attracting subspace for
  the sequence $\gamma_n$ under the dual action of $\Gamma$ on
  $(\R^d)^*$. So, by the previous argument, $\P(E_+^*)$ intersects
  $\overline{\Omega^*}$ nontrivially, which means $\P(E_+)$ cannot
  intersect $\Omega$.

  This shows that $\P(E_+)$ is a supporting subspace of $\Omega$ (and
  therefore $\P(E_-)$ is as well). To see that $\P(E_+)$ intersects
  $\Lambda$ nontrivially, note that since $\Lambda$ has nonempty
  convex hull in $\Omega$ and $\P(E_-)$ is a supporting subspace of
  $\Omega$, $\Lambda$ is not a subset of $\P(E_-) \cap \dee
  \Omega$. So, for some $x \in \Lambda$, $\gamma_n x$ accumulates to a
  point $y$ in $\P(E_+)$; since $\Lambda$ is $\Gamma$-invariant and
  closed, $y$ is in $\Lambda$ also.
\end{proof}

A straightforward consequence of Lemma
\ref{lem:limiting_subspaces_support} is part
(\ref{item:convergence_gp_action}) of Proposition
\ref{prop:geom_finite_action}:

\begin{prop}
  \label{prop:convergence_group_action}
  $\Gamma$ acts as a convergence group on $\Heq{\limset}$.
\end{prop}
\begin{proof}
  Let $\gamma_n$ be an infinite sequence in $\Gamma$. Since $\Gamma$
  is discrete, $\gamma_n$ is divergent, so we let $E_+, E_-$ be (the
  projectivizations of ) a pair of attracting and repelling subspaces,
  and extract a subsequence so that for any compact
  $K \subset \RP^{d-1} - E_-$, $\gamma_nK$ converges to a subset of
  $E_+$.

  Lemma \ref{lem:limiting_subspaces_support} implies that $E_+$ and
  $E_-$ are both supporting subspaces of $\Omega$ and both intersect
  $\limset$ nontrivially. The intersections $E_+ \cap \limset$ and
  $E_- \cap \limset$ are respectively the closures of subsets of a
  pair of faces $F_+, F_- \subset \limset$. By assumption, every
  face in $\limset$ containing a nontrivial projective segment lies
  in some $\Lambda_\Omega(H_i)$, so each face in $\limset$
  represents a single point of $\Heq{\limset}$. So we have
  \[
    \Heq{E_- \cap \limset} = a, \quad \Heq{E_+ \cap \limset} = b
  \]
  for (not necessarily distinct) points $a,b \in \Heq{\limset}$.

  Let $\Heq{K}$ be a compact subset of $\Heq{\limset} - \{a\}$, where
  $K$ is the preimage of $\Heq{K}$ in $\limset$. Because $\limset$ is
  compact, so is $K$. Moreover, $K$ cannot intersect $E_-$. So,
  $\gamma_n \cdot K$ converges to a subset of $E_+ \cap \limset$, and
  $\gamma_n \Heq{K}$ converges to $b$.
\end{proof}

\subsubsection{Topological properties of $\Heq{\limset}$}

Next, we will check that $\Heq{\limset}$ satisfies each of the
properties in part (\ref{item:topological_props}) of Proposition
\ref{prop:geom_finite_action}. The first, compactness, is
immediate from the compactness of $\limset$.

Showing that $\Heq{\limset}$ is metrizable is
equivalent to showing that it is Hausdorff, since it is a quotient of
a compact metrizable space.

Let
\[
  \pi_{\mc{H}}:\limset \to \Heq{\limset}
\]
be the quotient map. We will show that if $a$ is a point in
$\Heq{\limset}$, then we can find arbitrarily small
open neighborhoods of $\pi_{\mc{H}}^{-1}(a)$ in
$\limset$ which are of the form $\pi_{\mc{H}}^{-1}(U)$
for $U \subset \Heq{\limset}$.

Our first step is the following:
\begin{lem}
  \label{lem:open_set_line_bounded}
  Fix any metric $d_\P$ on projective space. Let
  $a \in \Heq{\limset}$.

  For any $\eps > 0$, there exists a subset
  $W(a, \eps) \subset \limset$ satisfying:
  \begin{enumerate}
  \item \label{item:line_bounded_1}
    $W(a, \eps) = \pi_{\mc{H}}^{-1}(V)$ for some
    $V \subset \Heq{\limset}$,
  \item \label{item:line_bounded_2} $W(a, \eps)$ contains an open
    neighborhood of $\pi_{\mc{H}}^{-1}(a)$ in $\limset$, and
  \item \label{item:line_bounded_3} For every $z \in W(a, \eps)$, we
    have
    \[
      d_\P(z, \pi_{\mc{H}}^{-1}(a)) < \eps.
    \]
  \end{enumerate}  
\end{lem}
\begin{proof}
  Let $X_a = \pi_{\mc{H}}^{-1}(a)$. For any open set $U$ in
  $\limset$ containing $X_a$, we let $W(U)$ be the set
  \[
    \pi_{\mc{H}}^{-1}(\Heq{U}) = U \cup \{x \in \Lambda_\Omega(H_i) :
    \Lambda_\Omega(H_i) \cap U \ne \emptyset\}.
  \]
  $W(U)$ is a subset of $\limset$ satisfying conditions
  (\ref{item:line_bounded_1}) and (\ref{item:line_bounded_2}). We
  claim that for any given $\eps > 0$, $W(U)$ also satisfies condition
  (\ref{item:line_bounded_3}) as long as $U$ is sufficiently small.
  
  We proceed by contradiction. Suppose otherwise, so that there is
  some $\eps > 0$ so that for a shrinking sequence of open
  neighborhoods $U_n$ of $X_a$, there is some $H_n \in \mc{H}$ such
  that
  \[
    \Lambda_\Omega(H_n) \cap U_n \ne \emptyset,
  \]
  and $\Lambda_\Omega(H_n)$ contains a point $z_n$ such that
  $d_\P(z_n, X_a) \ge \eps$.

  We write $\Lambda_n = \Lambda_\Omega(H_n)$. We can choose a
  subsequence so that in the topology on nonempty closed subsets of
  projective space, $\Lambda_n$ converges to some closed subset of
  $\limset$, which we denote $\Lambda_\infty$, and $z_n$ converges to
  $z_\infty \in \Lambda_\infty$ such that
  $d_\P(z_\infty, X_a) \ge \eps$.

  The set $\Lambda_\infty$ intersects every open subset of $\limset$
  containing $X_a$, and since $X_a$ is a closed subset of a metrizable
  space, this means $\Lambda_\infty$ intersects $X_a$. We will get a
  contradiction by showing that in fact $z_\infty \in X_a$.

  We consider two cases:
  \begin{aside}{Case 1: $\CH(\Lambda_\infty)$ is nonempty}
    Since the groups in $\mc{H}$ lie in only finitely many conjugacy
    classes, up to a subsequence, the $H_n$ are all conjugate to each
    other, and we may assume that $\Lambda_n = \gamma_n \Lambda_0$ for
    a sequence $\gamma_n \in \Gamma$.
    
    We can find a sequence $x_n \in \CH(\Lambda_n)$ converging to some
    $x_\infty \in \CH(\Lambda_\infty)$. Since the action of $H_0$ on
    $\CH(\Lambda_0)$ is cocompact, there is some fixed $R > 0$ so that
    every $H_0$-orbit in $\CH(\Lambda_0)$ intersects the Hilbert ball
    of radius $R$ about $x_0$. Since $H_n$ is a conjugate of $H_0$ by
    an isometry of the Hilbert metric on $\Omega$, the same is true
    (with the same $R$) for every $x_n$, $H_n$, and $\Lambda_n$.

    So, we can find a sequence
    \[
      \mu_n \in \gamma_nH_n\gamma_n^{-1}
    \]
    so that $\mu_n \gamma_n x_0$ lies in the Hilbert ball of radius
    $R$ about $x_n$. Since $x_n$ converges to $x_\infty \in \Omega$,
    and $\Gamma$ acts properly discontinuously on $\Omega$, this means
    that a subsequence of $\mu_n\gamma_n$ is eventually
    constant. Because $\mu_n\gamma_n\Lambda_0 = \gamma_n\Lambda_0$,
    this means we can assume there is some fixed $\gamma \in \Gamma$
    so that
    \[
      \Lambda_\infty = \gamma \Lambda_0 = \Lambda_\Omega(\gamma H_0
      \gamma^{-1}).
    \]
    But then since the limit sets $\Lambda_\Omega(H_i)$ are disjoint,
    we must have $X_a = \Lambda_\infty$, which means
    $z_\infty \in X_a$.
  \end{aside}

  \begin{aside}{Case 2: $\CH(\Lambda_\infty)$ is empty}
    In this case, $\CH(\Lambda_\infty)$ must be contained in the
    closure of some face $F$ of $\dee \Omega$. We may choose this face
    minimally, which means that $\CH(\Lambda_\infty)$ is not contained
    in $\dee F$. Then, because $\limset$ is boundary convex, we know
    that $\limset \cap F$ is nonempty. Because $\limset$ contains all
    of its faces this means that $F \subset \limset$.

    If $F$ is a single point, then its closure is also a singleton,
    hence $\Lambda_\infty$ is the singleton $\{z_\infty\}$. So in this
    case $z_\infty$ lies in $X_a$. If $F$ is \emph{not} a single
    point, it contains a nontrivial segment. By assumption, this
    segment lies in $\Lambda_\Omega(H_i)$ for some $H_i$; since
    $\Lambda_\Omega(H_i)$ is closed and contains its faces, all of
    $\overline{F}$ lies in $\Lambda_\Omega(H_i)$ as well. But then
    $\Lambda_\Omega(H_i)$ intersects both $X_a$ and $z_\infty$. Since
    $X_a = \pi_{\mc{H}}^{-1}(a)$ we must have
    $X_a = \Lambda_\Omega(H_i) = \Heq{z_\infty}$ and therefore
    $z_\infty \in X_a$ in this case as well.
  \end{aside}
  
\end{proof}

\begin{prop}
  \label{prop:limit_quotient_hausdorff}
  $\Heq{\limset}$ is Hausdorff.
\end{prop}
\begin{proof}
  Let $a, a'$ be distinct points in $\Heq{\limset}$,
  and let $X_a, X_a'$ be the preimages of $a$ and $a'$ in
  $\limset$.

  Since $X_a$ and $X_a'$ are closed disjoint subsets of the metrizable
  space $\limset$, there is some $\eps > 0$ such that
  for any $x \in X_a$, $x' \in X_a'$,
  \[
    d(x, x') > 2\eps.
  \]

  For each $n \in \N$, we define a sequence of sets $U_n$ containing
  $X_a$ as follows. We let $U_0 = X_a$. Then, for each $n > 0$, we
  take $U_n$ to be the set
  \[
    \bigcup_{b \in \Heq{U_{n-1}}} W(b, \eps / 2^n),
  \]
  where $W(b, \eps/2^n)$ is the set given by Lemma
  \ref{lem:open_set_line_bounded}. Note that each $U_n$ is a set of
  the form $\pi_{\mc{H}}^{-1}(V)$ for some
  $V \subset \Heq{\limset}$; moreover, if $z \in U_n$,
  then
  \[
    d(z, U_{n-1}) < \eps / 2^n.
  \]

  Consider the set $U = \bigcup_{n \in \N} U_n$. This set is the
  preimage of some $V \subset \Heq{\limset}$, and it must be contained
  in an $\eps$-neighborhood of $X_a$. In addition, $U$ is open in
  $\limset$: if $z$ is in $U_n$, then $U_{n+1}$ contains
  $W(\Heq{z}, \eps / 2^{n+1})$, which in turn contains an open
  neighborhood of $z$.
  
  This means that $\Heq{U}$ is an open set in $\Heq{\limset}$
  containing $a$. We can similarly construct an open set $\Heq{U'}$
  containing $a'$ such that $U'$ is contained in an
  $\eps$-neighborhood of $X_a'$. We know $U$ and $U'$ are disjoint, so
  $\Heq{U}$ and $\Heq{U'}$ separate $a$ and $a'$.
\end{proof}

Next we show that the space $\Heq{\limset}$ is perfect,
i.e. $\Heq{\limset}$ contains no isolated points. 

\begin{prop}
  \label{prop:limset_perfect}
  $\Heq{\limset}$ is perfect.
\end{prop}
\begin{proof}
  Fix $a \in \Heq{\limset}$ and a representative $x$ of $a$. Let
  $F = \strat(x)$.

  Let $x_n$ be a sequence of points in $\Cor(\Gamma)$ converging to
  $x$ in $\RP^{d-1}$. Convex cocompactness means that for some
  $\gamma_n \in \Gamma$, $\gamma_n^{-1}x_n \in C$ for a fixed compact
  $C \subset \Omega$.

  This means that (up to a subsequence) for fixed $x_0 \in \Omega$,
  $\gamma_n x_0$ converges to a point in $F$. And since $\gamma_n$
  acts by Hilbert isometries, Proposition
  \ref{prop:strata_hilbert_bounded} implies that if $B$ is any open
  ball with finite Hilbert radius about $x_0$, $\gamma_n B$ converges
  uniformly to a subset of $F$.

  $\gamma_n$ is divergent in $\PGL(d,\R)$, so let $E_+$ and $E_-$ be a
  pair of attracting and repelling (projective) subspaces for the
  sequence $\gamma_n$. We know that $E_+$ and $E_-$ are supporting
  subspaces of $\Omega$, and that
  \[
    \Heq{E_- \cap \limset}, \quad \Heq{E_+ \cap \limset}
  \]
  are single points in $\Heq{\limset}$. Moreover, since an open subset
  of projective space converges under $\gamma_n$ to $F$, $E_+$
  intersects $F$, and $\Heq{E_+ \cap \limset} = a$. Let
  $b = \Heq{E_- \cap \limset}$.

  Since we assume $\mc{H} \ne \{\Gamma\}$, $\Heq{\limset}$ cannot be a
  single point, and since $\Gamma$ is non-elementary, $\Heq{\limset}$
  contains at least three points. So, we can find a pair of points
  $c_1, c_2 \in \Heq{\limset}$ such that $\{b, c_1, c_2\}$ are
  pairwise distinct. Both $c_1$ and $c_2$ have a representative which
  does not lie in $E_-$, so both $\gamma_n c_1$ and $\gamma_n c_2$
  converge to $a$; since $c_1 \ne c_2$, $a$ cannot be isolated.
\end{proof}

\subsubsection{Parabolic points in $\Heq{\limset}$}

Our next task is to verify part (\ref{item:parabolic_pts}) of
Proposition \ref{prop:geom_finite_action}---that is, to show
that points stabilized by our candidate peripheral subgroups are
bounded parabolic points.

\begin{prop}
  Each point $\Heq{\Lambda_\Omega(H_i)}$ in $\Heq{\limset}$ is a
  parabolic point for the action of $\Gamma$, with stabilizer $H_i$.
\end{prop}
\begin{proof}
  The fact that $H_i$ is self-normalizing implies that $H_i$ is
  exactly the stabilizer of $\Heq{\Lambda_\Omega(H_i)}$ in $\Gamma$:
  for general $g \in \Aut(\Omega)$,
  \[
    g \cdot \Lambda_\Omega(H_i) = \Lambda_\Omega(g H_i g^{-1}),
  \]
  and since we assume that the full orbital limit sets of distinct
  groups in $\mc{H}$ are disjoint, $g \in \Gamma$ preserves
  $\Lambda_\Omega(H_i)$ if and only if $g$ normalizes $H_i$.
  
  So we just need to check that the groups $H_i$ are parabolic. Let
  $\gamma \in H_i$ be an infinite-order element, so that $\gamma^n$ is
  a divergent sequence in $\PGL(d, \R)$. We want to show that $\gamma$
  does not fix any point in $\Heq{\limset}$ other than
  $\Heq{\Lambda_\Omega(H_i)}$.

  Let $E_+$ and $E_-$ be a pair of attracting and repelling projective
  subspaces for the sequence $\gamma^n$. Lemma
  \ref{lem:limiting_subspaces_support} implies that both $E_+$ and
  $E_-$ support $\Omega$ and intersect $\Lambda_\Omega(H_i)$
  nontrivially.

  Let $b \in \Heq{\limset} - \{\Heq{\Lambda_\Omega(H_i)}\}$, let
  $y \in \Lambda_\Omega(H_i) \cap E_-$, and let $x \in \limset$ be a
  representative of $b$. Proposition
  \ref{prop:no_segment_near_limit_set} implies that $x$ cannot lie in
  $E_-$, so $\gamma^n x$ converges to a point in $\limset \cap
  E_+$. Then $\gamma^n b$ converges to $\Heq{\Lambda_\Omega(H_i)}$,
  and in particular $\gamma$ does not fix $b$.
\end{proof}

We still need to show that the parabolic points
$\Heq{\Lambda_\Omega(H_i)}$ are \emph{bounded} parabolic points,
i.e. that $H_i$ acts cocompactly on
\[
  \Heq{\limset} - \{\Heq{\Lambda_\Omega(H_i)}\}.
\]
Our strategy is to show that the set
\[
  \Lambda_i = \limset - \Lambda_\Omega(H_i)
\]
is a closed subset of the interior of some convex open set
$\Omega_{H_i}$, such that the ideal boundary of $\Lambda_i$ in
$\Omega_{H_i}$ is exactly $\Lambda_\Omega(H_i)$. Then, we can use the
fact that $H_i$ is uniformly expanding in supports at
$\Lambda_\Omega(H_i)$ to see that the action of $H_i$ on $\Lambda_i$
is cocompact.

If $H_i$ is irreducible (or more generally, if we know that $H_i$
contains a \emph{proximal element}), then as a consequence of
\cite[Proposition 3.1]{benoist2000automorphismes} (or
\cite[Proposition 4.5]{dgk2017convex}), we can simply take
$\Omega_{H_i}$ to be the unique $H_i$-invariant \emph{maximal}
properly convex domain $\Omega_{\max}$ in $\RP^{d-1}$. Since we do not
know if $H_i$ contains a proximal element in general, we do not know
if such a maximal domain exists. So, we will construct $\Omega_{H_i}$
directly.

To do so, we consider the \emph{dual} full orbital limit set
$\Lambda_{\Omega^*}(\Gamma)$ of a group $\Gamma$ acting on a properly
convex domain $\Omega$. i.e. the full orbital limit set in $\Omega^*$
of $\Gamma$ viewed as a subgroup of $\Aut(\Omega^*)$. Each element of
$\Lambda_{\Omega^*}(\Gamma)$ is an element of $\dee \Omega^*$, i.e. a
supporting hyperplane of $\Omega$.

\begin{prop}
  \label{prop:dual_limit_set_pairs_with_limit_set}
  Let $\Gamma$ be any subgroup of $\Aut(\Omega)$.
  \begin{enumerate}
  \item \label{item:dual_pairs_1} For every
    $x \in \Lambda_\Omega(\Gamma)$ there exists
    $w \in \Lambda_{\Omega^*}(\Gamma)$ such that $w(x) = 0$.
  \item For every $w \in \Lambda_{\Omega^*}(\Gamma)$ there exists
    $x \in \Lambda_\Omega(\Gamma)$ such that $w(x) = 0$.
  \end{enumerate}
\end{prop}
The statement follows from e.g. Proposition 5.6 in \cite{iz2019flat};
we provide an alternative proof for convenience.
\begin{proof}
  The two statements are dual to each other, so we only need to prove
  (\ref{item:dual_pairs_1}).

  Given a point $x \in \Omega$, and $W \in \Omega^*$, we consider a
  quantity
  \[
    \delta_\Omega(x, W)
  \]
  defined in \cite{dgk2017convex} as follows:
  \[
    \delta_\Omega(x, W) = \inf_{z \in W}\{\min\{[a_z, x; b_z,
    z], [b_z,x; a_z, z]\},
  \]
  where $a_z$ and $b_z$ are the points in $\dee \Omega$ such that
  $a_z, x, b_z, z$ lie on a projective line. The function
  $\delta_\Omega(x,W)$ can be thought of as an
  $\Aut(\Omega)$-invariant measure of how ``close'' $x$ is to
  $\dee \Omega$, relative to the projective hyperplane $W$: it takes
  on nonzero values for $x \in \Omega$, $W \in \Omega^*$, and for
  fixed $W \in \Omega^*$ and $x_n$ converging to $\dee \Omega$,
  $\delta_\Omega(x_n, W)$ converges to $0$.

  We now take $z \in \Lambda_\Omega(\Gamma)$, and choose
  $\gamma_n \in \Gamma$, $z_0 \in \Omega$ so that
  $\gamma_n z_0 \to z$. Fix some $W_0 \in \Omega^*$, and consider the
  sequence $\gamma_n W_0$. Up to a subsequence, this converges to some
  $W \in \Lambda_\Omega^*(\Gamma)$.

  \begin{figure}[h]
  \centering

  \def\svgwidth{2.3in}
\begingroup%
  \makeatletter%
  \providecommand\color[2][]{%
    \errmessage{(Inkscape) Color is used for the text in Inkscape, but the package 'color.sty' is not loaded}%
    \renewcommand\color[2][]{}%
  }%
  \providecommand\transparent[1]{%
    \errmessage{(Inkscape) Transparency is used (non-zero) for the text in Inkscape, but the package 'transparent.sty' is not loaded}%
    \renewcommand\transparent[1]{}%
  }%
  \providecommand\rotatebox[2]{#2}%
  \newcommand*\fsize{\dimexpr\f@size pt\relax}%
  \newcommand*\lineheight[1]{\fontsize{\fsize}{#1\fsize}\selectfont}%
  \ifx\svgwidth\undefined%
    \setlength{\unitlength}{685.20230973bp}%
    \ifx\svgscale\undefined%
      \relax%
    \else%
      \setlength{\unitlength}{\unitlength * \real{\svgscale}}%
    \fi%
  \else%
    \setlength{\unitlength}{\svgwidth}%
  \fi%
  \global\let\svgwidth\undefined%
  \global\let\svgscale\undefined%
  \makeatother%
  \begin{picture}(1,0.87303006)%
    \lineheight{1}%
    \setlength\tabcolsep{0pt}%
    \put(0,0){\includegraphics[width=\unitlength,page=1]{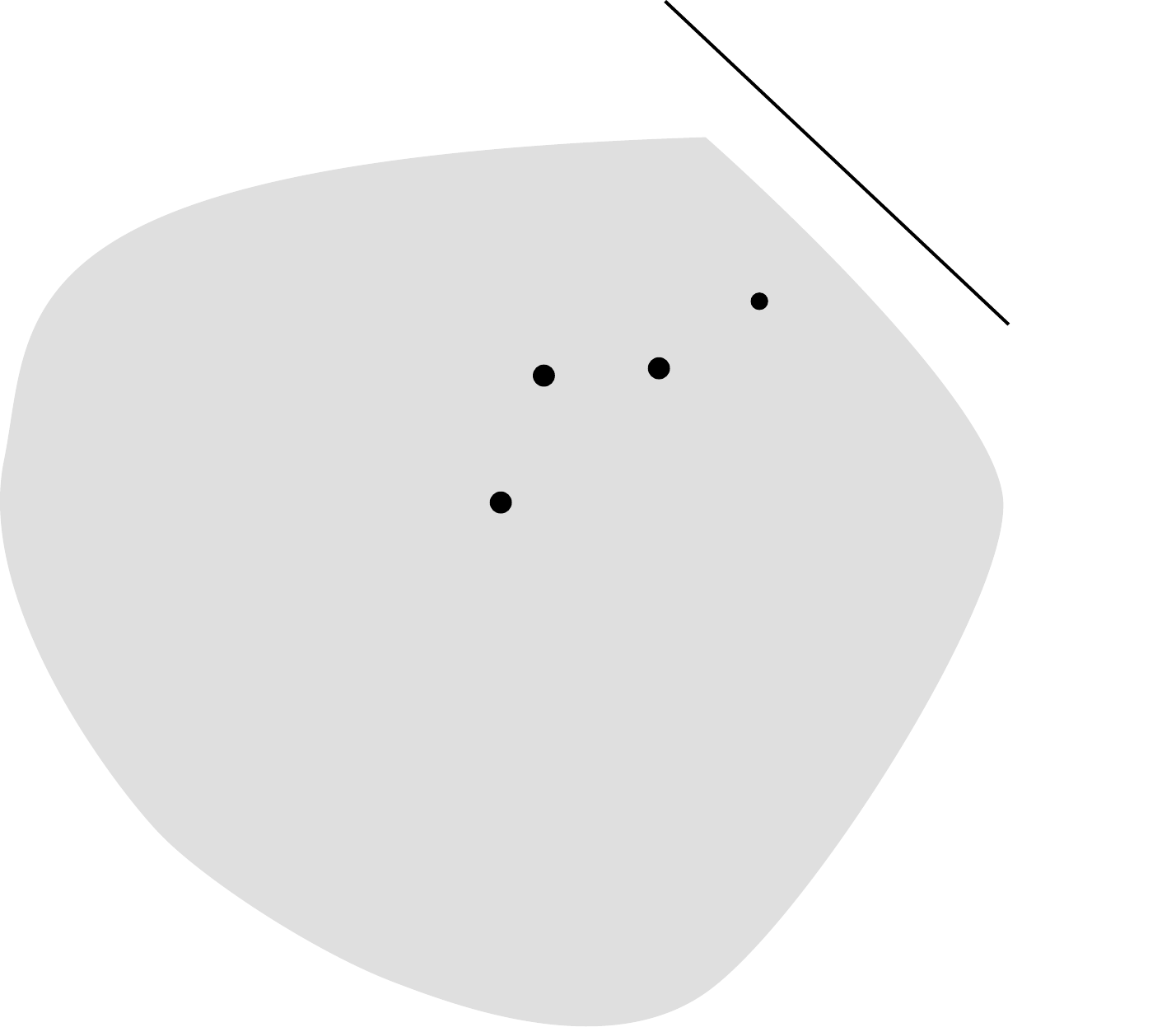}}%
    \put(0.21581833,0.02764793){\color[rgb]{0,0,0}\makebox(0,0)[lt]{\lineheight{1.25}\smash{\begin{tabular}[t]{l}$a_{y_n}$\end{tabular}}}}%
    \put(0.72317443,0.74720121){\color[rgb]{0,0,0}\makebox(0,0)[lt]{\lineheight{1.25}\smash{\begin{tabular}[t]{l}$y_n$\end{tabular}}}}%
    \put(0.63013462,0.56245864){\color[rgb]{0,0,0}\makebox(0,0)[lt]{\lineheight{1.25}\smash{\begin{tabular}[t]{l}$\gamma_n z_0$\end{tabular}}}}%
    \put(0.60635372,0.85292235){\color[rgb]{0,0,0}\makebox(0,0)[lt]{\lineheight{1.25}\smash{\begin{tabular}[t]{l}$\gamma_n W_0$\end{tabular}}}}%
    \put(0.59862954,0.69668038){\color[rgb]{0,0,0}\makebox(0,0)[lt]{\lineheight{1.25}\smash{\begin{tabular}[t]{l}$b_{y_n}$\end{tabular}}}}%
    \put(0,0){\includegraphics[width=\unitlength,page=2]{limit_set_dual.pdf}}%
  \end{picture}%
\endgroup%

  \caption{If $\gamma_n z_0$ approaches the boundary of $\Omega$, and
    $\delta_\Omega(\gamma_n z_0, \gamma_n W_0)$ is bounded away from
    $0$, $\gamma_n W_0$ must limit to a hyperplane containing the
    limit of $\gamma_n z_0$.}
  \label{fig:dual_limit_set_pairing}
\end{figure}

  Since $\delta_\Omega(x, W)$ is $\Gamma$-invariant, for any sequence
  \[
    y_n \in \gamma_n W_0,
  \]
  both of the cross-ratios
  \[
    [a_{y_n}, \gamma_n z_0; b_{y_n}, y_n], [b_{y_n}, \gamma_n z_0;
    a_{y_n}, y_n]
  \]
  remain bounded away from $0$ as $n \to \infty$. But since
  $\gamma_n z_0$ approaches $z \in \dee \Omega$, we can choose $y_n$
  so that exactly one of $a_{y_n}$, $b_{y_n}$ also approaches
  $z$. Thus, $y_n$ approaches $z$ as well, and so $W$ contains $z$.
\end{proof}

Next, we consider the \emph{dual convex core} for the $\Gamma$-action on
$\Omega$.
\begin{defn}
  \label{defn:dual_convex_hull}
  Let $\Omega \subset \RP^{d-1}$ be a properly convex domain, and let
  $\Gamma \subseteq \Aut(\Omega)$. The \emph{dual convex core}
  $\Cor^*(\Gamma)$ is the convex set
  \[
    [\mr{Hull}_{\Omega^*}(\Lambda_{\Omega^*}(\Gamma))]^*.
  \]
  Equivalently, $\Cor^*(\Gamma)$ is the unique connected component of
  \[
    \RP^{d-1} - \bigcup_{W \in \Lambda_{\Omega^*}(\Gamma)} W
  \]
  which contains $\Omega$.
\end{defn}

\begin{figure}[h]
  \centering
  \includegraphics[width=0.4\textwidth]{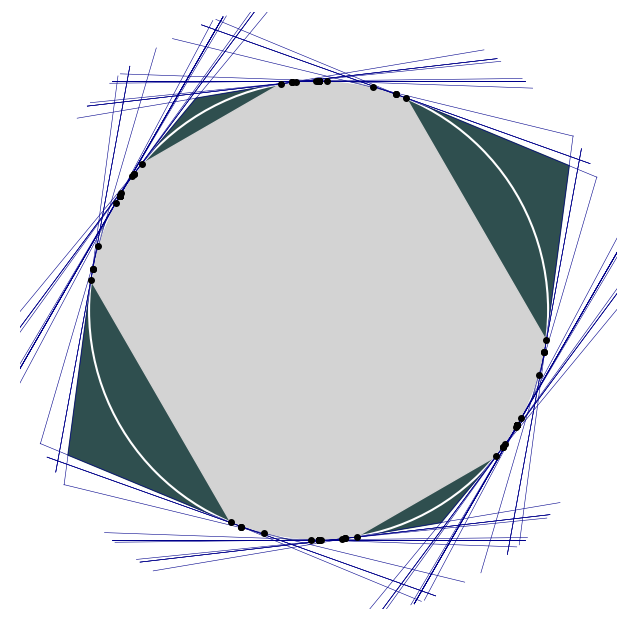}
  \caption{Part of the limit set and dual limit set for a group
    $\Gamma$ acting convex compactly on the projective model for
    $\H^2$ (the interior of the white circle). $\Cor(\Gamma)$ is the
    light region, and $\Cor^*(\Gamma)$ is the dark region.}
  \label{fig:convex_core}
\end{figure}

As long as $\Lambda_{\Omega^*}(\Gamma)$ contains at least two points,
$\Cor^*(\Gamma)$ does not contain all of $\RP^{d-1}$. It can be viewed
as an intersection of convex subspaces, so it is convex in the sense
of Definition \ref{defn:convex_domain}, but in general it is not
properly convex.

We can use the dual convex core to finish proving part
(\ref{item:parabolic_pts}) of Proposition
\ref{prop:geom_finite_action}.

\begin{prop}
  \label{prop:bounded_parabolic_points}
  The stabilizer of $\Heq{\Lambda_\Omega(H_i)}$ acts cocompactly on
  \[
    \Lambda_i = \Heq{\limset} - \{\Heq{\Lambda_\Omega(H_i)}\}.
  \]
\end{prop}
\begin{proof}
  Let $\Omega_{H_i} = \Cor^*(H_i)$ be the dual convex core of $H_i$ in
  $\Omega$. Proposition \ref{prop:dual_limit_set_pairs_with_limit_set}
  implies that $\Lambda_\Omega(H_i)$ lies in the boundary of
  $\Omega_{H_i}$.

  Moreover, the set $\limset - \Lambda_\Omega(H_i)$ lies in the
  interior of $\Omega_{H_i}$---for, every point in the boundary of
  $\Omega_{H_i}$ is contained in a projective hyperplane
  $W \in \Lambda_{\Omega^*}(H_i)$, and each such hyperplane supports
  some $x \in \Lambda_\Omega(H_i)$. Since $W$ is also a supporting
  hyperplane of $\Omega$, Proposition
  \ref{prop:no_segment_near_limit_set} implies that no
  $y \in \limset - \Lambda_\Omega(H_i)$ lies in $W$.

  $\limset$ is thus a closed subset of $\Omega_{H_i}$ whose ideal
  boundary in $\Omega_{H_i}$ is contained in
  $\Lambda_\Omega(H_i)$. Since $H_i$ acts convex cocompactly on
  $\Omega$, it is uniformly expanding in supports at
  $\Lambda_\Omega(H_i)$ by Theorem
  \ref{thm:convex_cocompact_equals_expansion}. Then Proposition
  \ref{prop:expansion_implies_cocpct_action} (applied to the convex
  domain $\Omega_{H_i}$) implies that the action of $H_i$ on
  $\limset - \Lambda_\Omega(H_i)$ is cocompact---which means that the
  $H_i$-action on the quotient
  $\Heq{\limset} - \{\Heq{\Lambda_\Omega(H_i)}\}$ is cocompact as
  well.
\end{proof}

\subsubsection{Conical limit points in $\Heq{\limset}$}

Finally we check part (\ref{item:conical_limit_pts}) of Proposition
\ref{prop:geom_finite_action}---that the remaining points in our
candidate Bowditch boundary are indeed conical limit points. We will
do this in two steps.

\begin{lem}
  \label{lem:bdry_cocompact_implies_north_south}
  Let $H_i \in \mc{H}$, let
  \[
    x_n \in \limset - \Lambda_\Omega(H_i)
  \]
  be a sequence approaching $x \in \Lambda_\Omega(H_i)$, and let
  $F = \strat(x)$. If $h_n$ is a sequence such that $h_n \Heq{x_n} $
  is relatively compact in
  \[
    \Heq{\limset} - \{\Heq{\Lambda_\Omega(H_i)}\},
  \]
  then for any compact
  \[
    K \subset \limset - \overline{F},
  \]
  $h_n$ sub-converges on $\Heq{K}$ to the constant map
  $\Heq{\Lambda_\Omega(H_i)}$.
\end{lem}
\begin{proof}
  Any such sequence $h_n$ must be divergent, so we let $E_+$ and $E_-$
  be a pair of attracting and repelling projective subspaces for
  $h_n$. We know from Lemma~\ref{lem:limiting_subspaces_support} that
  $E_+$ and $E_-$ are supporting subspaces of $\Omega$, each
  intersecting $\Lambda_\Omega(H_i)$
  nontrivially. Proposition~\ref{prop:no_segment_near_limit_set}
  implies that $E_+ \cap \limset \subset \Lambda_\Omega(H_i)$ and
  $E_- \cap \limset \subset \Lambda_\Omega(H_i)$. We can see that the
  subspace $E_-$ must contain $x$, since otherwise $h_nx_n$ would
  subconverge to a point in
  $E_+ \cap \dee \limset \subseteq \Lambda_\Omega(H_i).$

  But then $E_- \cap \dee \Omega$ is a subset of $\overline{F}$. Then,
  since $E_+$ and $E_-$ are a pair of attracting and repelling
  subspaces, if $K$ is compact in $\limset - E_-$, we know $h_nK$
  uniformly accumulates on
  $E_+ \cap \limset \subset \Lambda_\Omega(H_i)$.
\end{proof}

\begin{prop}
  Every element of the set
  \[
    \Heq{\limset} - \{\Heq{\Lambda_\Omega(H_i)} : H_i \in \mc{H}\}
  \]
  is a conical limit point for the action of $\Gamma$ on $\limset$.
\end{prop}
\begin{proof}
  By assumption, any point in this set has a unique representative
  $x \in \limset$ which is an extreme point in $\dee \Omega$. Fix a
  sequence $x_n \in \Omega$ limiting to $x$ along a line, and let
  $\gamma_n \in \Gamma$ be group elements taking $x_n$ back to some
  fixed compact in $\Omega$.

  Proposition \ref{prop:line_translation_implies_north_south} implies
  that there is a supporting subspace $E_+$ of $\Omega$, intersecting
  $\limset$, so that $\gamma_nx$ limits to some $x' \in \limset$ not
  intersecting $E_+$, and if $K$ is any compact subset of
  $\limset - x$, a subsequence of $\gamma_n K$ converges uniformly to
  a subset of $E_+ \cap \limset$. In particular, $\gamma_n$ converges
  uniformly on compacts in
  \[
    \Heq{\Lambda} - \{\Heq{x}\}
  \]
  to the constant map $\Heq{E_+ \cap \limset}$.

  If $\Heq{x'} \ne \Heq{E_+ \cap \limset}$, then we are done. However,
  it is also possible that $x'$ and $E_+ \cap \limset$ both lie in the
  same full orbital limit set of some convex cocompact subgroup $H_i$.

  In this case, we use the fact that $\Heq{\Lambda_\Omega(H_i)}$ is a
  bounded parabolic fixed point (Proposition
  \ref{prop:bounded_parabolic_points}) to find a sequence
  $h_n \in H_i$ such that $h_n \cdot \Heq{\gamma_n x}$ lies in a fixed
  compact set $C$ in
  \[
    \Heq{\limset} - \{\Heq{\Lambda_\Omega(H_i)}\},
  \]
  and consider the sequence of group elements $h_n\gamma_n$. We will
  show that $h_n\gamma_n$ is a conical limit sequence for $x$,
  i.e. that after taking a subsequence, for distinct
  $a, b \in \Heq{\limset}$, we have $h_n\gamma_n\Heq{x} \to a$ and
  $h_n\gamma_n\Heq{K} \to b$ for any compact
  $\Heq{K} \subset \Heq{\limset} - \Heq{x}$.

  So, fix an arbitrary compact subset $\Heq{K}$ of
  \[
    \Heq{\limset} - \{\Heq{x}\},
  \]
  where $K$ is the (compact) preimage of $\Heq{K}$ in
  $\limset - \{x\}$.

  After taking a subsequence, $\gamma_n K$ must converge to a compact
  subset of $E_+ \cap \limset$, which does not intersect $x'$. In
  fact, part (\ref{item:segments_converge}) of Proposition
  \ref{prop:line_translation_implies_north_south} implies that
  $\gamma_n K$ converges to a compact subset of
  $\limset - \overline{F'}$, where $F' = \strat(x')$. So there is a
  fixed compact
  \[
    K' \subset \limset - \overline{F'}
  \]
  so that for sufficiently large $n$, $\gamma_n K \subset K'$. Then
  Lemma \ref{lem:bdry_cocompact_implies_north_south} implies that
  \[
    h_n\gamma_n\Heq{K} \subseteq h_n \Heq{K'}
  \]
  sub-converges to $\Heq{\Lambda_\Omega(H_i)}$. But on the other hand,
  \[
    \Heq{(h_n\gamma_n) x_n} \in C
  \]
  sub-converges to some $b \ne \Heq{\Lambda_\Omega(H_i)}$, so
  $h_n\gamma_n$ gives us the sequence of group elements we need.
\end{proof}



\bibliographystyle{alpha}
\bibliography{references}

\newcommand{\etalchar}[1]{$^{#1}$}
\begin{thebibliography}{DGK{\etalchar{+}}21}

\bibitem[BDL15]{bdl2015convex}
Samuel Ballas, Jeffrey Danciger, and Gye-Seon Lee.
\newblock Convex projective structures on non-hyperbolic three-manifolds.
\newblock {\em Geometry \& Topology}, 22, 08 2015.

\bibitem[Ben60]{benzecri1960varietes}
Jean-Paul Benz\'{e}cri.
\newblock Sur les vari\'{e}t\'{e}s localement affines et localement
  projectives.
\newblock {\em Bull. Soc. Math. France}, 88:229--332, 1960.

\bibitem[Ben00]{benoist2000automorphismes}
Yves Benoist.
\newblock Automorphismes des c\^{o}nes convexes.
\newblock {\em Invent. Math.}, 141(1):149--193, 2000.

\bibitem[Ben03]{benoist2003convexes}
Yves Benoist.
\newblock Convexes hyperboliques et fonctions quasisym\'{e}triques.
\newblock {\em Publ. Math. Inst. Hautes \'{E}tudes Sci.}, (97):181--237, 2003.

\bibitem[Ben06a]{benoist2006convexes}
Yves Benoist.
\newblock Convexes divisibles. {IV}. {S}tructure du bord en dimension 3.
\newblock {\em Invent. Math.}, 164(2):249--278, 2006.

\bibitem[Ben06b]{benoist2006convexeshyp}
Yves Benoist.
\newblock Convexes hyperboliques et quasiisom\'{e}tries.
\newblock {\em Geom. Dedicata}, 122:109--134, 2006.

\bibitem[Ben08]{benoist2008survey}
Yves Benoist.
\newblock A survey on divisible convex sets.
\newblock In {\em Geometry, analysis and topology of discrete groups}, volume~6
  of {\em Adv. Lect. Math. (ALM)}, pages 1--18. Int. Press, Somerville, MA,
  2008.

\bibitem[BK553]{BK53}
Chapter iv: Projective metrics.
\newblock volume~3 of {\em Pure and Applied Mathematics}, pages 105--173.
  Elsevier, 1953.

\bibitem[{Bob}20]{bobb2020codimension}
Martin~D. {Bobb}.
\newblock {Codimension-$1$ Simplices in Divisible Convex Domains}.
\newblock {\em arXiv e-prints}, page arXiv:2001.11096, January 2020.

\bibitem[Bow98]{bowditch1998topological}
Brian~H. Bowditch.
\newblock A topological characterisation of hyperbolic groups.
\newblock {\em J. Amer. Math. Soc.}, 11(3):643--667, 1998.

\bibitem[Bow99]{bowditch1999convergence}
B.~H. Bowditch.
\newblock Convergence groups and configuration spaces.
\newblock In {\em Geometric group theory down under ({C}anberra, 1996)}, pages
  23--54. de Gruyter, Berlin, 1999.

\bibitem[Bow12]{bowditch2012relatively}
B.~H. Bowditch.
\newblock Relatively hyperbolic groups.
\newblock {\em Internat. J. Algebra Comput.}, 22(3):1250016, 66, 2012.

\bibitem[BPS19]{bps2019anosov}
Jairo Bochi, Rafael Potrie, and Andr\'{e}s Sambarino.
\newblock Anosov representations and dominated splittings.
\newblock {\em J. Eur. Math. Soc. (JEMS)}, 21(11):3343--3414, 2019.

\bibitem[BV23]{BV}
Pierre-Louis {Blayac} and Gabriele {Viaggi}.
\newblock {Divisible convex sets with properly embedded cones}.
\newblock {\em arXiv e-prints}, page arXiv:2302.07177, February 2023.

\bibitem[{Cho}10]{choi2011}
Suhyoung {Choi}.
\newblock {The convex real projective orbifolds with radial or totally geodesic
  ends: The closedness and openness of deformations}.
\newblock {\em arXiv e-prints}, page arXiv:1011.1060, November 2010.

\bibitem[CLM16]{clm2016convex}
Suhyoung Choi, Gye-Seon Lee, and Ludovic Marquis.
\newblock Convex projective generalized dehn filling.
\newblock {\em Annales scientifiques de l'École normale supérieure}, 53, 11
  2016.

\bibitem[CLT15]{clt2015convex}
D.~Cooper, D.~D. Long, and S.~Tillmann.
\newblock On convex projective manifolds and cusps.
\newblock {\em Adv. Math.}, 277:181--251, 2015.

\bibitem[CLT18]{clt2018deforming}
Daryl Cooper, Darren Long, and Stephan Tillmann.
\newblock Deforming convex projective manifolds.
\newblock {\em Geom. Topol.}, 22(3):1349--1404, 2018.

\bibitem[CM14]{cm2014finitude}
Micka\"{e}l Crampon and Ludovic Marquis.
\newblock Finitude g\'{e}om\'{e}trique en g\'{e}om\'{e}trie de {H}ilbert.
\newblock {\em Ann. Inst. Fourier (Grenoble)}, 64(6):2299--2377, 2014.

\bibitem[DGK17]{dgk2017convex}
Jeffrey {Danciger}, Fran{\c{c}}ois {Gu{\'e}ritaud}, and Fanny {Kassel}.
\newblock {Convex cocompact actions in real projective geometry}.
\newblock {\em arXiv e-prints}, page arXiv:1704.08711, April 2017.

\bibitem[DGK{\etalchar{+}}21]{dgklm2021convex}
Jeffrey {Danciger}, Fran{\c{c}}ois {Gu{\'e}ritaud}, Fanny {Kassel}, Gye-Seon
  {Lee}, and Ludovic {Marquis}.
\newblock {Convex cocompactness for Coxeter groups}.
\newblock {\em arXiv e-prints}, page arXiv:2102.02757, February 2021.

\bibitem[DS05]{ds2005tree}
Cornelia Dru\c{t}u and Mark Sapir.
\newblock Tree-graded spaces and asymptotic cones of groups.
\newblock {\em Topology}, 44(5):959--1058, 2005.
\newblock With an appendix by Denis Osin and Mark Sapir.

\bibitem[Fra91]{frankel89}
Sidney Frankel.
\newblock Applications of affine geometry to geometric function theory in
  several complex variables. {I}. {C}onvergent rescalings and intrinsic
  quasi-isometric structure.
\newblock In {\em Several complex variables and complex geometry, {P}art 2
  ({S}anta {C}ruz, {CA}, 1989)}, volume~52 of {\em Proc. Sympos. Pure Math.},
  pages 183--208. Amer. Math. Soc., Providence, RI, 1991.

\bibitem[Fre97]{freden1997properties}
Eric~M. Freden.
\newblock Properties of convergence groups and spaces.
\newblock {\em Conform. Geom. Dyn.}, 1:13--23 (electronic), 1997.

\bibitem[GM87]{gm1987discrete}
F.~W. Gehring and G.~J. Martin.
\newblock Discrete quasiconformal groups. {I}.
\newblock {\em Proc. London Math. Soc. (3)}, 55(2):331--358, 1987.

\bibitem[Gol88]{goldman88projective}
William~M. Goldman.
\newblock Projective geometry on manifolds, 1988.

\bibitem[GT87]{gt1987pinching}
M.~Gromov and W.~Thurston.
\newblock Pinching constants for hyperbolic manifolds.
\newblock {\em Invent. Math.}, 89(1):1--12, 1987.

\bibitem[GW12]{gw2012anosov}
Olivier Guichard and Anna Wienhard.
\newblock Anosov representations: domains of discontinuity and applications.
\newblock {\em Invent. Math.}, 190(2):357--438, 2012.

\bibitem[{Isl}19]{islam2019rank}
Mitul {Islam}.
\newblock {Rank-One Hilbert Geometries}.
\newblock {\em arXiv e-prints}, page arXiv:1912.13013, December 2019.

\bibitem[IZ19]{iz2019convex}
Mitul {Islam} and Andrew {Zimmer}.
\newblock {Convex co-compact actions of relatively hyperbolic groups}.
\newblock {\em arXiv e-prints}, page arXiv:1910.08885, October 2019.

\bibitem[IZ20]{iz2020convex}
Mitul {Islam} and Andrew {Zimmer}.
\newblock {Convex co-compact representations of 3-manifold groups}.
\newblock {\em arXiv e-prints}, page arXiv:2009.05191, September 2020.

\bibitem[IZ21]{iz2019flat}
Mitul Islam and Andrew Zimmer.
\newblock A flat torus theorem for convex co-compact actions of projective
  linear groups.
\newblock {\em J. Lond. Math. Soc. (2)}, 103(2):470--489, 2021.

\bibitem[IZ22]{iz2022structure}
Mitul {Islam} and Andrew {Zimmer}.
\newblock {The structure of relatively hyperbolic groups in convex real
  projective geometry}.
\newblock {\em arXiv e-prints}, page arXiv:2203.16596, March 2022.

\bibitem[Kap07]{kapovich2007convex}
Michael Kapovich.
\newblock Convex projective structures on {G}romov-{T}hurston manifolds.
\newblock {\em Geom. Topol.}, 11:1777--1830, 2007.

\bibitem[KL18]{kl2018relativizing}
Michael {Kapovich} and Bernhard {Leeb}.
\newblock {Relativizing characterizations of Anosov subgroups, I}.
\newblock {\em arXiv e-prints}, page arXiv:1807.00160, June 2018.

\bibitem[KLP17]{klp2017}
Michael Kapovich, Bernhard Leeb, and Joan Porti.
\newblock Anosov subgroups: dynamical and geometric characterizations.
\newblock {\em Eur. J. Math.}, 3(4):808--898, 2017.

\bibitem[KLP18]{klpdomains}
Michael Kapovich, Bernhard Leeb, and Joan Porti.
\newblock Dynamics on flag manifolds: domains of proper discontinuity and
  cocompactness.
\newblock {\em Geom. Topol.}, 22(1):157--234, 2018.

\bibitem[KS58]{ks1958curvature}
Paul Kelly and Ernst Straus.
\newblock Curvature in {H}ilbert geometries.
\newblock {\em Pacific J. Math.}, 8:119--125, 1958.

\bibitem[Lab06]{labourie2006anosov}
Fran\c{c}ois Labourie.
\newblock Anosov flows, surface groups and curves in projective space.
\newblock {\em Invent. Math.}, 165(1):51--114, 2006.

\bibitem[Mar14]{marquis2013around}
Ludovic Marquis.
\newblock Around groups in {H}ilbert geometry.
\newblock In {\em Handbook of {H}ilbert geometry}, volume~22 of {\em IRMA Lect.
  Math. Theor. Phys.}, pages 207--261. Eur. Math. Soc., Z\"{u}rich, 2014.

\bibitem[Sul79]{sullivan1979density}
Dennis Sullivan.
\newblock The density at infinity of a discrete group of hyperbolic motions.
\newblock {\em Inst. Hautes \'{E}tudes Sci. Publ. Math.}, (50):171--202, 1979.

\bibitem[Sul85]{sullivan1985quasiconformal}
Dennis Sullivan.
\newblock Quasiconformal homeomorphisms and dynamics. {II}. {S}tructural
  stability implies hyperbolicity for {K}leinian groups.
\newblock {\em Acta Math.}, 155(3-4):243--260, 1985.

\bibitem[Tuk98]{tukia1998conical}
Pekka Tukia.
\newblock Conical limit points and uniform convergence groups.
\newblock {\em J. Reine Angew. Math.}, 501:71--98, 1998.

\bibitem[Vin71]{vinberg1971discrete}
\`E.~B. Vinberg.
\newblock Discrete linear groups that are generated by reflections.
\newblock {\em Izv. Akad. Nauk SSSR Ser. Mat.}, 35:1072--1112, 1971.

\bibitem[{Wei}22]{Weisman22}
Theodore {Weisman}.
\newblock {An extended definition of Anosov representation for relatively
  hyperbolic groups}.
\newblock Preprint, \href{http://arxiv.org/abs/2205.07183}{arXiv:2205.07183},
  2022.
\newblock [math.GR].

\bibitem[Yam04]{yaman2004topological}
Asli Yaman.
\newblock A topological characterisation of relatively hyperbolic groups.
\newblock {\em J. Reine Angew. Math.}, 566:41--89, 2004.

\bibitem[Zhu21]{zhu2019relatively}
Feng Zhu.
\newblock Relatively dominated representations.
\newblock {\em Ann. Inst. Fourier (Grenoble)}, 71(5):2169--2235, 2021.

\bibitem[{Zim}20]{zimmer2020higher}
Andrew {Zimmer}.
\newblock {A higher rank rigidity theorem for convex real projective
  manifolds}.
\newblock {\em arXiv e-prints}, page arXiv:2001.05584, January 2020.

\bibitem[Zim21]{zimmer2017projective}
Andrew Zimmer.
\newblock Projective {A}nosov representations, convex cocompact actions, and
  rigidity.
\newblock {\em J. Differential Geom.}, 119(3):513--586, 2021.

\bibitem[ZZ22]{ZZ}
Feng {Zhu} and Andrew {Zimmer}.
\newblock {Relatively Anosov representations via flows I: theory}.
\newblock {\em arXiv e-prints}, page arXiv:2207.14737, July 2022.

\end{thebibliography}

\end{document}